\documentclass[10pt,a4paper]{article}
\usepackage[utf8]{inputenc}
\usepackage[english]{babel}
\usepackage{amsmath}
\usepackage{amsthm}
\usepackage{float}
\usepackage{amsfonts}
\usepackage{amssymb}
\usepackage{graphicx}
\usepackage{enumerate}
\usepackage[margin=1in]{geometry}
\usepackage{appendix}
\usepackage{authblk}
\usepackage{subcaption}
\usepackage{accents}

\usepackage[backend=biber, style=numeric]{biblatex}
\usepackage{placeins}
\usepackage{hyperref}
\usepackage{xcolor}
\hypersetup{
  colorlinks=true,
  linkcolor={blue!80!black},
  citecolor=black,
  urlcolor={blue!80!black}
}

\theoremstyle{plain}
\newtheorem{theorem}{Theorem}[section]
\newtheorem{lemma}[theorem]{Lemma}
\newtheorem{corollary}[theorem]{Corollary}
\newtheorem{proposition}[theorem]{Proposition}
\newtheorem{assumption}{Assumption}

\theoremstyle{definition}
\newtheorem{definition}[theorem]{Definition}

\theoremstyle{remark}
\newtheorem{remark}[theorem]{Remark}

\numberwithin{equation}{section}

\newcommand{\Space}{X}
\newcommand{\E}{{\mathbb E}}
\newcommand{\D}{\mathcal{D}}

\newcommand{\tr}{{\rm tr}}

\newcommand{\R}{{\mathbb R}}
\newcommand{\norm}[1]{\left\|#1\right\|}

\newcommand{\N}{\mathcal{N}}
\newcommand{\Natural}{\mathbb{N}}
\newcommand{\KL}[2]{{\rm KL}\left(#1 \,\|\, #2\right)}

\newcommand{\argmin}{\operatorname*{argmin}}
\newcommand{\dom}{\operatorname*{dom}}
\newcommand{\dist}{\operatorname*{dist}}
\newcommand{\Var}{\text{Var}}
\newcommand{\Cov}{\text{Cov}}
\newcommand{\cl}{\text{cl}}

\newcommand{\Reg}{{\mathcal R}}

\newcommand{\Rsym}{\mathrm{Sym}(m, \R)}

\newcommand{\ubar}[1]{\underaccent{\bar}{#1}}
\newcommand{\F}{\mathcal{F}}

\newcommand{\Sm}{{\mathcal S}(m)}

\newcommand{\wA}{\widehat{A}}
\newcommand{\OAE}{{\mathcal J}}

\title{Optimal approximation error approach in Bayesian inverse problems}
\author[1]{Fabian Schneider\thanks{\href{mailto:fabian.schneider@lut.fi}{fabian.schneider@lut.fi}}}
\author[1]{Tapio Helin\thanks{\href{mailto:tapio.helin@lut.fi}{tapio.helin@lut.fi}}}
\affil[1]{Department of Computational Engineering, LUT University, Lappeenranta, Finland}

\begin{document}

\maketitle

\begin{abstract}
In many inverse problems, evaluating the forward operator is computationally expensive, making it common practice to accelerate inference with a less costly approximate model. However, such approximations introduce unwanted bias and must be carefully accounted for to ensure the reliability of the solution. This phenomenon is particularly pronounced in Bayesian inference, where accurate posterior characterisation typically requires a large number of forward model evaluations.
In this paper, we focus on the \emph{approximation error model} (AEM), which compensates for the modelling error by adjusting the likelihood density. More precisely, AEM incorporates the model error into the additive Gaussian noise distribution. A central challenge in this approach is how to choose the modified noise mean and covariance in a principled way. To address this, we propose a new framework, \emph{optimal approximation error} (OAE) which determines the noise model by minimising the expected information loss, measured via the Kullback-Leibler (KL) divergence between the true and approximate posterior distributions.
We show that the standard \emph{enhanced error approach} (EEA) is recovered as the minimiser of the expected divergence between likelihoods rather than posteriors and develop OAE as a direct alternative targeting posterior accuracy. We establish bi-convexity of the OAE objective, introduce a regularised, provably convergent proximal alternating minimisation scheme, and derive closed-form expressions in the Gaussian linear setting. We further show that OAE estimates computed from empirical training data are statistically consistent. Numerical experiments on a linear CT problem and a nonlinear Darcy-flow problem show that OAE improves posterior accuracy relative to both the uncorrected and EEA-corrected likelihoods.
\end{abstract}

\section{Introduction}

We consider an inverse problem
\begin{equation}\label{eq: inverse problem}
    y = A(x) + \delta \epsilon,
\end{equation}
where $y \in \R^m$ is the observed data, $x\in \Space$ is the unknown, taking values in a separable Banach space $\Space$ and $A: \Space \rightarrow \R^m$ is the, possibly nonlinear, forward map. Moreover, the observational noise vector $\epsilon \in \R^m \sim \N(0, I_m)$ has white noise statistics and $\delta>0$ represents the noise level.

We adopt the Bayesian approach \cite{kaipio2006statistical, Stuart2010Inverse} and assume to have some prior knowledge of the distribution of the unknown available, modelled by the probability measure $\mu_0$. The posterior distribution $\mu(\cdot|y)$, i.e., the conditional distribution of $x$ given the observation $y$, then blends the prior distribution with the likelihood induced by \eqref{eq: inverse problem} via Bayes' formula, expressed through the Radon-Nikodym derivative
\begin{eqnarray*}
         \frac{d\mu}{d\mu_0}(x | y) = \frac{1}{Z(y)} \exp \left( -\frac{1}{2 \delta^2} \norm{A(x) -y}_2^2\right) \quad \text{in}\; L^1(\mu_0)
\end{eqnarray*}
with normalising constant
\begin{eqnarray*}
    Z(y) = \int_{\Space}  \exp \left( -\frac{1}{2 \delta^2} \norm{A(x) -y}_2^2\right) \mu_0(dx).
\end{eqnarray*}
For more details, see e.g. \cite{Stuart2010Inverse}. 
In what follows, we write $p(y|x) = (2\pi\delta^2)^{-m/2} \exp\left(-\tfrac{1}{2\delta^2}\norm{A(x)-y}^2\right)$ for the likelihood density and $p(y) = (2\pi\delta^2)^{-m/2} Z(y)$ for the marginal density of $y$.

In many practical settings, the forward operator is either not precisely known or computationally expensive to evaluate. Therefore, it is often convenient to study the problem using a less computationally intensive approximate model, say, by replacing the partial differential equation (PDE) solver associated with $A$ by one that employs a coarser domain discretisation. However, such an approximation induces an unwanted error to the posterior distribution. 
Consequently, the trade-off between induced posterior error and computational gains has been intensively studied during the last two decades or so and a number of surrogate modelling approaches have been proposed; among them, we highlight \cite{Kennedy2001Bayesian, Marzouk2007Stochastic, Brynjarsdottir2014Learning, Bui-Thanh2008ModelReduction, Approximation_error_2006}. In the inverse problems literature, such questions are also tackled in the context of stability analysis related to the Bayesian posterior \cite{stuart2016posterior}.

In this paper, we focus on an early likelihood-based surrogate model for inverse problems, introduced in \cite{Approximation_error_2006} and known as the approximation error model (AEM).
The central idea is that, given a computationally efficient forward approximation $\wA$, one compensates for the resulting posterior error by recalibrating the mean and covariance of the likelihood distribution. Given an approximation $\wA$, the AEM approach approximates this with
\begin{equation}
    \label{eq: surrogate likelihood}
    \hat p(y | x; \xi, \Gamma) = {\mathcal N}(\wA(x) + \xi, \Gamma),
\end{equation}
where $\xi \in \R^m$ and $\Gamma \in \Sm \subset \R^{m\times m}$, with $\Sm$ standing for the symmetric positive definite matrices, are free parameters. We denote by $ \hat p(y; \xi, \Gamma)$ the marginal distribution of $y$ in \eqref{eq: surrogate likelihood}. The approximate likelihood and the prior give rise to an approximative posterior 
\begin{equation}
	\label{eq: surrogate posterior}
	\frac{d\hat\mu}{d\mu_0}(x | y; \xi, \Gamma) = \frac{\hat p(y | x; \xi, \Gamma)}{\hat p(y; \xi, \Gamma)} \quad \text{in}\; L^1(\mu_0).
\end{equation}
Further, parameters $(\xi,\Gamma)$ are then calibrated to reduce the resulting posterior error with $\hat \mu(\cdot | y)$ approximating $\mu(\cdot | y)$.

In our view, the AEM approach has sparked significant attention in the application domain (see e.g. \cite{Approximation_error_2006, Kolemainen2011Marginalization, Nicholson_2018Estimation, Nissinen_2009Compensation, lunz2020learnedoperatorcorrectioninverse, Sahlstrom2020Modeling,Heinio2005Compensation,Yamashita2016Multi-subject}) due to three primary reasons: first, the easy-to-use formulation combining a Gaussian likelihood with the given arbitrary approximative scheme for the forward mapping. Second, intuitive interpretability of parameters $\xi$ and $\Gamma$ being descriptive of the probability distribution of the sum of observational noise and model error (see discussion in \cite{kaipio2006statistical}). Third, the empirical evidence of improved posterior modelling reported in many works \cite{ Approximation_error_2006, kaipio2006statistical, Lahivaara2015Estimation, Rimpilinen_2017_Bayesian, Banasiak2012Improving}.

The immediate fundamental question in AEM is how to effectively choose the parameters $\xi$ and $\Gamma$. The typical strategy in the literature is the \emph{enhanced error approach} (EEA) \cite{kaipio2006statistical}. 
In EEA, the parameters $(\xi, \Gamma)$ are chosen according to the moments of the residual between the measurement and the surrogate averaged over the true Bayesian joint distribution:
\begin{eqnarray}\label{eq: enhanced error approach}
    \xi_* = \E^\lambda\left(y - \wA(x) \right), \qquad \Gamma_* = \E^\lambda\left[\left(y - \wA(x)-\xi_*\right) \otimes \left(y - \wA(x)-\xi_*\right)\right],
\end{eqnarray}
where $\lambda$ is the joint law of $(x,y)$ on $X \times \R^m$ under \eqref{eq: inverse problem} and the prior $\mu_0$.
In particular, this formulation of the EEA enables straightforward computational estimation of $\xi_*$ and $\Gamma_*$ using supervised training data generated from the true model \eqref{eq: inverse problem}.

In this paper, we adopt a novel perspective on AEM by exploring the \emph{optimal} choice of parametrisation $(\xi,\Gamma)$ in \eqref{eq: surrogate posterior}. We show that, as expected, the EEA can be interpreted as fitting the best possible likelihood approximation when averaged over the prior. This leaves open the key question of our paper: how to directly optimise the \emph{posterior error} given prior information, whether in the form of a closed-form distribution or empirical data?

\subsection{Our contribution}

The main contribution of this paper is to propose and analyse an optimality condition for the surrogate posterior by minimising the \emph{expected information loss}, i.e., we solve
\begin{equation}
    \label{eq:AEM_optimization}
    (\hat{\xi},\widehat{\Gamma}) = \argmin_{\xi, \Gamma}\, \E^p\, \KL{\mu(\cdot | y)}{\hat \mu(\cdot | y; \xi, \Gamma)},
\end{equation}
where the \textit{Kullback-Leibler} (KL) divergence is defined in \eqref{eq: KL div} ahead. We refer to the resulting criterion as the \emph{optimal approximation error} (OAE) approach, and to any corresponding minimiser $(\hat{\xi}, \widehat{\Gamma})$ as the OAE parameters.
Notice that similar Bayes cost is frequently utilised in the context of Bayesian optimal experimental design \cite{Huan_2024}. 
In what follows, we write $G := \Gamma^{-1}$ for the precision matrix, which we use as the optimisation variable.
In particular, in the spirit of variational Bayes, the optimisation task in \eqref{eq:AEM_optimization} reduces to
\begin{equation*}
    (\hat{\xi},\widehat{G}) = \argmin_{\xi, G}\, \OAE(\xi,G), \quad \text{where} \quad \OAE(\xi, G) = -\E^\lambda \log \left(\frac{d \hat \mu(\cdot | y; \xi, G)}{d \mu_0}(x)\right),
\end{equation*}
see Lemma \ref{lem: decompositon of obj fct}.

Our key findings are as follows:
\begin{itemize}
    \item We prove in Theorem~\ref{thm: enhanced error approach} that the EEA parametrisation \eqref{eq: enhanced error approach} is the unique global minimiser of the prior-averaged KL divergence between the true and the approximate likelihood, i.e.,
    	\begin{equation*}
    		(\xi_*, G_*) = \argmin_{\xi,G} \E^{\mu_0}\, \KL{p(\cdot|x)}{\hat p(\cdot|x; \xi,G)}.
    	\end{equation*}
	This confirms that EEA is optimal in a likelihood-averaged sense, but leaves open whether directly optimising the posterior error in \eqref{eq:AEM_optimization} yields further benefits.
    \item In Section~\ref{sec: OAE} we establish the fundamental properties of the functional $\OAE$ needed for practical optimisation. Under Assumption~\ref{ass: assumption} we show in Theorem~\ref{thm: Expected KL divergence analytic} that it is twice continuously differentiable. In Theorem~\ref{thm: biconvexity} we show that $\OAE$ is bi-convex in $(\xi, G)$ and characterise its critical points. 
     \item In the special case of a linear surrogate $\wA$ and a Gaussian prior $\mu_0$ (Section~\ref{sec: Gaussian linear}), we derive $\OAE$ in closed form and show in Theorem~\ref{lem: Gaussian case; xi} that a minimiser for $\xi$ is given by $\hat \xi = \E^\lambda(y - \wA(x))$, coinciding with the EEA choice \eqref{eq: enhanced error approach}.
    \item For practical optimisation, in Section~\ref{sec:numerical_aspects} we first exhibit that the infimum is in general not attained and further the level sets of $\OAE$ may be unbounded. We hence introduce the regularised functional $\OAE_\alpha := \OAE + \alpha \Reg$, $\alpha>0$, where $\Reg(G) = \delta^2\tr(G) - \log\det G$ is convex and coercive and penalises deviations of $\Gamma = G^{-1}$ from the observational noise covariance (Lemma~\ref{lem:regularisation}). Under a mild non-degeneracy condition (Assumption~\ref{ass: assumption conv of AM}) we show $\OAE_\alpha$ is Morse and coercive and that all critical points are local minima, so that the proximal alternating minimisation scheme of Definition~\ref{def:proximal_alternating} converges to a critical point of $\OAE_\alpha$ (Proposition~\ref{prop:iterates_bounded}), building on the general convergence theory of \cite{attouch2013proximal}.
    \item In Section~\ref{subsec:empirical} we devise two plug-in estimators from empirical data, a single-loop estimator for the Gaussian linear case and a double-loop estimator for the general case, and prove in Theorems~\ref{thm: M estimation} and~\ref{cor: consistency for double loop} that they concentrate on the set of minimisers of the regularised objective $\OAE_\alpha$ as the sample size $M \to \infty$.
    \item Finally, we validate the proposed approach on a linear CT-based and a nonlinear PDE-based (Darcy flow) inverse problem in Section~\ref{sec: numerics}, where OAE is shown to reduce the expected posterior KL divergence and improve mean-squared error relative to both the unmitigated and EEA-corrected likelihoods.
\end{itemize}

Bayesian inverse problems often seek to infer quantities such as coefficients of PDEs or integral equations that are inherently infinite dimensional. While any numerical solution requires a discretisation of the problem, there has been a significant effort to develop methods that are independent of the choice of discretisation \cite{Lassas2004CANONE, Lassas2009Discretization, Lehtinen1989LinearIP}. 
In line with \cite{Stuart2010Inverse}, we seek discretisation-invariant methods 
and avoid discretisation until the last possible moment. 
The core results of Section~\ref{sec: OAE} hold for unknowns taking values in a separable, potentially infinite-dimensional, Banach space $\Space$, while the observed data $y \in \R^m$ and the parameters $(\xi, \Gamma)$ remain finite-dimensional. Section~\ref{sec: Gaussian linear} specialises to Hilbert spaces and Section~\ref{subsec:empirical} to $\Space = \R^n$.

\subsection{Literature overview}
Variational inference is a classical technique that approximates a posterior by a member of a tractable family of distributions. It rests on writing the log-evidence as the sum of the so-called \emph{evidence lower bound} and the KL divergence from the approximating distribution to the posterior. Since the log-evidence does not depend on the approximating distribution, maximising the former is equivalent to minimising the latter. See e.g. \cite[Section 9.4]{bishop_pattern_2006}. Our objective \eqref{eq:AEM_optimization} has the same structure (Lemma~\ref{lem: decompositon of obj fct}), but with two differences: the approximating family is the set of AEM posteriors \eqref{eq: surrogate posterior}, parametrised by $(\xi,\Gamma)$, and the divergence is taken in the opposite, mean-seeking direction $\KL{\mu}{\hat\mu}$ and averaged over the data.
Optimal approximation of measures with respect to KL divergence is studied in the literature, we mention in particular the works of \cite{PinskiKLApproximation2015, YulongG2017GaussianApproximation} where mode-seeking approximation of general measures in a class of  Gaussian or Gaussian mixtures are explored.

The literature on surrogate modelling in inverse problems is very vast. The importance of accounting for a misspecified forward model is well-studied, e.g. in the works of \cite{Brynjarsdottir2014Learning, Kennedy2001Bayesian, kaipio2006statistical, Smigth2013UQ}.
Classical strategies include the complete error model \cite{KAIPIO2007Statistical}, which infers the unknown and the model error jointly, at the cost of a significantly higher-dimensional inference problem.
Another idea presented in the literature is based on a two-stage sampling algorithm. Therein samples generated utilising the surrogate are used to speed up an MCMC algorithm employing the exact forward map \cite{ZENG2012Asparse, Laloy2013Efficient}.

We mention a few of the many further strategies discussed in the literature.
An adaptive scheme that constructs local approximations of the posterior distribution depending on measurement-data is proposed in \cite{Jinglai2014Adaptive}. The authors of \cite{Calvetti_2018} propose an iterative scheme updating the posterior densities through the model error.
A framework that enables embedded internal model adjustment is proposed in \cite{Sargsyan2015Onthestatistical, Sargsyan_2019Embeddedmodelerror}. 
Therein an additive modification of model parameters using a Polynomial Chaos representation is utilised.
In \cite{lunz2020learnedoperatorcorrectioninverse}, a learned model correction utilising measurement data is discussed based on a variational approach. 
Perhaps the closest contribution to ours is given in the recent work of \cite{cvetkovic2024choosing}, where the introduced distance of an approximate posterior to the posterior distribution in a KL divergence based metric is studied. The authors analyse multiple posterior approximations including EEA and the complete error model. Further the choice of the observation operator is analysed from the viewpoint of mitigating model error.

Convex and non-convex minimisation problems based on the Kurdyka-\L{}ojasiewicz property have been studied in \cite{attouch2013proximal, Xu2013ABlockCoordinate}, including convergence guarantees for numerical algorithms. These results are based on the Kurdyka-\L{}ojasiewicz inequalities \cite{Kurdyka1998,Ljasiewicz2013UneProperiete} and have since been applied in numerous works (see e.g. \cite{CORNELIO2015Aconvergent, Tao2017Convergence, Chen2019Nonnegtive, Li2019Alternating}).

\emph{Structure of the paper.} Section~\ref{sec: OAE} introduces the OAE functional formally, relates it to the EEA, and establishes smoothness (Theorem~\ref{thm: Expected KL divergence analytic}) and bi-convexity (Theorem~\ref{thm: biconvexity}). 
Section~\ref{sec: Gaussian linear} specialises to linear surrogates $\wA$ with a Gaussian prior, where $\OAE$ and its minimiser over $\xi$ admit closed form (Theorem~\ref{lem: Gaussian case; xi}).
Section~\ref{sec:numerical_aspects} addresses the practical optimisation of $\OAE$: we exhibit non-attainment of the infimum and unbounded level sets. We introduce a regularised functional $\OAE_\alpha$, and prove convergence of a proximal alternating minimisation scheme to a critical point (Proposition~\ref{prop:iterates_bounded}). Section~\ref{subsec:empirical} studies stability of the method: we establish consistency of an empirical estimator as the number of training samples grows (Theorems~\ref{thm: M estimation} and \ref{cor: consistency for double loop}). Section~\ref{sec: numerics} validates the approach numerically on a linear CT-based inverse problem and a nonlinear PDE-based inverse problem governed by Darcy's law. Section~\ref{sec: conclusions} draws concluding remarks and gives an outlook on future work.

\section{Optimal approximation error}
\label{sec: OAE}

For two probability measures $\mu$ and $\nu$ on $\Space$ such that $\mu \ll \nu$, the KL divergence is defined by 
\begin{equation}\label{eq: KL div}
	\KL{\mu}{\nu} = \E^\mu \log\left(\frac{d\mu}{d\nu}(x)\right).
\end{equation}
Notice that KL divergence is not symmetric;
in the following we use the arrangement above with the aim to optimise $\nu$ to match $\mu$. This approach is also called \emph{mean-seeking}. Depending on the arrangement, different behaviour of the approximation is enforced. If $\mu$ is a fixed  Gaussian mixture and $\nu$ is optimised over the set of all Gaussian distributions, optimal $\nu$ has the same mean and covariance as $\mu$. If the reverse ordering is chosen, the objective tends to have local minima matching the modes of $\mu$ (see \cite{bishop_pattern_2006, YulongG2017GaussianApproximation}).

For convenience, in what follows, we reparametrise the likelihood approximation in \eqref{eq: surrogate likelihood} 
with respect to the mean-shift $\xi$ and the precision matrix denoted by $G:= \Gamma^{-1}$. Moreover, below we utilise notation
$\norm{y}_G^2 = \langle y, G y\rangle$.

Let us next study the EEA method from the perspective of a minimisation problem. Let us define a functional ${\mathcal K} :\R^m \times \Sm \rightarrow \R$ by setting
\begin{equation}
	{\mathcal K}(\xi, G) = \E^{\mu_0}\, \KL{p(\cdot|x)}{\hat p(\cdot|x; \xi, G)}
\end{equation}

In order to prove that the objective functional $\mathcal{K}$ is  well-defined, we need to impose further assumptions on the problem.
\begin{assumption}\label{ass: assumption}
The operators $\mathcal{A} \in \{A, \wA\}$ together with the prior distribution $\mu_0$ satisfy following conditions:
    \begin{itemize}
        \item[(A1)] (Lipschitz) There exists $L_1 > 0$ such that
        \begin{eqnarray*}
            \norm{\mathcal{A}(x) - \mathcal{A}(x')} \le L_1 \norm{x - x'}
        \end{eqnarray*}
        \item[(A2)] (sub-Gaussian prior) There exists $L_2 > 0$ such that
        \begin{eqnarray*}
            \E^{\mu_0} \exp\left( L_2 \norm{x}^2\right) < \infty
        \end{eqnarray*}
    \end{itemize}
\end{assumption}

The assumption ensures the boundedness of several key expectations and moment terms, see Theorem~\ref{thm: enhanced error approach} and Lemma~\ref{lem: assumptions} ahead. We point out now, that Assumptions~\ref{ass: assumption} are inherited from \cite{helin2025bayesianoptimalexperimentaldesign}. Indeed, condition (iii) of \cite[Assumption 3.1]{helin2025bayesianoptimalexperimentaldesign} follows from (A1). Hence the auxiliary estimates of \cite[Section 3.1]{helin2025bayesianoptimalexperimentaldesign}, which we use throughout, apply. We emphasise that these estimates, unlike the stability theorems of that work, do not require any coupling between $L_1$ and $L_2$.

\begin{remark}
        We point out that \eqref{eq: enhanced error approach} implies $\Gamma_* \succeq \delta^2 I$ by the law of total variance. For $e = y - \widehat{A}(x) - \xi_*$ we have
    \begin{eqnarray*}
        \Gamma_*& =& \Cov(e) = \Cov\left(\E (e|x)\right) + \E \left(\Cov (e | x) \right) \\
        & = &\Cov\left(\E (e|x)\right) + \delta^2 I_m \succeq \delta^2 I.
    \end{eqnarray*}
    This is a natural condition for $\Gamma_*$, by which the introduction of a surrogate cannot decrease the optimal noise covariance.
\end{remark}

\begin{theorem}\label{thm: enhanced error approach}
Let $\xi_* \in \R^m$ and $\Gamma_* \in \Sm$ be given by \eqref{eq: enhanced error approach}, and let $G_* = \Gamma_*^{-1}$. If Assumption \ref{ass: assumption} holds, then $(\xi_*, G_*)$ is the unique global minimiser of ${\mathcal K}$. 
\end{theorem}

\begin{proof}[Proof of Theorem \ref{thm: enhanced error approach}:]
\begin{eqnarray*}
    {\mathcal K}(\xi, G) &=& \int_\Space\int_{\R^m} \log \frac{p(y|x)}{\hat p(y|x; \xi, G)} p(y|x) dy \mu_0(dx) \\
    & = & \mathrm{const}
     - \E^\lambda \log \hat p(y|x; \xi, G).
\end{eqnarray*}
where $\mathrm{const} := \int_\Space\int_{\R^m} \log p(y | x) p(y | x) dy \mu_0(dx)$ is constant with respect to $(\xi, G)$. Notice that $\mathcal{K}$ is finite under Assumption \ref{ass: assumption}.
 Hence, using the definition of the Gaussian likelihood, it holds
\begin{eqnarray*}
     {\mathcal K}(\xi, G)  
     &=&  \mathrm{const} + \frac{m}{2}\log 2 \pi -  \frac12 \log \det G +  \tfrac{1}{2} \, \E^\lambda \norm{\wA(x) -y + \xi}_G^2.
\end{eqnarray*}
We optimise over $\xi$ and $G$ separately. Since $G$ is strictly positive definite, the term 
\begin{eqnarray*}
    \tfrac{1}{2} \, \E^\lambda \norm{\wA(x) -y + \xi}_G^2    
\end{eqnarray*}
is minimised in $\xi$ if and only if $\mathbb{E}^\lambda[y - \widehat{A}(x) - \xi] = 0$, resulting in $\xi_*$ given by \eqref{eq: enhanced error approach}.
We now optimise over $G$. Set $e = y - \widehat{A}(x) - \xi_*$. Then
    \begin{align*}
        \mathbb{E}^\lambda\|e\|_G^2
        &= \mathbb{E}^\lambda \tr(G\, e e^\top)
        = \tr(G\,\Cov^\lambda(e)).
    \end{align*}
    To differentiate $\mathcal{K}$, we notice that $G \mapsto -\log \det G$ is twice continuously differentiable on $\Sm$. Denote by $D_G[\Delta_G]$ the Fréchet derivative evaluated in a direction $\Delta_G \in \Rsym$. It holds by \cite{Matrixcookbook} that
    \begin{eqnarray*}
        D_G[\Delta_G] \log \det G = \tr (G^{-1} \Delta_G).
    \end{eqnarray*}
    Plugging these observations into the definition of $\mathcal{K}$,
    \begin{eqnarray*}
        D_G[\Delta_G] \mathcal{K}(\xi_*, G) = \tfrac12\tr \left([\Cov^\lambda(e) -G^{-1}] \Delta_G\right).
    \end{eqnarray*}
    This derivative vanishes for $G_*$ in \eqref{eq: enhanced error approach}, so it remains to show that the second derivative $D^2_G[\Delta_G, \Delta_G]$ evaluated at $\Delta_G$ is positive. Using \cite[Section 2.2]{Matrixcookbook} we get for $\Delta_G \neq 0$
    \begin{eqnarray*}
        D_G^2[\Delta_G, \Delta_G] \mathcal{K}(\xi_*, G) = \tfrac12 \tr(G^{-1} \Delta_G G^{-1} \Delta_G) >0,
    \end{eqnarray*}
    for every $G \in \Sm$ and hence $\mathcal{K}$ is strictly convex in $G$ when $\xi$ is fixed. Since $\xi_*$ is independent of $G$, $(\xi_*, G_*)$ is the unique global minimiser.
\end{proof}

What Theorem \ref{thm: enhanced error approach} indicates is that the EEA approach indeed provides best possible likelihood approximation when averaged over the prior. However, as we will later see, this can differ from the best possible approximation in terms of posterior distributions. 
We now define the objective functional $\OAE$ on $\R^m \times \Sm $ through
\begin{equation}
	\label{eq:OAE}
    \OAE: \R^m \times \Sm \to \R, \quad \OAE(\xi, G) = -\E^\lambda \log \left(\frac{d \hat \mu(\cdot | y; \xi, G)}{d \mu_0}(x)\right).
\end{equation}
The functional $\OAE$ satisfies the following identity.

\begin{lemma}\label{lem: decompositon of obj fct}
    Let Assumption \ref{ass: assumption} hold. It holds for any $(\xi, G) \in \R^m \times \Sm$ that
    \begin{eqnarray}\label{eq: from KL divergence to our objective}
        \E^p\, \KL{\mu(\cdot | y)}{\hat \mu(\cdot | y; \xi, G)} = \E^p\, \KL{\mu(\cdot|y)}{\mu_0} + \OAE(\xi, G).
    \end{eqnarray}
\end{lemma}
    \begin{proof}[Proof of Lemma \ref{lem: decompositon of obj fct}:]
By definition of KL divergence:
    \begin{eqnarray}\label{eq: KL div applied to our poblem}
    \E^p \KL{\mu(\cdot | y)}{\hat \mu(\cdot | y;\xi, G)} 
     = \E^\lambda \log \left(\frac{d\mu(\cdot| y)}{d \hat \mu(\cdot| y;\xi, G)}(x)  \right).
\end{eqnarray}
Note that all measures $\mu_0, \mu(\cdot|y), \hat \mu(\cdot| y; \xi, G)$ on $\Space$ are equivalent for any $y \in \R^m, \xi \in \R^m$ and $G \in \Sm$. We can write
     \begin{eqnarray}\label{eq: composition of kl term}
       \log \left(\frac{d\mu(\cdot| y)}{d \hat \mu(\cdot| y;\xi, G)}(x) \right)= \log \left(\frac{d \mu(\cdot|y)}{d \mu_0}(x)  \right) -\log \left(\frac{d \hat \mu(\cdot|y; \xi, G)}{d \mu_0}(x)\right).
    \end{eqnarray}
    Plugging into \eqref{eq: KL div applied to our poblem} yields the claim. By Lemma \ref{lem: assumptions} ahead, both quantities on the right-hand side are finite under Assumption \ref{ass: assumption}.
\end{proof}

Since the term $\E^p\, \KL{\mu(\cdot|y)}{\mu_0}$ is independent of $\xi$ and $G$, we have that the optimal parameters with respect to the posterior error on the left-hand side of \eqref{eq: from KL divergence to our objective} also minimise $\OAE$, and we can focus our study to $\OAE$ in what follows. We show now that $\OAE$ is finite under Assumption \ref{ass: assumption}.

\begin{lemma}[$\OAE$ is well-defined]\label{lem: assumptions}
Suppose Assumption \ref{ass: assumption} holds. Then it holds for every $(\xi, G) \in \R^m \times \Sm$ that
\begin{equation*}
   -\infty < - \E^p\, \KL{\mu(\cdot|y)}{\mu_0}  \le   \OAE(\xi, G)  < \infty.
\end{equation*}
\end{lemma}
\begin{proof}
    We begin with the lower bound. By \eqref{eq: composition of kl term} we have
    \begin{align*}
        \OAE(\xi, G) \ge - \E^p\, \KL{\mu(\cdot|y)}{\mu_0}.
    \end{align*}
    By straightforward calculations
    \begin{eqnarray*}
        \E^p\, \KL{\mu(\cdot|y)}{\mu_0} &=& \iint_{\Space \times \R^m } \log \frac{p(y|x)}{p(y)} \lambda(dx, dy) \\
        & = & -\frac m2 \log (2\pi \delta^2) -\frac{1}{2\delta^2}\E^\lambda \norm{A(x) -y}^2 - \int_{\R^m} p(y) \log p(y) dy.
    \end{eqnarray*}
    The second term equals $-\tfrac{m}{2}$ by \eqref{eq: inverse problem} and the final term of right-hand side is finite by \cite[Lemma 3.11]{helin2025bayesianoptimalexperimentaldesign}. It is used that $p(y) = c_m Z(y)$.
    For the upper bound, we make the observation 
    \begin{eqnarray}\label{eq: objective function; easy}
        \OAE(\xi, G) = \tfrac{1}{2}\E^\lambda \norm{\wA(x) -y + \xi}_G^2 + \E^p \log \hat Z(y; \xi, G),
    \end{eqnarray}
    where 
    \begin{eqnarray*}
        \hat Z(y; \xi, G) = \E^{\mu_0} \exp\left(-\frac{1}{2} \norm{\wA(x') -y + \xi}_G^2\right).
    \end{eqnarray*}
    The first term in \eqref{eq: objective function; easy} is bounded by an application of triangle inequality in conjunction with \cite[Corollary 3.8]{helin2025bayesianoptimalexperimentaldesign}, making the observation that $\lambda$ has finite expectation for  polynomials in $(x, y)$.  
    The latter term can be bounded through \cite[Lemma 3.11]{helin2025bayesianoptimalexperimentaldesign}.
    
\end{proof}

The following theorem establishes that the objective function is smooth.

\begin{theorem}[$\OAE$ is smooth]\label{thm: Expected KL divergence analytic}
Suppose Assumption \ref{ass: assumption} holds. The objective functional $\OAE$ is twice continuously differentiable and, in particular, satisfies
\begin{equation*}
    (D^{\gamma}_{\eta} \OAE )(\xi, G) = - \E^\lambda\left[ D^{\gamma}_\eta \log \left(\frac{d \hat \mu(\cdot | y; \xi, G)}{d \mu_0}(x)\right) \right],
\end{equation*}
where $\gamma = 1, 2$ and $\eta \in \{\xi, G\}^\gamma$.
\end{theorem}

The following lemma will be helpful to establish dominating functions to justify differentiation under the integral (see \cite[Theorem 2.27]{folland2013real}) in our setting.
    \begin{lemma}\label{lem: poly moments}
        Let Assumption \ref{ass: assumption} hold. Then it holds that for any compact $U \subset \R^m \times \Sm$ and $k \in \Natural$ that there exists $ C(U, k) < \infty$ such that 
        \begin{eqnarray*}
            \sup_{(\xi, G) \in U} \E^p \E^{\hat\mu(\cdot| y; \xi, G)} \left( \norm{x}^k + \norm{y}^k \right) \le C(U, k).
        \end{eqnarray*}
    \end{lemma}
    \begin{proof}
        We show that this follows from the arguments of \cite{schneider2026convexapproximationframeworkneural}. Notice first by \cite[Lemma 3.10]{helin2025bayesianoptimalexperimentaldesign} applied to $(y-\xi)$ that for $(\xi, G) \in U$ and $\varepsilon > 0$, there exists $C(U)$ and $\kappa(\varepsilon, U)>0$ with 
        \begin{eqnarray*}
        \norm{\wA(x) -y + \xi}_G^2 \le C(U) \, (1 + \norm{x}^2_2 + \norm{y-\xi}^2_2) \le 2C(U) \, (1 + \norm{\xi}^2 + \norm{x}^2_2 + \norm{y}^2_2).
        \end{eqnarray*}
        By \cite[Lemma 3.10]{helin2025bayesianoptimalexperimentaldesign} and Young's inequality
        \begin{eqnarray*}
             -C(U) - \kappa(\varepsilon, U) \norm{\xi}^2 - \varepsilon \norm{x}^2_2 + \tfrac{\kappa(\varepsilon, U)}{2} \norm{y}^2 \le -C(U) - \varepsilon \norm{x}^2_2 + \kappa(\varepsilon, U) \norm{y-\xi}^2 \le \norm{\wA(x) -y + \xi}_G^2.
        \end{eqnarray*}
        Hence the Assumptions of \cite[Lemma A.3]{schneider2026convexapproximationframeworkneural} are satisfied and we may use that there exists $q>1$ such that
        \begin{eqnarray*}
            \sup_{(\xi, G) \in U} \E^{\lambda} \left(\frac{d \hat \mu(\cdot| y;\xi, G)}{d  \mu(\cdot| y)}\right)^q < \infty.
        \end{eqnarray*}
        We can now prove the claim by using Hölder's inequality for $q$ given above and $r = \tfrac{q}{q-1}$:
        \begin{eqnarray*}
             \sup_{(\xi, G) \in U} \E^p \E^{\hat\mu(\cdot| y; \xi, G)} \left(\norm{x}^k+\norm{y}^k\right) &=&  \sup_{(\xi, G) \in U} \E^\lambda \left( (\norm{x}^k+\norm{y}^k)   \frac{d \hat \mu(\cdot| y;\xi, G)}{d  \mu(\cdot| y)}\right) \\
             &\le& \left( \E^\lambda \left(\norm{x}^{k}+ \norm{y}^{k} \right)^r\right)^{1/r} \sup_{(\xi, G) \in U} \left( \E^{\lambda} \left(\frac{d \hat \mu(\cdot| y;\xi, G)}{d  \mu(\cdot| y)}\right)^q\right)^{1/q}.
        \end{eqnarray*}
    \end{proof}
    \begin{proof}[Proof of Theorem \ref{thm: Expected KL divergence analytic}]
    We compute the Fréchet derivatives $D^{\gamma}_\eta[\Delta_\eta]$ evaluated in a direction $\Delta_\eta$ of 
    \begin{equation*}
        f_{x,y} (\xi, G) := -\log \left(\frac{d \hat \mu(\cdot | y; \xi, G)}{d \mu_0}(x)\right) = \tfrac{1}{2} \norm{\wA(x) -y + \xi}_G^2 +  \log \E^{\mu_0} \exp\left(-\tfrac{1}{2}\norm{\wA(x') -y + \xi}_G^2 \right).
    \end{equation*}
    Introduce the notation 
    \begin{eqnarray*}
        g_{x,y}(\xi, G; \Delta_\xi) &=& (\wA(x) -y + \xi )^\top G \Delta_\xi \\
        h_{x,y}(\xi, G; \Delta_G) &=& \frac{1}{2} \norm{\wA(x) - y + \xi}_{\Delta_{G}}^2 \\
        k_{x,y}(\xi; \Delta_\xi, \Delta_G) &=& (\wA(x) -y + \xi)^\top \Delta_{G} \Delta_\xi
    \end{eqnarray*}
    and notice that for $C$ constant
    \begin{align}
    \begin{aligned}\label{eq: polynomial moments}
    & |g_{x,y}(\xi, G; \Delta_\xi)| + |h_{x,y}(\xi, G; \Delta_G)| + |k_{x,y}(\xi; \Delta_\xi, \Delta_G)| + \norm{\wA(x)-y+\xi}_G^2 \\
    & \quad \le C (1 + \norm{\Delta_G}+ \norm{\Delta_\xi} + \norm{\Delta_G}\norm{\Delta_\xi}) \left( (\norm{x} + 1)^2 + (\norm{y} + 1)^2 + (\norm{\xi} + 1)^2\right),
    \end{aligned}
    \end{align}
    where \cite[Corollary 3.8]{helin2025bayesianoptimalexperimentaldesign} was used to bound $|\wA(x)|$ in terms of $\norm{x}+1$.
     We notice first that 
    \begin{eqnarray*}
        D_\xi[\Delta_\xi] \, \hat Z(y; \xi, G) & = & -\E^{\mu_0}\left[ g_{x',y}(\xi, G; \Delta_\xi) \exp(-\frac{1}{2} \norm{\wA(x') - y + \xi}_G^2) \right] \\
        D_G[\Delta_G] \, \hat Z(y; \xi, G) & = & -\E^{\mu_0} \left[ h_{x',y}(\xi, G; \Delta_G) \exp(-\frac{1}{2} \norm{\wA(x') - y + \xi}_G^2) \right],
    \end{eqnarray*}
    Notice that the exponential term is bounded by one and $g, h$ can be bounded by polynomials of order $2$ by \eqref{eq: polynomial moments}. We can find a dominating function (locally in $(\xi, G)$) which is $\mu_0-$ integrable by Assumption \ref{ass: assumption}.

    The following derivatives may be taken under the integral with respect to $\lambda$ by using \eqref{eq: polynomial moments} together with Lemma \ref{lem: poly moments} and since $\lambda$ is sub-Gaussian:
    \begin{align}\label{eq: xi, G derivative}
    \begin{aligned}        
        D_\xi[\Delta_\xi] f_{x,y} (\xi, G) 
        & = g_{x,y}(\xi, G; \Delta_\xi) - \E^{\hat\mu(\cdot| y; \xi, G)} g_{x',y}(\xi, G; \Delta_\xi)  \\
        D_G[\Delta_G] f_{x,y}(\xi, G) 
         & = h_{x,y}(\xi, G; \Delta_G) - \E^{\hat\mu(\cdot| y; \xi, G)} h_{x',y}(\xi, G; \Delta_G).
    \end{aligned}
    \end{align}
    We move on to the second order derivatives and justify differentiation under the integral by using the estimate \eqref{eq: polynomial moments} giving a polynomial bound for the integrands together with Lemma \ref{lem: poly moments}.
   \begin{align}\label{eq: second derivtaives in first step OAE}
       \begin{aligned}
        D^2_{\xi, \xi} [\Delta_\xi, \Delta_\xi] f_{x,y} (\xi, G) & = \, \Var^{\hat\mu(\cdot| y; \xi, G)}g_{x',y}(\xi, G; \Delta_\xi) \\
        D^2_{G, G} [\Delta_G, \Delta_G] f_{x,y} (\xi, G) & = \, \Var^{\hat\mu(\cdot| y; \xi, G)}h_{x',y}(\xi, G; \Delta_G)  \\
        D^2_{\xi, G} [\Delta_\xi, \Delta_G] f_{x,y} (\xi, G) & = \,  k_{x,y}(\xi; \Delta_\xi, \Delta_G) - \E^{\hat\mu(\cdot| y; \xi, G)} k_{x',y}( \xi; \Delta_\xi, \Delta_G)  \\
        & \quad + \Cov^{\hat\mu(\cdot| y; \xi, G)}\left(g_{x',y}( \xi,  G; \Delta_\xi), h_{x',y}(\xi,  G; \Delta_G)\right).
       \end{aligned}
   \end{align}
    Differentiation under the integral follows from the arguments above.
\end{proof}

Theorem \ref{thm: Expected KL divergence analytic} gives a framework to express and analyse both the first- and second-order optimality conditions through the study of stationary points. In addition to local optima, we next establish a convexity property of $\OAE$. For this purpose, recall that a function of two variables is bi-convex if it is convex in each variable when the other is held fixed.

\begin{theorem}\label{thm: biconvexity}
    Suppose Assumption \ref{ass: assumption} holds. The functional $\OAE$ is bi-convex in $(\xi, G)$. Further the Hessian at any critical point $(\hat \xi, \hat G) \in \R^m\times \Sm$ is positive semi-definite. 
    A critical point is attained at any $(\hat \xi, \hat G)$ satisfying
    \begin{eqnarray*}
        \E^{\mu_0} \left[\wA(x)\right] &=& \E^p \E^{\hat\mu(\cdot|y;\hat \xi, \hat G)} \left[\wA(x')\right] \\ 
        \E^\lambda \left[\wA(x) - y + \hat \xi \right]\otimes \left[\wA(x) - y + \hat\xi \right] &=& \E^p \E^{\hat\mu(\cdot|y;\hat \xi, \hat G)} \left[\wA(x') - y + \hat \xi\right] \otimes \left[\wA(x') - y + \hat \xi \right].
    \end{eqnarray*}
\end{theorem}

\begin{proof}
    Bi-convexity of $\OAE$ follows from the conditions in \eqref{eq: second derivtaives in first step OAE}: the second derivatives $D^2_{\xi, \xi}$ and $D^2_{G, G}$ of $\OAE$ are non-negative as integrals of non-negative functions.
    The form of a particular critical point follows from setting the derivatives $D_{\xi}$ and $D_{G}$ of $\OAE$ to zero and plugging in \eqref{eq: xi, G derivative}:
    \begin{align}\label{eq: derivative xi of OAE}
    \begin{aligned}
    D_\xi[\Delta_\xi] \OAE(\xi, G) &= \E^{\lambda} \left[(\wA(x) -y + \xi )^\top G \Delta_\xi \right] - \E^p \E^{\hat\mu(\cdot|y;\xi,  G)} \left[ (\wA(x') -y + \xi )^\top G \Delta_\xi \right] \\
    & =   \left( \E^{\mu_0} \left[\wA(x)\right] -\E^{p} \left[y \right] + \xi - \E^p \E^{\hat\mu(\cdot|y; \xi,  G)} \left[\wA(x')\right] +\E^{p} \left[y \right] - \xi \right)^\top G \Delta_\xi \\
    & =   \left( \E^{\mu_0} \left[\wA(x)\right] - \E^p \E^{\hat\mu(\cdot|y; \xi, G)} \left[\wA(x')\right] \right)^\top G \Delta_\xi \stackrel{!}{=} 0.
    \end{aligned}
    \end{align}
    For the other derivative:
    \begin{align}\label{eq: derivative G of OAE}
    \begin{aligned}
    D_G[\Delta_G] \OAE(\xi, G) &=  \frac{1}{2} \left( \E^{\lambda} \norm{\wA(x) - y + \xi}_{\Delta_{G}}^2 - \E^p \E^{\hat\mu(\cdot|y; \xi, G)}  \norm{\wA(x') - y + \xi}_{\Delta_{G}}^2\right) \stackrel{!}{=} 0.
    \end{aligned}
    \end{align}
We show now that for fixed $G \in \Sm$ and a critical point $\hat \xi$ of $\OAE(\cdot, G)$, the second order directional derivatives of $\OAE$ in any direction $(\Delta_\xi, \Delta_G) \in \R^m \times \Rsym$ are non-negative.
    The second directional derivative $D^2[\Delta_\xi, \Delta_G]$ of $\OAE$ at $(\hat \xi, G)$ is given by 
        \begin{align*}
            \begin{aligned}
            D^2[\Delta_\xi, \Delta_G] \OAE(\hat \xi, G) =& D^2_{\xi, \xi}[\Delta_\xi, \Delta_\xi] \OAE(\hat \xi, G) + D^2_{G, G}[\Delta_G, \Delta_G] \OAE(\hat \xi, G) + 2 D^2_{\xi, G}[\Delta_\xi, \Delta_G] \OAE(\hat \xi, G) \\
            = & \E^p \, \Var^{\hat\mu(\cdot| y; \hat \xi, G)}g_{x',y}(\hat \xi, G; \Delta_\xi) + \E^p \, \Var^{\hat\mu(\cdot| y; \hat \xi, G)}h_{x',y}(\hat \xi, G; \Delta_G) \\
            & + 2\, \E^\lambda \, k_{x,y}(\hat \xi; \Delta_\xi, \Delta_G) - 2 \, \E^p \,\E^{\hat\mu(\cdot| y; \hat \xi, G)} k_{x',y}( \hat \xi; \Delta_\xi, \Delta_G)  \\
        & + 2 \,\E^p \, \Cov^{\hat\mu(\cdot| y; \hat \xi, G)}\left(g_{x',y}( \hat \xi,  G; \Delta_\xi), h_{x',y}(\hat \xi,  G; \Delta_G)\right).
            \end{aligned}
        \end{align*}
    From the condition \eqref{eq: derivative xi of OAE} and the definition of $k$, it follows that for the critical point $\hat \xi$:
    \begin{eqnarray*}
         \E^\lambda k_{x,y}(\hat \xi; \Delta_\xi, \Delta_G) = \E^p  \E^{\hat\mu(\cdot| y; \hat \xi, G)} k_{x',y}( \hat \xi; \Delta_\xi, \Delta_G).
    \end{eqnarray*}
    For this reason, the second derivative (evaluated at $(\hat \xi, G)$) simplifies to
    \begin{align}\label{eq: second derivative at xi critical}
    \begin{aligned}
    D^2[\Delta_\xi, \Delta_G] \OAE(\hat \xi, G) =&  \E^p \, \Var^{\hat\mu(\cdot| y; \hat \xi, G)}g_{x',y}(\hat \xi, G; \Delta_\xi) + \E^p \, \Var^{\hat\mu(\cdot| y; \hat \xi, G)}h_{x',y}(\hat \xi, G; \Delta_G) \\
    & + 2 \,\E^p \, \Cov^{\hat\mu(\cdot| y; \hat \xi, G)}\left(g_{x',y}( \hat \xi,  G; \Delta_\xi), h_{x',y}(\hat \xi,  G; \Delta_G)\right) \\
    & = \E^p \, \Var^{\hat\mu(\cdot| y; \hat \xi, G)} \left[ g_{x',y}( \hat \xi,  G; \Delta_\xi) + h_{x',y}(\hat \xi,  G; \Delta_G) \right] \ge 0.
    \end{aligned}
    \end{align}
\end{proof}

Theorem \ref{thm: biconvexity} provides an interesting characterisation of the critical points. In the special case of linear models $\wA$ with a Gaussian prior $\mu_0$, one could derive explicit polynomial-type formulas for $\hat \xi$ and $\hat G$. In fact, in Theorem \ref{lem: Gaussian case; xi} it is shown that
\begin{eqnarray*}
    \hat{\xi} = \E^\lambda \left(y - \wA(x)\right),
\end{eqnarray*}
and hence $\hat{\xi}$ optimising $\OAE$ coincides with the optimiser of the EEA objective functional $\mathcal{K}$.
However, in the general case, the task is remarkably harder and is beyond the scope of this paper.

\section{Linear approximations with Gaussian prior}\label{sec: Gaussian linear}
In this section we specialise to the case of a Gaussian prior together with a linear surrogate model $\wA$. This setting is of particular interest for three reasons. First, it is arguably the most widely used case in applications, since linear-Gaussian models underlie classical Kalman-filter-type and Laplace-approximation methods, and often serve as a reasonable local approximation even when the true forward map is nonlinear. Second, in this setting the normalising constant $\hat Z$ admits a closed-form Gaussian expression, which allows for a simplified estimator in an empirical setting, see Section~\ref{subsec:empirical} ahead for details.
Finally, the linear-Gaussian structure allows us to say more: in Theorem~\ref{lem: Gaussian case; xi} below we derive a closed-form expression for $\OAE$ itself and, moreover, identify the global optimiser $\hat\xi$ explicitly. 

Assume now that $\wA$ is linear and continuous, $\Space=H$ a Hilbert space and $\mu_0 = \N(\overline{x}, C_0)$, with $\bar x \in H, C_0 \in L_+^1(H)$, where $L_+^1(H)$ denotes the space of positive, symmetric and injective operators of trace class. Then, the surrogate posterior distribution is given by $\hat \mu(\cdot | y; \xi, G) = {\mathcal N}(\hat{m}(y; \xi, G), \widehat{C}(G))$, where
\begin{align}
    \begin{aligned}
        \label{eq:linGauss_posterior}
\hat{m}(y;\xi, G) &= \bar x + C_0 \wA^* \left( G^{-1} + \wA C_0 \wA^*\right)^{-1}(y-\wA \bar x - \xi)
\\
\widehat{C}(G) &= C_0 - C_0\wA^*\left( G^{-1} + \wA C_0 \wA^*\right)^{-1} \wA C_0.
    \end{aligned}
\end{align}

In the following, we derive $\OAE$ for general Hilbert spaces $X = H$ and identify its optimiser $\hat \xi$.
\begin{theorem}\label{lem: Gaussian case; xi}
Let the Assumption \ref{ass: assumption} hold and the assumptions of this section and abbreviate $\textbf{M}(G) := G^{-1} + \wA C_0 \wA^*$ and $\textbf{B}(G):=C_0^{1/2}\wA^* \textbf{M}(G)^{-1}\wA C_0^{1/2}$. Then the objective functional can be written as
\begin{eqnarray*}
    \OAE(\xi, G) &=& \tfrac{1}{2} \log \det\left( I -  \textbf{B}(G) \right) + \tfrac{1}{2}\, \E^\lambda \norm{\wA\big(x - \hat m(y; \xi, G)\big)}_G^2 \\
    && +\; \tfrac{1}{2} \, \E^p \norm{C_0^{1/2}\wA^* \textbf{M}(G)^{-1}(y-\wA \bar x - \xi)}_{H}^2 \;-\; \E^\lambda \big\langle \wA^* \textbf{M}(G)^{-1} y,\; x - \bar x \big\rangle_H.
\end{eqnarray*}
For every fixed $G \in \Sm$, the minimum over $\xi \in \R^m$ is attained at
\begin{eqnarray}\label{eq: global opt xi}
    \hat \xi = \E^{\lambda} \left(y- \wA x\right),
\end{eqnarray}
independently of $G$, and the set of minimisers is the affine space $\hat\xi + \textbf{M}(G) \ker \wA^*$. In particular, $\hat \xi$ is the unique minimiser if and only if $\mathrm{Range}(\wA) = \R^m$.
\end{theorem}
Since the optimiser for $\xi$ of $\OAE$ coincides with the minimiser for the EEA, any advantage that OAE can achieve must stem from choosing the precision $G$ in a preferable way.
\begin{proof}[Proof of Theorem~\ref{lem: Gaussian case; xi}]
    By the Feldman-Hájek theorem, $\hat \mu(\cdot|y;\xi, G)$ and $\mu_0$ are equivalent with
    \begin{eqnarray*}
        \frac{d \hat \mu(\cdot|y;\xi, G)}{d \mu_0}(x) = \frac{d\N(\hat m, C_0)}{d \N(\bar x, C_0)}(x) \, \frac{d \N(\hat m, \hat C)}{d \N(\hat m, C_0)}(x),
    \end{eqnarray*}
    $\mu_0$ almost everywhere. We find that $\mu_0$ almost everywhere 
    \begin{eqnarray*}
        \frac{d\N(\hat m, C_0)}{d \N(\bar x, C_0)}(x) &= &\exp\left(-\tfrac{1}{2}\norm{\hat m-\bar x}_{C_0}^2 + \langle \hat m -\bar x, x - \bar x \rangle_{C_0}\right) 
    \end{eqnarray*}
    where $\langle \cdot, \cdot \rangle_{C_0}$ is defined as an extension of the white-noise mapping, see \cite{DaPratoStochasticEq} for details. 
    By using the property $\langle C_0 z, \cdot \rangle_{C_0} = \langle z, \cdot \rangle_H$ (see the arguments in \cite[Section 3.2]{Schneider2026diss}), $\mu_0$ almost everywhere in $z$, we are able to simplify
    \begin{eqnarray*}
        \frac{d\N(\hat m, C_0)}{d \N(\bar x, C_0)}(x) &= &\exp\big(-\tfrac{1}{2}\norm{C_0^{1/2} \wA^* M^{-1}(y-\wA \bar x - \xi)}_{H}^2 + \langle  \wA^* M^{-1}(y-\wA \bar x - \xi), x - \bar x \rangle_H\big).
    \end{eqnarray*}
    Notice that $\textbf{B}(G)$ is positive and of trace class. Further all eigenvalues of $\textbf{B}(G)$ are smaller than one since, by a brief calculation,
    \begin{eqnarray}\label{eq: I-BG inv}
        \left( I -  \textbf{B}(G) \right)^{-1} = C_0^{1/2} \wA^* G \wA C_0^{1/2}+ I,
    \end{eqnarray}
    is strictly positive definite.
    By \cite[Proposition 1.3.11]{DaPratoGiuseppe2002Sopd}, we can write for $\textbf{D}(G) := \textbf{B}(G) (I-\textbf{B}(G))^{-1} = (I-\textbf{B}(G))^{-1} -I$ that
    \begin{eqnarray}\label{eq: Feldman Hajek for covs}
         \frac{d \N(\hat m, \hat C)}{d \N(\hat m, C_0)}(x)  &=& \det\left( I -  \textbf{B}(G) \right)^{-1/2} \exp\left(-\tfrac{1}{2}\langle  \textbf{D}(G) C_0^{-1/2}(x-\hat m), C_0^{-1/2}(x-\hat m)\rangle_H\right),
    \end{eqnarray}
    where $\det$ denotes the Fredholm determinant. By plugging \eqref{eq: I-BG inv} into \eqref{eq: Feldman Hajek for covs} we find 
    \begin{eqnarray*}
        -\E^\lambda \log \left( \frac{d \N(\hat m, \hat C)}{d \N(\hat m, C_0)}(x) \right)  = \tfrac{1}{2} \log \det\left( I -  \textbf{B}(G) \right) +  \tfrac{1}{2} \E^\lambda \norm{\wA\left( x-\hat m \right)}_G^2.
    \end{eqnarray*}
    Introducing the notations $r = y - \wA \bar x - \xi$ and $\bar r = \E^p y - \wA \bar x - \xi$, we have
    \begin{eqnarray*}
        \OAE(\xi, G) = \tfrac{1}{2} \log \det\left( I -  \textbf{B}(G) \right) +  \tfrac{1}{2} \E^\lambda \norm{\wA\left( x-\hat m \right)}_G^2 + \tfrac{1}{2} \, \E^p \norm{C_0^{1/2}\wA^* M^{-1} r}_{H}^2  -  \E^\lambda \big\langle \wA^* M^{-1} y,\; x - \bar x \big\rangle_H.
    \end{eqnarray*}
    We first note that
    \begin{eqnarray*}
         \E^\lambda \wA\left( x-\hat m \right) = \E^\lambda \wA\left( x- \bar x - C_0 \wA^* M^{-1} r\right) &= - \wA C_0 \wA^* M^{-1} \bar r 
    \end{eqnarray*}
    By isolating the $\xi-$dependent terms and using the bias-variance decomposition we get
    \begin{eqnarray*}
        \E^p \norm{C_0^{1/2}\wA^* \textbf{M}(G)^{-1} r}_{H}^2 + \E^\lambda \norm{\wA\left( x-\hat m \right)}_G^2 &=& \norm{C_0^{1/2}\wA^* M^{-1} \bar r}_{H}^2 + \norm{\wA C_0 \wA^* M^{-1} \bar r}_G^2 \\
        & & + \tr\left( \Cov^\lambda(C_0^{1/2}\wA^* \textbf{M}(G)^{-1} y)\right) + \tr(G \Cov^\lambda \wA\left( x-\hat m \right)).
    \end{eqnarray*}
     The minimisers are given by $\bar r =  \E^\lambda (y-\wA(x)) - \xi \in \ker(\wA^* M^{-1})$, meaning that for $\hat \xi =\E^\lambda (y- \wA(x) )$ and any $G \in \Sm$ it holds
     \begin{eqnarray*}
         \argmin_{\xi \in \R^m } \OAE(\xi, G) = \hat \xi +  M \ker(\wA^*).
     \end{eqnarray*}
\end{proof}

The following remark can be helpful in evaluating the objective function in the case where prior and exact and approximate posterior distribution admit a known density even in the non-Gaussian-linear setting. 
\begin{remark}\label{rem: OAE for densities}
    Assume that the prior $\mu_0$ and the posterior distributions $\mu(\cdot|y), \hat \mu(\cdot|y; \xi, G)$ corresponding to \eqref{eq: inverse problem} and \eqref{eq: surrogate posterior} all admit densities, denoted by $p_0$, $q(\cdot|y)$ and $\hat q(\cdot|y; \xi, G)$, respectively. Then $\OAE$ satisfies
    \begin{eqnarray*}
        \OAE(\xi, G) &=& \E^\lambda \left[\log \left(\frac{d \mu(\cdot|y)}{d \hat \mu(\cdot|y; \xi, G)}(x)\right) - \log \left(\frac{d \mu(\cdot|y)}{d \mu_0}(x)  \right) \right] \\
        & = & \E^\lambda \left[ \log \left( \frac{ q(x|y) p_0(x)} {\hat q(x|y; \xi, G) q(x|y)}\right) \right] = \E^{\mu_0} \log p_0(x) - \E^\lambda \log \hat q(x|y; \xi, G).
    \end{eqnarray*}
\end{remark}
In particular, in the Gaussian-linear setting and when $X=\R^n$, the objective functional equals
\begin{equation}\label{eq: empirical expected KL divergence, linear surrogate}
      \OAE(\xi, G) = \tfrac{1}{2}\E^\lambda  \norm{x-\hat{m}(y; \xi, G)}^2_{\widehat{C}(G)^{-1}}   +\tfrac12 \log \det  \widehat{C}(G) + \E^{\mu_0} \log p_0(x) +\tfrac n2\log(2\pi).
\end{equation}
In finite dimensions the posterior mean and covariance satisfy
\begin{align}\label{eq: fin dim lin post}
    \begin{aligned}
\hat{m}(y;\xi, G) = \widehat{C}(G) \left( C_0^{-1} \bar x + \wA^\top G (y- \xi)\right), \quad 
\widehat{C}(G) = \left( C_0^{-1} + \wA^\top G \wA \right)^{-1}.
    \end{aligned}
\end{align}

\section{Optimisation of the regularised objective functional}
\label{sec:numerical_aspects}

This section considers practical optimisation of the functional $\OAE$. We first demonstrate the ill-posed nature of the optimisation problem, introduce suitable regularisation and convergent numerical method to minimise this regularised objective functional. 

\subsection{Regularised objective functional}
\label{subsec:regularisation}

We immediately observe three potential limitations in the problem formulation:
\begin{itemize}
    \item[(i)] \emph{Non-attainment of the infimum.} Let $X = \R$, $A = I_1
 = -\wA$, $\mu_0 = \N (0 , \tfrac 12)$, and $\delta = 1$. A direct computation 
    yields
    \begin{eqnarray*}
       \OAE(\xi, G) = \frac{3G (G+1) + \xi^2 G^2}{2(2+G)} + \frac{\log 2 }{2} - \frac{\log(2+G)}{2},
    \end{eqnarray*}
    from where it follows that
    $\argmin_{G \in \cl(\mathcal{S}(1))} \OAE(0, G) = 0$. Since $0 \notin \mathcal{S}(1)$, no minimiser exists over the admissible domain.

    \item[(ii)] \emph{Unbounded level sets for $G$.} Suppose $X = \R^k$ with $k<m$, $\mu_0$ is Gaussian and $\wA \in \R^{m\times k}$ is linear. 
    Clearly, $\wA$ cannot be surjective and $\wA^\top$ has a nontrivial null-space. It follows from \eqref{eq: fin dim lin post}, that for some $z\neq 0$ we have $\wA^\top z= 0$ and for every $t \ge 0$ and $\xi \in \R^m$,
    \begin{eqnarray}\label{eq: posterior const wrt G}
        \OAE(\xi, G + t z\otimes z) = \OAE(\xi, \, G).
    \end{eqnarray}
    Thus $\OAE$ is constant along the unbounded ray $\{G + t z\otimes z : t \ge 0\}$,
    and the level sets are unbounded. 
    \item[(iii)]\emph{Unbounded level sets for $\xi$.} Suppose again $X = \R^k$ with $k<m$, $\mu_0$ is Gaussian and $\wA \in \R^{m\times k}$ is linear. The previous degeneracy propagates to $\xi$: for fixed $G$ and every $t\in\R$,
\[
    \OAE(\xi + t\, G^{-1}z, \, G) = \OAE(\xi, \, G).
\]
Thus $\OAE(\cdot, G)$ is likewise constant along the unbounded line $\{\xi + t\,G^{-1}z : t\in\R\}$, and its level sets are unbounded as well.
\end{itemize}

The difficulties in (i), (ii) are resolved by adding an appropriate regularisation term for $G$, as we now show. We are able to resolve (iii) by imposing a suitable non-degeneracy condition, requiring that $\wA$ is non-constant under $\mu_0$ in every direction of $\R^m$. Further, we explain in Remark \ref{rem: why no reg in xi} that regularising $\xi$ may not be desirable.
We introduce the following regularisation functional, which is convex and coercive on $\Sm$ and penalises deviations from the true noise covariance.


\begin{lemma}\label{lem:regularisation}
The function
\[
    \Reg : \Sm \to \R, \qquad G \mapsto \delta^2 \tr (G) - \log \det G,
\]
has the following properties:
\begin{itemize}
    \item[(i)] $\Reg$ is convex and twice continuously differentiable on $\Sm$.
    \item[(ii)] Let $\{G_n\}_{n=1}^\infty \subset \Sm$ be a sequence such that $\gamma_{max}(G_n)\to \infty$ or $\gamma_{min}(G_n) \to 0$, where $\gamma_{max}$ and $\gamma_{min}$ stand for the maximal and the minimal eigenvalue, respectively.
    Then it follows that  
    $\Reg(G_n) \to \infty$.
    \item[(iii)] $\Reg$ has a unique minimiser $\delta^{-2} I_m \in \Sm$.
\end{itemize}
\end{lemma}

\begin{proof}
\emph{(i) Convexity.} The map $G \mapsto \tr G$ is linear, hence convex. The map
$G \mapsto -\log \det G$ is convex on $\Sm$ (it is a standard result that $-\log \det$
is convex; see, e.g., \cite{Boyd2004ConvexOptimization}). As a sum of convex functions,
$\Reg$ is convex and twice continuously differentiable. Differentiability follows from the differentiability of both summands on the open set
$\Sm$.

\emph{(ii) Coercivity.} Write $\Reg(G) = \sum_{i=1}^m s(\gamma_i(G))$, where $\gamma_i(G)$ are the eigenvalues of $G$ and $s(\gamma) := \delta^2 \gamma - \log \gamma$. The function $s$ is bounded from below by $ s(\delta^{-2})$ and satisfies $s(\gamma) \to \infty$ both as $\gamma \to 0^+$ and as $\gamma \to \infty$. Hence, if $\gamma_{max}(G_n) \to \infty$, then $\Reg(G_n) \ge (m-1)s(\delta^{-2}) + s(\gamma_{max}(G_n)) \to \infty$, and analogously if $\gamma_{min}(G_n) \to 0$.

\emph{(iii) Unique minimiser.} Setting the derivative to zero, the first-order optimality
condition for $G \mapsto \delta^2 \tr G - \log \det G$ is $\delta^2 I - G^{-1} = 0$,
giving $G = \tfrac{1}{\delta^2} I$. Strict convexity of $-\log \det$ on $\Sm$ (its
Hessian is positive definite) guarantees this critical point is the unique global minimiser.
\end{proof}

We emphasise that $\Reg(G) = 2 \, \KL{\N(0, \delta^2I)}{\N(0, G^{-1})} + m ( 1+ \log \delta^2)$ and it is therefore straightforward to generalise the regulariser to settings where observational noise in the inverse problem is non-isotropic.

For $\alpha > 0$, define the regularised objective function
\[
    \OAE_\alpha : \R^m \times \Sm \to \R, \qquad
    \OAE_\alpha(\xi, G) \;=\; \OAE(\xi, G) + \alpha\, \Reg(G).
\]

The regularisation parameter $\alpha > 0$ controls the strength of the regularisation. By Lemma~\ref{lem:regularisation}(ii), the sublevel sets of $\OAE_\alpha$ are bounded in the $G$-direction and bounded away from the boundary of $\Sm$, thereby resolving the first two issues identified in the beginning of this section. Boundedness in the $\xi$-direction requires the non-degeneracy condition of Assumption~\ref{ass: assumption conv of AM} below and is established in Lemma~\ref{lem: coercive decomp}.

\subsection{Proximal Alternating Minimisation}
\label{subsec:proximal_alternating}

We adopt a proximal alternating minimisation strategy, which interleaves exact minimisation
over each block variable with a quadratic proximal penalty that stabilises the iterates.
Convergence of such algorithms under the Kurdyka-\L{}ojasiewicz (K\L{}) property has been
established in \cite{attouch2013proximal, Xu2013ABlockCoordinate}; we recall the relevant
definitions before stating our results.

Denote by $\partial f$ the sub-differential set of $f$ and by $\dom(\partial f)$ the domain of the sub-differential i.e. the points where the sub-differential is defined. In what follows, our objective functional is extended outside its natural domain by $\infty$ and differentiable, in which case $\dom(\partial f)$ coincides with its natural domain, see Remark~\ref{rem: lsc extension}.
\begin{definition}[Kurdyka-\L ojasiewicz property {\cite{attouch2013proximal}}]
\label{def:KL}
A lower semi-continuous function $f : \R^\ell \to \R \cup \{+\infty\}$
satisfies the \emph{K\L{} property} if for every $z^* \in \dom(\partial f)$ there
exist a constant $\eta \in (0, +\infty)$, a neighbourhood $\mathcal{U}$ of $z^*$,
and a continuous concave function $\varphi\colon [0, \eta) \to \R_{\geq 0}$ such that:
\begin{itemize}
    \item[(i)] $\varphi(0) = 0$,
    \item[(ii)] $\varphi \in C^1(0, \eta)$ with $\varphi' > 0$ on $(0, \eta)$,
    \item[(iii)] for all $z \in \mathcal{U}$ with
    $f(z^*) < f(z) < f(z^*) + \eta$, the \emph{K\L{} inequality} holds:
    \[
        \varphi'\!\bigl(f(z) - f(z^*)\bigr)
        \cdot \dist\bigl(0,\, \partial f(z)\bigr) \;\geq\; 1.
    \]
\end{itemize}
\end{definition}

The K\L{} inequality quantifies how fast the function values separate from the critical value, and is the key analytic tool in the convergence proofs of \cite{attouch2013proximal}.

\begin{remark}\label{rem:KL_examples}
The following classes of functions satisfy the K\L{} property and will be relevant below.
\begin{itemize}
    \item[(a)] \emph{Morse functions.} A $C^2$ function $f$ is called a \emph{Morse function}
    if the Hessian $D^2 f(\bar x)$ is invertible at every critical point $\bar x$. This holds in particular when $D^2 f(\bar x)$ is strictly positive definite at every critical point, which we will verify.
    Every Morse function satisfies the K\L{} property \cite{attouch2013proximal}.
    \item[(b)] \emph{Non-critical points.} Any proper lower semi-continuous function
    satisfies the K\L{} property at every non-critical point \cite{attouch2013proximal}.
\end{itemize}
\end{remark}

We now define the proximal alternating minimisation sequence. 
\begin{definition}[Proximal alternating minimisation]\label{def:proximal_alternating}
Let $f : \R^{\ell_1} \times \R^{\ell_2} \to \R \cup \{+\infty\}$
be an objective functional. 
Given an initial state $(z_0, w_0) \in \R^{\ell_1} \times \R^{\ell_2}$, the proximal alternating minimisation generates the sequence $(z_k, w_k)_{k \geq 1}$ defined
recursively by
\begin{align*}
    z_k &\;\in\; \argmin_{z \in \R^{\ell_1}}
        \left[ f(z, w_{k-1}) + \frac{1}{\zeta_k^1} \|z - z_{k-1}\|^2 \right], \\
    w_k &\;\in\; \argmin_{w \in \R^{\ell_2}}
        \left[ f(z_k, w) + \frac{1}{\zeta_k^2} \|w - w_{k-1}\|^2 \right],
\end{align*}
where the proximal parameters satisfy $\zeta_k^1, \zeta_k^2 \in (r_-, r_+)$ for fixed constants
$0 < r_- < r_+$.
\end{definition}

\begin{remark}\label{rem:proximal_interpretation}
The proximal terms $\frac{1}{\zeta_k^1}\|z - z_{k-1}\|^2$ and $\frac{1}{\zeta_k^2}\|w - w_{k-1}\|^2$
serve a dual purpose: they ensure that each subproblem has a unique solution, and they
prevent large, potentially destabilising steps. For large $\zeta_k^1, \zeta_k^2$ the proximal
penalty is negligible and the algorithm reduces to standard alternating minimisation. For
moderate values, the algorithm resembles an alternating proximal gradient method; see the
discussion in \cite{attouch2013proximal}.
\end{remark}

The central result of \cite{attouch2013proximal} states that, if the objective functional satisfies the K\L{} property at every point of its domain, then (under mild assumptions) any bounded sequence generated by the proximal alternating minimisation algorithm converges to a critical point of $f$. We include this result here for the reader's convenience.
\begin{theorem}[{\cite[Theorems 8--10]{attouch2013proximal}}]\label{thm:convergence_dichotomy}
Let $f :\R^{\ell_1} \times \R^{\ell_2} \to \R \cup \{+\infty\}$ be of the form $f(z, w) = f_1(z) + Q(z, w) + f_2(w)$, where $f_1, f_2$ are proper and lower semi-continuous and $Q$ is continuously differentiable with gradient Lipschitz continuous on bounded sets, and assume that $f$ is bounded from below and $f(\cdot, w_0)$ is proper. Let $(z_k, w_k)$ be the sequence generated by Definition~\ref{def:proximal_alternating} with $\zeta_k^1, \zeta_k^2 \in (r_-, r_+)$.
\begin{itemize}
    \item If $(z_k, w_k)$ is bounded and $f$ has the K\L{} property at every cluster point of $(z_k, w_k)$, then $\sum_{k=1}^\infty \|(z_k, w_k) - (z_{k-1}, w_{k-1})\| < \infty$ and $(z_k, w_k)$ converges to a critical point of $f$. 
    \item Furthermore, if $(\hat z, \hat w)$ is a global minimiser of $f$ at which $f$ has the K\L{} property, then any sequence initialised sufficiently close to $(\hat z, \hat w)$ converges to a global minimiser of $f$.
\end{itemize}
\end{theorem}
The first statement is the content of the proof of \cite[Theorem~9]{attouch2013proximal}, which applies \cite[Theorem~8]{attouch2013proximal} at a cluster point of the iterates; the second is \cite[Theorem~10]{attouch2013proximal}.

We note that convergence to the global optimiser can only be guaranteed if the iterates are initialised  close to the global optimiser. Since the objective function in the following may have multiple local optima, we will not be able to ensure convergence to global optimisers. 

As we see in the following lemma, when $f$ is coercive, the case where the sequence diverges can be ruled out (compare with \cite[Remark 5]{attouch2013proximal}). 
\begin{lemma}\label{lem:iterates_bounded}
Let the assumptions of Theorem \ref{thm:convergence_dichotomy} hold and assume that $f$ is coercive.
Then the sequence $(z_k, w_k)$ generated by Definition~\ref{def:proximal_alternating} and with $\zeta_k^1, \zeta_k^2 \in (r_-, r_+)$ is bounded.
\end{lemma}
\begin{proof}
    By \cite[Lemma 5]{attouch2013proximal}, the sequence $f(z_k, w_k)$ is non-increasing.
   If the sequence $(z_k, w_k)$ were unbounded, then $\|(z_{n_k}, w_{n_k}) \| \to \infty$ along some subsequence $n_k$.
   Since $f$ is coercive, $f(z_{n_k}, w_{n_k}) \rightarrow \infty$, which is a contradiction with the sequence being non-increasing. 
\end{proof}

\begin{remark}\label{rem: lsc extension}
Theorem~\ref{thm:convergence_dichotomy} is stated for functions on $\R^{\ell_1} \times \R^{\ell_2}$ with $Q \in C^1$ everywhere, whereas $\OAE_\alpha$ is defined only on the open set $\R^m \times \Sm$. The idea is that the regulariser confines the iterates to a compact subset of $\Sm$, so that $\OAE$ may be modified at will away from this set to meet the hypotheses of the theorem without changing the sequence the algorithm generates.

We reduce to that setting as follows. Throughout, we identify $\Rsym$ with $\R^{m'}$, $m' := m(m+1)/2$, via the upper-triangular entries, so that $\Sm$ is an open convex subset of $\R^{m'}$. Set $L_0 := \OAE_\alpha(\xi_0, G_0)$, $\kappa := \E^p\KL{\mu(\cdot|y)}{\mu_0}$ and $\mathcal B := \{G \in \Sm : \alpha\Reg(G) \le L_0 + \kappa\}$, which is compact in $\Sm$ by Lemma~\ref{lem:regularisation}(ii). Since the objective values are non-increasing along the iteration \cite{attouch2013proximal} and $\OAE \ge -\kappa$ by Lemma~\ref{lem: assumptions}, every iterate $G_k$ lies in $\mathcal B$; in particular the iterates never approach $\partial\Sm$. Choose a smooth cut-off $\chi \in C^\infty(\R^{m'})$ with compact support in $\Sm$ and $\chi \equiv 1$ on a neighbourhood $\mathcal N$ of $\mathcal B$, and set $f_1 \equiv 0$, $Q := \chi\OAE - (1-\chi)\kappa \in C^2(\R^m \times \R^{m'})$, and $f_2 := \alpha\Reg$ extended by $+\infty$ outside $\Sm$, which is proper and lower semi-continuous by Lemma~\ref{lem:regularisation}(ii). Then $Q \ge -\kappa$, $Q = \OAE$ on $\R^m \times \mathcal N$, and $Q + f_2 > L_0$ on $\R^m \times (\Sm \setminus \mathcal B)$, so the proximal subproblems for $f_1 + Q + f_2$ and for $\OAE_\alpha$ have the same minimisers, and the common iterates have all cluster points in $\R^m \times \mathcal B$, where the two functionals coincide. We suppress this reduction in what follows.
\end{remark}

 To establish the Morse property for $\OAE_\alpha$ and hence convergence of the proximal alternating scheme, we require the following mild non-degeneracy condition, which in turn will ensure that the second derivatives in \eqref{eq: second derivative at xi critical} do not vanish.

\begin{assumption}\label{ass: assumption conv of AM}
No non-zero vector $v\in\R^m$ exists such that $v \cdot \wA$ is $\mu_0$-almost surely constant.
\end{assumption}

\begin{remark}\label{rem:linear_gaussian}
If $\wA\colon \R^n \to \R^m$ is linear, then Assumption~\ref{ass: assumption conv of AM} holds if and only if
$\mathrm{Range}(\wA) = \R^m$, i.e., $\wA$ is surjective. This is a natural identifiability condition: it requires that the observation model $\wA$ can distinguish all directions in $\R^m$.
\end{remark}
\begin{remark}\label{rem:existence_minimum}
Under the Assumptions~\ref{ass: assumption} and~\ref{ass: assumption conv of AM}, the function $\OAE_\alpha$ is coercive. Coercivity of $\OAE_\alpha(\,\cdot\,,G)$ in
$\xi$ is shown in Lemma~\ref{lem: coercive decomp} ahead and coercivity of $\Reg$ in $G$ together imply that $\OAE_\alpha(\xi, G) \to \infty$ whenever $\|(\xi, G)\| \to \infty$ or $G$ approaches the boundary of $\Sm$. This can be used to rule out the case of divergence of $(\xi_k, G_k)$ as follows. If the iterates $(\xi_k, G_k)$ were unbounded, then also $\OAE_\alpha(\xi_k, G_k)$ would have an accumulation point at infinity, a contradiction to the estimates $\OAE_\alpha(\xi_k, G_k)$ being non-increasing.

Even further, being continuous and coercive, $\OAE_\alpha$ attains its minimum, so a global minimiser exists (it need not be unique). Moreover, each subproblem in Definition~\ref{def:proximal_alternating} is strongly convex and coercive, since $\OAE_\alpha$ is bi-convex, so the subproblem minimisers $\xi_k$ and $G_k$ are uniquely defined at every step.
\end{remark}
\begin{remark}
    If Assumption \ref{ass: assumption conv of AM} is violated, then we will show in Proposition~\ref{prop:iterates_bounded}, that there exists an unbounded sequence $\xi_k \subset \R^m$ with the property $\OAE(\xi_k, G) = \OAE(\xi_0, G)$ for given $G$, meaning that no optimal $\xi$ over the entirety of $\R^m$ can be identified through $\OAE$. 
\end{remark}

\begin{lemma}\label{lem: cov strictly psd}
The covariance matrix
\[
    C = \mathbb{E}^{\mu}\bigl[(\wA(x) - \mathbb{E}^{\mu} \wA(x))\otimes (\wA(x) - \mathbb{E}^{\mu} \wA(x))\bigr],
\]
where $\mu = \hat \mu(\cdot | y ; \xi, G)$,
is positive definite if and only if no non-zero vector $v\in\R^m$ exists such that $v \cdot \wA$
is $\mu_0$-almost surely constant.
\end{lemma}
 
\begin{proof}
First observe that the Radon-Nikodym derivative
\begin{equation}
\frac{d\mu}{d\mu_0}(x | y; \xi, G) = \frac 1{\widehat{Z}(y; \xi, G)} \exp\left(-\tfrac12 \norm{\wA(x) - y + \xi}_G^2\right) > 0
\end{equation}
for any $x\in X$, $y,\xi\in \R^m$ and $G\in \Sm$. Therefore, if for any Borel measurable set ${\mathcal B} \subset X$, we have $\mu({\mathcal B} | y; \xi, G) = 0$, it follows that $\mu_0({\mathcal B}) = 0$ and, consequently, $\mu_0$ is absolutely continuous with respect to $\mu(\cdot | y; \xi, G)$.
Therefore, the assumption of the claim is equivalently stated for $\mu(\cdot | y; \xi, G)$ instead of $\mu_0$.

For any $v \in \mathbb{R}^m$,
\[
    v^\top C v = \mathrm{Var}^{\mu}(v \cdot \wA(x))
    = \| v \cdot \wA(x) - \mathbb{E}^\mu[v\cdot \wA(x)] \|_{L^2(\mu)}^2.
\]
This vanishes if and only if $v \cdot \wA$ is $\mu$-a.s.\ constant.
Hence $v^\top Cv > 0$ for all $v \neq 0$, i.e.\ $C \succ 0$,
if and only if no nontrivial $v \cdot \wA$ is $\mu$-a.s.\ constant.
\end{proof}

\begin{lemma}\label{lem: Morse property}
Let Assumptions \ref{ass: assumption} and \ref{ass: assumption conv of AM} be satisfied. Then the functional $\OAE_\alpha$ satisfies the assumptions of Theorem \ref{thm:convergence_dichotomy} and every critical point of $\OAE_\alpha$ is a local minimum.
\end{lemma}
\begin{proof}
\hfill
\begin{itemize}
    \item  By Theorem \ref{thm: Expected KL divergence analytic} and Lemma \ref{lem:regularisation}, $\OAE_\alpha = \OAE +\alpha \Reg $ is twice continuously differentiable as a sum of twice continuously differentiable functionals.
    \item  $\OAE_\alpha$ is globally bounded from below since $\OAE$ is bounded from below by Lemma \ref{lem: assumptions} and $\Reg$ is bounded from below as a coercive and convex function (see Lemma \ref{lem:regularisation}).
    \item  $\OAE_\alpha$ is Morse, which implies the K\L{} property by Remark \ref{rem:KL_examples}. 
    Let $(\hat \xi, \hat G) \in \R^m \times \R^{m'}$ be a critical point of $\OAE_\alpha$. We show that the second derivative in a direction $(\Delta_\xi, \Delta_G)$ is strictly positive. In the case where $\Delta_G = 0$, we notice that $\OAE(\cdot, G) - \OAE_\alpha(\cdot, G)$ is constant for any fixed $G \in \Sm$. Hence under utilisation of \eqref{eq: second derivative at xi critical} we find
    \begin{eqnarray*}
        D^2[\Delta_\xi, 0] \OAE_\alpha(\hat \xi, \hat G) &=& D^2[\Delta_\xi, 0] \OAE(\hat \xi, \hat G) \\
        & = & \E^p \, \Var^{\hat\mu(\cdot| y; \hat \xi, \hat G)} \left[ \wA(x)^\top \hat G \Delta_\xi \right].
    \end{eqnarray*}
    By Lemma \ref{lem: cov strictly psd} applied with $v = \hat G \Delta_\xi \neq 0$, in connection with Assumption \ref{ass: assumption conv of AM}, $D^2[\Delta_\xi, 0] \OAE_\alpha(\hat \xi, \hat G) > 0$ for $\Delta_\xi \in \R^m \setminus \{0\}$, which is precisely what we needed to show. In the case where $\Delta_G \neq 0$,
    \begin{eqnarray*}
        D^2[\Delta_\xi, \Delta_G] \OAE_\alpha(\hat \xi, \hat G) &=& D^2[\Delta_\xi, \Delta_G] \OAE(\hat \xi, \hat G) + \alpha D^2[\Delta_G, \Delta_G]\Reg(\hat G) \\
        & = & \E^p \, \Var^{\hat\mu(\cdot| y; \hat \xi, \hat G)} \left[ \wA(x)^\top \hat G \Delta_\xi + \tfrac{1}{2} \norm{\wA(x)-y+\hat \xi}^2_{\Delta_G}\right] +  \alpha D^2[\Delta_G, \Delta_G]\Reg(\hat G) \\
        & \ge & \alpha D^2[\Delta_G, \Delta_G]\Reg(\hat G)  > 0,
    \end{eqnarray*}
 where it was used that $\alpha D^2[\Delta_G, \Delta_G]\Reg(\hat G)  > 0$ for $\Delta_G \neq 0$, since $\Reg$ is strictly convex. So indeed $\OAE_\alpha$ is Morse.
\end{itemize}
\end{proof}

\begin{remark}\label{rem: why no reg in xi}
    The proof of Lemma~\ref{lem: Morse property} crucially relies on the fact that the cross term involving $k$ vanishes at critical points of $\OAE$ in $\xi$, see \eqref{eq: second derivative at xi critical}. Regularising $\xi$ would shift these critical points, at which the cancellation identity \eqref{eq: second derivative at xi critical} no longer holds. Hence the cross term may contribute negatively to the Hessian, and positive definiteness, and hence the Morse property, can fail.
\end{remark}

We can now define a proximal alternating minimisation sequence corresponding to the minimisation of $\OAE_\alpha$ as follows. Let $\xi_0, G_0 \in \R^m \times \Sm$ be given and fix $0 <r_- <r_+ $ such that $\zeta_k^1, \zeta_k^2 \in (r_-, r_+)$. For $k \ge 1$ define
\begin{align*} 
    \xi_k &\;\in\; \argmin_{\xi \in \R^{m}}
        \left[ \OAE_\alpha(\xi, G_{k-1}) + \frac{1}{\zeta_k^1} \|\xi - \xi_{k-1}\|^2 \right], \\
    G_k &\;\in\; \argmin_{G \in \Sm}
        \left[\OAE_\alpha(\xi_k, G) + \frac{1}{\zeta_k^2} \|G - G_{k-1}\|^2 \right].
\end{align*}

Since we may use Theorem~\ref{thm:convergence_dichotomy}, the sequence $(\xi_k, G_k)$ converges to a critical point if it is bounded. By Lemma \ref{lem:iterates_bounded}, we can show that the objective function is coercive to ensure boundedness of the iterates and hence convergence to a critical point.. For $G$, this follows from the regularisation, but coercivity in $\xi$ requires a bit more work, as we will see in the following lemma. We will present the statement in a general version, which will be utilised for the measures $(p, \lambda)$ but also their empirical versions such as $(\tfrac{1}{M} \sum_{i=1}^M \delta_{y_i}, \tfrac{1}{M} \sum_{i=1}^M \delta_{x_i, y_i})$ for suitable data, as we will see in Section~\ref{subsec:empirical}.

\begin{lemma}[Coercive decomposition]\label{lem: coercive decomp}
Let $\tilde{\lambda}_M$ be a sub-Gaussian probability measure on $\Space \times \R^m$ and let $\tilde{p}_M$ denote its marginal on $\R^m$. Let $\ubar{\wA}^{\tilde{\lambda}_M}:= \E^{\tilde{\lambda}_M}\wA(x)$, denote by $\langle \cdot, \cdot \rangle^+$ the positive part of $\langle \cdot, \cdot \rangle$ and
    \begin{eqnarray*}
        C_*(\tilde \lambda_M) :=  \inf_{ \norm{v}_G^2 = 1, \, \gamma(G) \in [\gamma_{min}, \gamma_{max}] }
  \E^{\mu_0}\!\left[
    \exp\left( -\tfrac{1}{2} \norm{ \wA(x) - \ubar{\wA}^{\tilde{\lambda}_M}}_{G}^2\right)\langle \ubar{\wA}^{\tilde{\lambda}_M} -\wA(x) ,\, v\rangle_G^+ 
  \right].
    \end{eqnarray*}
For $G \in \Sm$ with eigenvalues in $[\gamma_{min}, \gamma_{max}]$ it holds
\begin{align}
\begin{aligned}\label{eq: tilde OAE}
        \tilde \OAE^M(\xi, G) &:=& \tfrac{1}{2}\E^{\tilde{\lambda}_M} \norm{\wA(x) -y + \xi}_G^2 + \E^{\tilde{p}_M} \log \E^{\mu_0} \exp\left( -\tfrac{1}{2} \norm{\wA(x) -y + \xi}_G^2 \right) \\
    & \ge &  - \E^{\tilde{\lambda}_M} \langle \wA(x) - \ubar{\wA}^{\tilde{\lambda}_M}, y\rangle_G + \E^{\tilde{p}_M} \log \left(\norm{\ubar{\wA}^{\tilde{\lambda}_M}-y+\xi}_G C_*(\tilde \lambda_M) \right).
\end{aligned}
\end{align}
\end{lemma}
\begin{proof}
    We have
\begin{equation*}
    \norm{\wA(x) -y + \xi}_G^2 =  \norm{\wA(x) -\ubar{\wA}^{\tilde{\lambda}_M}}_G^2 + \norm{\ubar{\wA}^{\tilde{\lambda}_M}-y+\xi}_G^2 +2\langle \wA(x) -\ubar{\wA}^{\tilde{\lambda}_M}, \ubar{\wA}^{\tilde{\lambda}_M}-y+ \xi\rangle_G.
\end{equation*}
Taking expectations and using that $\norm{\ubar{\wA}^{\tilde{\lambda}_M}-y+\xi}_G^2$ does not depend on $x$ and may therefore be taken out of the inner expectation over $\mu_0$, we obtain
\begin{eqnarray*}
    \tilde \OAE^M(\xi, G)  &=& \tfrac{1}{2} \E^{\tilde{\lambda}_M} \norm{ \wA(x) -\ubar{\wA}^{\tilde{\lambda}_M}}_G^2 +  \underbrace{\E^{\tilde{\lambda}_M} \langle \wA(x) - \ubar{\wA}^{\tilde{\lambda}_M}, \ubar{\wA}^{\tilde{\lambda}_M} + \xi \rangle_G}_{=0} - \E^{\tilde{\lambda}_M} \langle \wA(x) - \ubar{\wA}^{\tilde{\lambda}_M}, y\rangle_G\\
    & & + \tfrac12 \E^{\tilde{\lambda}_M} \norm{\ubar{\wA}^{\tilde{\lambda}_M}-y+\xi}_G^2 - \tfrac12 \E^{\tilde{p}_M} \norm{\ubar{\wA}^{\tilde{\lambda}_M}-y+\xi}_G^2 \\
    & & + \E^{\tilde{p}_M}\log \E^{\mu_0} \exp\left( -\tfrac{1}{2} \norm{ \wA(x) -\ubar{\wA}^{\tilde{\lambda}_M}}_G^2  - \langle \wA(x) -\ubar{\wA}^{\tilde{\lambda}_M}, \ubar{\wA}^{\tilde{\lambda}_M}-y+\xi\rangle_G \right).
\end{eqnarray*}
Since $\tilde{p}_M$ is the marginal of $\tilde{\lambda}_M$, the two terms in the second line cancel, and the first term is non-negative. Hence
\begin{eqnarray*}
    \tilde \OAE^M(\xi, G)  &\ge& - \E^{\tilde{\lambda}_M} \langle \wA(x) - \ubar{\wA}^{\tilde{\lambda}_M}, y\rangle_G  \\
    & & + \E^{\tilde{p}_M}\log \E^{\mu_0} \exp\left( -\tfrac{1}{2} \norm{ \wA(x) -\ubar{\wA}^{\tilde{\lambda}_M}}_G^2  - \langle \wA(x) -\ubar{\wA}^{\tilde{\lambda}_M}, \ubar{\wA}^{\tilde{\lambda}_M}-y+\xi\rangle_G \right).
\end{eqnarray*}
We introduce
\begin{eqnarray*}
    (\ubar{\wA}^{\tilde{\lambda}_M}-y+\xi) = r\hat{v}, \quad r := \norm{\ubar{\wA}^{\tilde{\lambda}_M}-y+\xi}_G, \quad \hat{v} := (\ubar{\wA}^{\tilde{\lambda}_M}-y+\xi)/\norm{\ubar{\wA}^{\tilde{\lambda}_M} -y+ \xi}_G.
\end{eqnarray*}
By using $e^t \ge t^+$ for all $t \in \R$, we get
\begin{align}
\begin{aligned}\label{eq: coercive step2}
     &\log \E^{\mu_0} \exp\left( -\tfrac{1}{2} \norm{ \wA(x) -\ubar{\wA}^{\tilde{\lambda}_M}}_G^2  - \langle \wA(x) -\ubar{\wA}^{\tilde{\lambda}_M}, \ubar{\wA}^{\tilde{\lambda}_M}-y+\xi\rangle_G \right) \\
     & \ge  \log \E^{\mu_0} \exp\left( -\tfrac{1}{2} \norm{ \wA(x) -\ubar{\wA}^{\tilde{\lambda}_M}}_G^2 \right) \langle \ubar{\wA}^{\tilde{\lambda}_M}-\wA(x) , \ubar{\wA}^{\tilde{\lambda}_M}-y+\xi\rangle_G^+ \\
    & = \log \left[ r \, \E^{\mu_0} \exp\left( -\tfrac{1}{2} \norm{ \wA(x) -\ubar{\wA}^{\tilde{\lambda}_M}}_G^2 \right) \langle \ubar{\wA}^{\tilde{\lambda}_M}-\wA(x) , \hat{v}\rangle_G^+ \right] \ge \log r + \log C_*(\tilde \lambda_M). 
\end{aligned}
\end{align}
\end{proof}

\begin{proposition}\label{prop:iterates_bounded}
Suppose Assumption \ref{ass: assumption} is satisfied. If Assumption \ref{ass: assumption conv of AM} is additionally satisfied, then the iterates $(\xi_k, G_k)$ remain bounded and converge to a critical point of $\OAE_\alpha$. If Assumption \ref{ass: assumption conv of AM} is violated, then for any $\xi_0 \in \R^m$ there exists an unbounded sequence $(\xi_n)_{n \in \Natural}$ such that $\OAE(\xi_n, G) = \OAE(\xi_0, G)$.
\end{proposition}
\begin{proof}
Let first Assumption \ref{ass: assumption conv of AM} hold. By Lemma \ref{lem:iterates_bounded}, we need to show that $\OAE_\alpha$ is coercive, i.e. it holds $\OAE_\alpha(\xi_k, G_k) \to \infty$ for any $(\xi_k, G_k) \in \R^m \times \Sm$ with $(\xi_k, G_k) \to \infty$. 
\begin{itemize}
    \item If $\gamma_{max}(G_k) \to \infty$ or $\gamma_{min}(G_k) \to 0$, then $\Reg(G_k) \rightarrow \infty$ by Lemma~\ref{lem:regularisation}(ii). By Lemma \ref{lem: assumptions}, it holds
    \begin{eqnarray*}
        \alpha \Reg(G_k) - \E^p \KL{\mu(\cdot | y)}{\mu_0} \le \OAE_\alpha(\xi_k, G_k) \to \infty.
    \end{eqnarray*}
    \item Assume now $\xi_k \to \infty$ and $G_k$ are bounded, such that the eigenvalues of $G_k$ for all $k\in \Natural$ lie in the interval $[\gamma_{min}, \gamma_{max}]$. Observe carefully that the eigenvalues of $G_k$ are uniformly bounded from below, since otherwise the regulariser goes to infinity. 
    Recall that
    \begin{eqnarray*}
        C_*(\lambda) =  \inf_{ \norm{v}_G^2 = 1, \, \gamma(G) \in [\gamma_{min}, \gamma_{max}] }
  \E^{\mu_0}\!\left[
    \exp\left( -\tfrac{1}{2} \norm{ \wA(x) - \ubar{\wA}^{\mu_0}}_{G}^2\right)\langle \ubar{\wA}^{\mu_0} -\wA(x) ,\, v\rangle_G^+ 
  \right].
    \end{eqnarray*}
We show now that $C_*>0$. Notice first that the infimum of continuous functions over compact sets are attained. Hence $C_* = 0$ would imply that there exists $v\in \R^m$ such that $\langle \ubar{\wA}^{\mu_0}- \wA(x) ,\, v\rangle_G \le 0$ almost surely. Since also $\E^{\mu_0} \langle \ubar{\wA}^{\mu_0} -\wA(x) ,\, v\rangle_G =0$, this implies $\langle \ubar{\wA}^{\mu_0} -\wA(x) ,\, v\rangle_G = 0$ almost surely, a contradiction with Assumption \ref{ass: assumption conv of AM}. 

By Lemma~\ref{lem: coercive decomp} for $(p, \lambda)$, we get that for any $\xi_k \to \infty$ 
\begin{eqnarray*}
\OAE_\alpha(\xi_k, G_k) \ge \alpha \Reg(\delta^{-2} I) + \E^p \log \left(C_*(\lambda)  \sqrt{\gamma_{min}} \norm{\ubar{\wA}^{\mu_0}-y+ \xi_k}_2\right)  - \E^\lambda \langle \wA(x) -\ubar{\wA}^{\mu_0}, y\rangle_{G_k}.
\end{eqnarray*}
The last term is bounded uniformly in $k$, since the eigenvalues of $G_k$ lie in $[\gamma_{min}, \gamma_{max}]$ and $\lambda$ has finite second moments. For the logarithmic term, on the event $$E_k := \left\{\norm{y - \ubar{\wA}^{\mu_0}}_2 \le \norm{\xi_k}_2/2\} \right\}$$ we have $\norm{\ubar{\wA}^{\mu_0}-y+ \xi_k}_2 \ge \norm{\xi_k}_2/2$, while $\sup_{c \in \R^m} \E^p \log^- \norm{y - c}_2 < \infty$ because the density $p$ is bounded. Since $p(E_k) \to 1$, it follows that $\E^p \log \norm{\ubar{\wA}^{\mu_0}-y+ \xi_k}_2 \to \infty$, and hence $\OAE_\alpha(\xi_k, G_k) \to \infty$.
\end{itemize}
If Assumption \ref{ass: assumption conv of AM} is violated, then there exists $v \neq 0$ such that $\mu_0$ almost surely, $\langle \wA(x) - \ubar{\wA}^{\mu_0}, v \rangle_G = 0$. If we choose $\xi_r:=r v+ \xi_0$ for $r>0$, then by the decomposition in the proof of Lemma \ref{lem: coercive decomp}, applied to $(p, \lambda)$,
\begin{align*} 
        \OAE_\alpha(\xi_r, G) &= \alpha \Reg(G) + \tfrac{1}{2}\E^\lambda \norm{\wA(x) -\ubar{\wA}^{\mu_0}}_G^2  - \E^\lambda \langle \wA(x) -\ubar{\wA}^{\mu_0}, y\rangle_G \\
        & +\E^p \log \E^{\mu_0} \exp\left(-\frac{1}{2} \norm{\wA(x) -\ubar{\wA}^{\mu_0}}_G^2 - \langle \wA(x) -\ubar{\wA}^{\mu_0}, \ubar{\wA}^{\mu_0}- y+ \xi_0\rangle_G - r \langle \wA(x) - \ubar{\wA}^{\mu_0}, v \rangle_G  \right) \\
        &= \OAE_\alpha(\xi_0, G).
\end{align*}
\end{proof}

\begin{proposition}[Empirical objectives]\label{prop: empirical convergence}
Let $\tilde\lambda$ be a probability measure on $\R^n \times \R^m$ with finite support, $\tilde p$ its marginal on $\R^m$, and for each $y$ in the support of $\tilde p$ let $\tilde\mu_0^{(y)}$ be a probability measure on $\R^n$ with finite support. Consider
\begin{equation*}
    \tilde\OAE_\alpha(\xi, G) := \alpha \Reg(G) + \tfrac12 \E^{\tilde\lambda}\norm{\wA(x) - y + \xi}_G^2 + \E^{\tilde p} \log \E^{\tilde\mu_0^{(y)}} \exp\left(-\tfrac12 \norm{\wA(x') - y + \xi}_G^2\right).
\end{equation*}
Assume
\begin{itemize}
    \item[(E1)] for every $y \in \mathrm{supp}\,\tilde p$, no $v \neq 0$ exists such that $v \cdot \wA$ is constant on $\mathrm{supp}\,\tilde\mu_0^{(y)}$;
    \item[(E2)] $\ubar{\wA}^{\tilde\lambda} := \E^{\tilde\lambda} \wA(x)$ lies in the interior of $\mathrm{conv}\,\wA(\mathrm{supp}\,\tilde\mu_0^{(y)})$ for every $y \in \mathrm{supp}\,\tilde p$.
\end{itemize}
Then $\tilde\OAE_\alpha$ is bi-convex, Morse and coercive, and the proximal alternating minimisation sequence of Definition~\ref{def:proximal_alternating} converges to a critical point of $\tilde\OAE_\alpha$, which is a local minimum.
\end{proposition}
\begin{proof}
The proofs of Theorem~\ref{thm: biconvexity} and Lemma~\ref{lem: Morse property} use only that $\tilde p$ is the $y$-marginal of $\tilde\lambda$ (for the cancellation of the $y$-terms in \eqref{eq: derivative xi of OAE} and of the cross term in \eqref{eq: second derivative at xi critical}) and Lemma~\ref{lem: cov strictly psd} for the inner measure, which is (E1). Lemma~\ref{lem: coercive decomp} holds with $\mu_0$ replaced by $\tilde\mu_0^{(y)}$ inside the $\E^{\tilde p}$-expectation, and (E2) is equivalent to $C_*(\tilde\lambda) > 0$: $C_* = 0$ iff some half-space with $\ubar{\wA}^{\tilde\lambda}$ on its boundary contains $\wA(\mathrm{supp}\,\tilde\mu_0^{(y)})$. Coercivity and convergence then follow as in Proposition~\ref{prop:iterates_bounded}.
\end{proof}
\begin{remark}\label{rem: alternating optimization works for empirical}
For the single-loop estimator of Section~\ref{subsec:empirical}, where $\tilde\mu_0^{(y)} = \mu_0$ is Gaussian and $\wA$ is linear and surjective, (E1) is Assumption~\ref{ass: assumption conv of AM} and (E2) holds automatically, since $\wA_\# \mu_0$ has full support. 

For the double-loop estimator, where $\tilde\mu_0^{(y_j)}$ is the empirical measure of the inner samples $(x_{i,j})_{i=1}^{L}$, (E1) and (E2) are conditions on the sampled data: if $\wA(\mathrm{supp}\,\mu_0)$ affinely spans $\R^m$ and $\ubar{\wA}^{\mu_0}$ lies in the interior of $\mathrm{conv}\,\mathrm{supp}(\wA_\#\mu_0)$, both hold with probability tending to one as $L, M \to \infty$, but for fixed $L$ they may fail, in which case $\tilde\OAE_\alpha$ is constant along a line in $\xi$ as in Proposition~\ref{prop:iterates_bounded} and boundedness of the iterates cannot be guaranteed.
\end{remark}

\section{Concentration under empirical approximations}
\label{subsec:empirical}

In many applications, the full Bayesian joint distribution $\lambda$ is unavailable and, instead, some empirical or training data is given. This could be obtained, for instance, from physical measurements of the system or from a limited number of evaluations of a high-fidelity computational model.

Recall from \eqref{eq: objective function; easy} that
\begin{equation*}
    \OAE(\xi, G) = \tfrac12 \, \E^\lambda \norm{\wA(x)-y+\xi}_G^2 + \E^p \log \hat Z(y;\xi,G), \qquad \hat Z(y;\xi,G) = \E^{\mu_0} \exp\left(-\tfrac12\norm{\wA(x)-y+\xi}_G^2\right).
\end{equation*}
An obstruction to a direct sample-average approximation of $\OAE$ is the term $\E^p \log \hat Z(y;\xi,G)$. While the first term is an expectation with respect to $\lambda$ and hence immediately amenable to Monte Carlo estimation from the pairs $(x_j,y_j)$, the normalising constant $\hat Z(y;\xi,G)$ is itself an integral over $\mu_0$ and is generally intractable. We present now estimators for the Gaussian linear case, when the normalisation $\hat Z(y;\xi,G)$ can be computed explicitly. For the sake of simplicity, we will focus on the finite-dimensional case $\Space = \R^n$.

\begin{definition}[Single-loop estimators]
Given i.i.d.\ samples $(x_j,y_j)_{j=1}^M \sim \lambda$, we define the empirical regularised objective functional for a Gaussian prior $\mu_0 \sim \N(\bar x, C_0)$ on $\R^n$ with density $p_0$ and linear $\wA$ by
\begin{equation*}
    \widehat{\OAE}_\alpha^{M, \infty}(\xi, G) :=
     \alpha \Reg(G) + \tfrac12 \log \det  \widehat{C}(G) + \frac{1}{2 M} \sum_{j= 1}^M \norm{x_j-\hat{m}(y_j; \xi, G)}^2_{\widehat{C}(G)^{-1}} + \frac{1}{M} \sum_{j= 1}^M\log p_0(x_j) + \tfrac n2\log(2\pi),
\end{equation*}
which serves as an approximation of \eqref{eq: empirical expected KL divergence, linear surrogate}. Indeed, by the definition of $\hat q$ it holds point-wise in $(x, y) \in \R^n \times \R^m$,
\begin{eqnarray*}
    \tfrac{1}{2} \norm{\wA(x)-y +\xi}_G^2 + \log \hat Z(y;\xi, G) = \log p_0(x) - \log \hat q(x|y;\xi, G),
\end{eqnarray*}
\end{definition}
so that $\widehat{\OAE}_\alpha^{M, \infty}$ is of the form \eqref{eq: tilde OAE} for $(\tilde{p}_M, \tilde \lambda_M) = \frac{1}{M} \sum_{j=1}^M \delta_{y_j}, \frac{1}{M}\sum_{j=1}^M \delta_{x_j, y_j}$.
Above the term $ \tfrac{1}{M} \sum_{j= 1}^M\log p_0(x_j)$ does not depend on $\xi, G$ and can hence be omitted from the minimisation problem. We recall that posterior mean and covariance satisfy
\begin{align*}
    \begin{aligned}
\hat{m}(y;\xi, G) = \widehat{C}(G) \left( C_0^{-1} \bar x + \wA^\top G (y- \xi)\right), \quad 
\widehat{C}(G) = \left( C_0^{-1} + \wA^\top G \wA \right)^{-1}.
    \end{aligned}
\end{align*}

We recall that convergence of the alternating optimisation algorithm is guaranteed for empirical estimators as pointed out in Proposition \ref{prop: empirical convergence} and Remark \ref{rem: alternating optimization works for empirical}.
Our objective is to show that a sequence of optimisers of the empirical single-loop problems $\widehat{\OAE}_\alpha^{M, \infty}$ is consistent, that is, the estimators converge to an optimum of the limiting objective. 
To do so, we shall rely on the following well-known result, which establishes the desired result up to containment of the empirical estimators in a compact set with high probability. We point out that the original statement is formulated for maximisation rather than minimisation and hence the signs have been reversed compared to the original reference.

\begin{theorem}[{\cite[Theorem 5.14]{Vaart_1998}}]\label{thm: M estimation general thm}
    Denote by $\Theta$ a parameter space with metric $d$ and assume that $l: X\times \R^m \times \Theta \rightarrow \R$ is lower semi-continuous in $\theta \in \Theta$ for fixed $(x, y)\in X\times \R^m$ and satisfies for any sufficiently small ball $U \subset \Theta$
    \begin{eqnarray}\label{eq: bounded away from minus inf}
        \E^\lambda \left[\inf_{\theta \in U} l_\theta(x, y)  \right] > - \infty.
    \end{eqnarray}
    Assume that the set of minimisers $\Theta_0$ of $\theta \mapsto \E^\lambda l_\theta(x, y)$ is non-empty, and let $(\hat \theta_M)_{M \in \Natural}$ be estimators satisfying, for some $\theta_0 \in \Theta_0$,
    \begin{eqnarray}\label{eq: near minimiser}
        \frac{1}{M} \sum_{j=1}^M l_{\hat \theta_M}(x_j, y_j) \le \frac{1}{M} \sum_{j=1}^M l_{\theta_0}(x_j, y_j) + o_P(1).
    \end{eqnarray}
    Then for any $\varepsilon > 0 $ and compact set $K\subset \Theta$
    \begin{eqnarray*}
        \mathbb{P}\left(d(\hat \theta_M, \Theta_0) > \varepsilon \land \hat \theta_M \in K\right) \rightarrow 0.
    \end{eqnarray*}
\end{theorem}
Exact minimisers of $\theta \mapsto \frac1M \sum_{j=1}^M l_\theta(x_j, y_j)$ satisfy \eqref{eq: near minimiser} trivially; the weaker form is needed for the plug-in estimator of Corollary~\ref{cor: gaussian linear consistency} below.
To ensure the required containment of the empirical estimators in a compact set with high probability, we will use the following lemma.

\begin{lemma}\label{lem: emp opt in compact set}
    Let Assumptions \ref{ass: assumption} and \ref{ass: assumption conv of AM} hold. Assume that $\mu_0 = \N(\bar x, C_0)$ with density $p_0$ and $\wA: \R^n \to \R^m$. Let
    \begin{eqnarray}\label{eq: optimal emp est}
        (\hat \xi_M, \hat G_M) \in \argmin_{(\xi, G) \in \R^m \times \Sm} \widehat{\OAE}_\alpha^{M, \infty}(\xi, G).
    \end{eqnarray}
    Then for every $\kappa > 0$ there exists a \emph{deterministic} compact set $K_\kappa \subset \R^m \times \Sm$ such that
        \begin{eqnarray*}
            \mathbb{P}\Big( (\hat \xi_M, \hat G_M) \in K_\kappa \ \text{ for all } M \in \Natural \Big) \ \ge \ 1 - \kappa.
        \end{eqnarray*}
\end{lemma}

\begin{proof}
    Write $\ell_M := \tfrac{1}{M}\sum_{j=1}^M \log p_0(x_j)+ \tfrac{n}{2}\log(2\pi)$, which converges almost surely to $\E^{\mu_0}\log p_0(x) > -\infty$ by the law of large numbers and is therefore almost surely bounded in $M$. Moreover,
    \begin{equation*}
        \widehat{\OAE}_\alpha^{M, \infty}(0, I_m) - \ell_M \to  \alpha \Reg(I_m) + \tfrac12 \log \det  \widehat{C}(I_m) + \frac{1}{2 } \E^\lambda \norm{x-\hat{m}(y; 0, I_m)}^2_{\widehat{C}(I_m)^{-1}} < \infty
    \end{equation*}
    almost surely, so that by optimality
    \begin{equation*}
        B := \sup_{M \in \Natural} \Big( \widehat{\OAE}_\alpha^{M, \infty}(\hat \xi_M, \hat G_M) - \ell_M \Big) \le \sup_{M \in \Natural} \Big( \widehat{\OAE}_\alpha^{M, \infty}(0, I_m) - \ell_M \Big) < \infty \quad \text{almost surely}.
    \end{equation*}
    
    \textbf{Confinement of $G$:}
        Since the quadratic term is non-negative and $\log \det \widehat C(G) = - \log\det (C_0^{-1} + \wA^\top G \wA)$, it holds for every $(\xi, G) \in \R^m \times \Sm$
        \begin{eqnarray*}
            \widehat{\OAE}_\alpha^{M, \infty}(\xi, G) - \ell_M \ge \alpha \Reg(G) - \tfrac12 \log \det\left( C_0^{-1} + \wA^\top G \wA \right).
        \end{eqnarray*}
        Since $\wA^\top G \wA \preceq \norm{\wA}_F^2 \tr(G) I_n$, we have that 
        \begin{eqnarray*}
            \det \left(C_0^{-1}( I + C_0^{1/2}\wA^\top G \wA C_0^{1/2})\right) \le \det C_0^{-1} \left( 1+ \norm{C_0^{1/2}}_F^2\norm{\wA}_F^2 \tr(G)\right)^n
        \end{eqnarray*} 
        and hence, for all $M \in \Natural$,
        \begin{eqnarray*}
            B \ge \alpha \Reg(\hat G_M) + \log \det \left(C_0^{1/2} \right) - \tfrac{n}{2} \log\left(1+ \norm{C_0^{1/2}}_F^2\norm{\wA}_F^2 \tr(\hat G_M) \right).
        \end{eqnarray*}
    By Lemma~\ref{lem:regularisation}(ii), $\Reg(\hat G_M)$ grows at least linearly in $\gamma_{max}(\hat G_M)$ and tends to infinity as $\gamma_{min}(\hat G_M) \to 0$, whereas the last term grows only logarithmically in $\gamma_{max}(\hat G_M)$. Hence there exist almost surely finite random constants $0 < \gamma_{min} \le \gamma_{max} < \infty$ such that the eigenvalues of $\hat G_M$ lie in $[\gamma_{min}, \gamma_{max}]$ for all $M\in \Natural$.
    
    \textbf{Confinement of $\xi$:}
    For fixed $G \in \Sm$, the map $\xi \mapsto \widehat{\OAE}_\alpha^{M, \infty}(\xi, G)$ is a convex quadratic. Using $\partial_\xi \hat m(y; \xi, G) = - \widehat C(G) \wA^\top G$, its gradient is
    \begin{eqnarray*}
        \nabla_\xi \widehat{\OAE}_\alpha^{M, \infty}(\xi, G) = G \wA \, \frac{1}{M} \sum_{j=1}^M \big( x_j - \hat m(y_j; \xi, G) \big) = G \wA \Big( \bar x_M - \widehat C(G) \big( C_0^{-1} \bar x + \wA^\top G (\bar y_M - \xi) \big) \Big),
    \end{eqnarray*}
    where $\bar x_M := \tfrac1M \sum_{j=1}^M x_j$ and $\bar y_M := \tfrac1M \sum_{j=1}^M y_j$. By Assumption~\ref{ass: assumption conv of AM} and Remark~\ref{rem:linear_gaussian}, $\wA$ is surjective, so $\wA \widehat C(G) \wA^\top$ is positive definite and the minimiser in $\xi$ is unique and given by
    \begin{eqnarray}\label{eq: explicit xi hat M}
        \hat \xi_M = \bar y_M - \big( \wA \widehat C(\hat G_M) \wA^\top \hat G_M \big)^{-1} \wA \Big( \bar x_M - \widehat C(\hat G_M) C_0^{-1} \bar x \Big).
    \end{eqnarray}
    The sample means converge almost surely by the law of large numbers, and $G \mapsto (\wA \widehat C(G) \wA^\top G)^{-1}$ is continuous on the compact set $\{G \in \Sm : \gamma(G) \in [\gamma_{min}, \gamma_{max}]\}$. Hence $\sup_{M \in \Natural} \norm{\hat \xi_M}_2 < \infty$ almost surely.

\textbf{From almost-sure to uniform confinement.} Set
\begin{eqnarray*}
    \rho := \max\Big\{ \sup_{M \in \Natural} \gamma_{max}(\hat G_M), \ \sup_{M\in \Natural} \gamma_{min}(\hat G_M)^{-1}, \ \sup_{M\in \Natural} \norm{\hat \xi_M}_2 \Big\}.
\end{eqnarray*}
We have shown that $\rho < \infty$ almost surely. Since the events $\{\rho \le n\}$ increase to $\{\rho < \infty\}$, continuity of $\mathbb P$ from below gives $\mathbb P(\rho \le n) \to 1$ as $n \to \infty$ (See e.g. \cite[Theorem 1.8]{folland2013real}). Hence for every $\kappa > 0$ there is $\rho_\kappa < \infty$ with $\mathbb P(\rho > \rho_\kappa) < \kappa$. The set
\begin{eqnarray*}
    K_\kappa := \left\{ (\xi, G) \in \R^m \times \Sm \ : \ \norm{\xi}_2 \le \rho_\kappa, \ \ \rho_\kappa^{-1} I_m \preceq G \preceq \rho_\kappa I_m \right\}
\end{eqnarray*}
is a deterministic compact subset of $\R^m \times \Sm$, and $\{\rho \le \rho_\kappa\} \subset \{(\hat\xi_M, \hat G_M) \in K_\kappa \ \forall M\}$,
and hence
\begin{eqnarray*}
    \mathbb{P}\big( (\hat \xi_M, \hat G_M) \in K_\kappa \ \forall M \big) \ge \mathbb P(\rho \le \rho_\kappa) > 1-\kappa.
\end{eqnarray*}
\end{proof}

\begin{theorem}\label{thm: M estimation}
    Let $d$ be a metric on $\R^m\times \Sm$ inducing the standard topology, and let Assumptions \ref{ass: assumption} and \ref{ass: assumption conv of AM} hold. Assume that $\mu_0 = \N(\bar x, C_0)$ with density $p_0$ and $\wA: \R^n \to \R^m$. The set
    \begin{equation*}
        \Theta_0 := \argmin_{(\xi,G)\in\R^m\times\Sm} \OAE_\alpha(\xi,G)
    \end{equation*}
    is non-empty (Remark \ref{rem:existence_minimum}), and the estimator $(\hat\xi_M,\hat G_M)$ defined in \eqref{eq: optimal emp est} for the objective $\widehat{\OAE}_\alpha^{M, \infty}$ satisfies, for every $\varepsilon>0$,
    \begin{equation*}
        \mathbb{P}\left[ d\big((\hat\xi_M,\hat G_M), \Theta_0\big) > \varepsilon \right] \;\longrightarrow\; 0, \qquad M\to\infty.
    \end{equation*}
\end{theorem}

\begin{proof}
    We apply Theorem~\ref{thm: M estimation general thm} with $\theta = (\xi, G)$ and
    \begin{eqnarray*}
        l_\theta(x, y) := \alpha \Reg(G) + \tfrac12 \log \det \widehat C(G) + \tfrac12 \norm{x - \hat m(y; \xi, G)}^2_{\widehat C(G)^{-1}} + \log p_0(x),
    \end{eqnarray*}
    so that $\widehat{\OAE}_\alpha^{M, \infty}(\theta) = \frac1M \sum_{j=1}^M l_\theta(x_j, y_j)$ and $\E^\lambda l_\theta(x,y) = \OAE_\alpha(\theta)$ up to the additive constant of \eqref{eq: empirical expected KL divergence, linear surrogate}; hence $\Theta_0$ is the set of minimisers of $\E^\lambda l_\theta$. The map $\theta \mapsto l_\theta(x, y)$ is continuous, and since the quadratic term is non-negative, $l_\theta(x,y) \ge \alpha \Reg(G) + \tfrac12 \log\det \widehat C(G) + \log p_0(x)$, where the first two terms are continuous in $\theta$ and hence bounded on small balls, while $\E^{\mu_0} \log p_0(x) > -\infty$. This gives \eqref{eq: bounded away from minus inf}, and \eqref{eq: near minimiser} holds since $(\hat \xi_M, \hat G_M)$ is an exact minimiser. We notice that for any $\kappa > 0$ 
    \begin{eqnarray*}
         \lim_{M \to \infty} \mathbb{P}\left[ d\big((\hat\xi_M,\hat G_M), \Theta_0\big) > \varepsilon \right] &\le&  \lim_{M \to \infty} \mathbb{P}\left[ d\big((\hat\xi_M,\hat G_M), \Theta_0\big) > \varepsilon,  (\hat\xi_M,\hat G_M) \in K_\kappa \right] \\
         & & +  \lim_{M \to \infty} \mathbb{P}\left[ (\hat\xi_M,\hat G_M) \notin K_\kappa \right]
    \end{eqnarray*}
    By applying Theorem~\ref{thm: M estimation general thm}, the first term  equals zero, whereas by Lemma~\ref{lem: emp opt in compact set}, the second term is bounded by $\kappa>0$, which is arbitrary.
\end{proof}

The following corollary highlights a particular way to drop Assumption \ref{ass: assumption conv of AM} by cleverly defining estimators of $\xi$.
\begin{corollary}\label{cor: gaussian linear consistency}
   Let Assumption \ref{ass: assumption} hold. Assume that $\mu_0 = \N(\bar x, C_0)$ with density $p_0$ and $\wA: \R^n \to \R^m$, and introduce the estimators
    \begin{eqnarray*}
        \hat \xi_M= \frac{1}{M} \sum_{j=1}^M \left( y_j - \wA(x_j) \right), \quad 
        \hat G_M =  \argmin_{G \in \Sm}\widehat{\OAE}_\alpha^{M, \infty}(\hat \xi_M, G) 
    \end{eqnarray*}
    Then     \begin{equation*}
        \mathbb{P}\left[ d\big((\hat\xi_M,\hat G_M), \Theta_0\big) > \varepsilon \right] \;\longrightarrow\; 0, \qquad M\to\infty.
    \end{equation*}
\end{corollary}
\begin{proof}
    By the law of large numbers, $\hat \xi_M \to \hat \xi := \E^\lambda(y - \wA(x))$ almost surely, so $(\hat \xi_M)_{M \in \Natural}$ is almost surely bounded. By Theorem~\ref{lem: Gaussian case; xi}, $\hat \xi$ minimises $\OAE(\cdot, G)$ for every $G \in \Sm$, so $\min_{(\xi, G)} \OAE_\alpha = \min_{G} \OAE_\alpha(\hat \xi, G)$, and the latter minimum is attained at a unique $\hat G$ because $\OAE_\alpha(\hat \xi, \cdot)$ is strictly convex and coercive. In particular $(\hat \xi, \hat G) \in \Theta_0$. We verify \eqref{eq: near minimiser} with $\theta_0 = (\hat \xi, \hat G)$ and $l_\theta$ as in the proof of Theorem~\ref{thm: M estimation}. By definition of $\hat G_M$,
    \begin{eqnarray*}
        \widehat{\OAE}_\alpha^{M, \infty}(\hat \xi_M, \hat G_M) \le \widehat{\OAE}_\alpha^{M, \infty}(\hat \xi_M, \hat G),
    \end{eqnarray*}
    and $\xi \mapsto \widehat{\OAE}_\alpha^{M, \infty}(\xi, \hat G)$ is a quadratic whose coefficients are sample means of $x_j$, $y_j$ and their products, which converge almost surely; together with $\hat \xi_M \to \hat \xi$ this gives $\widehat{\OAE}_\alpha^{M, \infty}(\hat \xi_M, \hat G) - \widehat{\OAE}_\alpha^{M, \infty}(\hat \xi, \hat G) \to 0$ almost surely. Containment of $\hat G_M$ in a compact set with high probability follows as in Lemma~\ref{lem: emp opt in compact set}, using $(\hat\xi_M, I_m)$ as comparison point: $\widehat{\OAE}^{M,\infty}_\alpha(\hat\xi_M, I_m) - \ell_M$ is a quadratic in $\hat\xi_M$ with almost surely convergent coefficients, hence almost surely bounded in $M$, and the confinement argument for $G$ does not use Assumption~\ref{ass: assumption conv of AM}.
\end{proof}

More generally, if $\mu_0$ can be sampled and $\wA$ evaluated cheaply, $\hat Z$ can be approximated by Monte Carlo, resulting in a double-loop estimator for $\OAE$, as follows. The single-loop estimator exploits the Gaussian-linear structure, under which the normalising constant $\hat{Z}(y;\xi,G)$ is available in closed form. In many applications of interest, however, the prior $\mu_0$ is non-Gaussian, or the surrogate $\wA$ is nonlinear, making $\hat{Z}$ analytically intractable. Provided that $\mu_0$ can be sampled and $\wA$ evaluated efficiently, a natural remedy is to approximate $\hat{Z}(y;\xi,G)$ by a Monte Carlo average over draws from $\mu_0$. This leads to a double-loop estimator, in which the outer loop runs over the training pairs $(x_j,y_j)$ and the inner loop approximates the normalisation at each $y_j$.
\begin{definition}[Double-loop estimator]
    Given a coercive function $L: \Natural \to \Natural$ and i.i.d.\ samples $(x_j,y_j)_{j=1}^M \sim \lambda$ and $((x_{i, j})_{i=1}^{L(M)})_{j=1}^M \sim \mu_0$, we introduce the estimator. 
\begin{align*}
\begin{aligned}
    \widehat{\OAE}_\alpha^{M, L}(\xi, G) = \alpha \Reg(G) + \frac{1}{2M} \sum_{j=1}^M \norm{\wA(x_j)-y_j+\xi}_G^2 + \frac{1}{M} \sum_{j=1}^M \log \frac{1}{L(M)} \sum_{i=1}^{L(M)} \exp\left( -\tfrac12 \norm{\wA(x_{i,j})-y_j +\xi}_G^2 \right)
\end{aligned}
\end{align*}
\end{definition}
We recall the discussion in Proposition~\ref{prop: empirical convergence} and Remark~\ref{rem: alternating optimization works for empirical}, by which the previously developed theory is applicable under the data-dependent conditions.

For double-loop estimators we will make the following assumption of differentiability of $A, \wA$ with bounded derivatives.
\begin{assumption}\label{ass: assumption diff}
     Assume that $\Space = \R^n$ and there exists $L_3 >0$ such that for $A, \wA \in C^1(\R^n)$ it holds
     \begin{eqnarray*}
        \norm{\nabla \wA(x)} + \norm{\nabla A(x)} \le L_3 (1+\norm{x}).
     \end{eqnarray*}
\end{assumption}

Consistency of the double loop estimator is established in \cite{schneider2026convexapproximationframeworkneural} for very general function-valued likelihood approximations for functions $f \in \mathcal{F}$ contained in a function class $\mathcal{F}$ being a bounded subset of suitable weighted $C^1$-spaces of functions on $\mathbb{R}^n\times\mathbb{R}^m$. In the work, the authors analyse objectives of the form 
\begin{eqnarray*}
        \E^\lambda f(x, y) + \E^p \log \E^{\mu_0} \exp\left[- f(x, y) \right],
\end{eqnarray*}
which coincides with \eqref{eq: objective function; easy} for $f = f_{\xi, G}$ below.
In particular the work subsumes our approximations when taking 
\begin{eqnarray*}
    f_{\xi, G}(x, y) = \frac{1}{2}\norm{\wA(x) - y +\xi}_G^2, \quad \mathcal{F} =  \left\{  f_{\xi, G}, \, (\xi, G) \in K\right\} \subset C^2(\R^n \times \R^m),
\end{eqnarray*}
whenever $\xi$ and $G$ range over a compact subset $K$ of $\R^m\times\Sm$ and both the forward map $A$ and the approximating forward maps are uniformly bounded in $C^1$.

A metric measuring the discrepancy between potentials after accounting for their $y$-dependent normalisation, and hence identifies approximations that induce the same likelihood up to this normalisation is then given by
\begin{eqnarray*}
        d_\Phi\left((\xi_1, G_1), (\xi_2,G_2)\right) = \norm{  f_{\xi_1, G_1}  - f_{\xi_2, G_2} +\log \E^{\mu_0} \exp\left(-f_{\xi_1, G_1}\right)  - \log \E^{\mu_0} \exp\left(-f_{\xi_2, G_2}\right) }_{L^1(H \times \R^m, \lambda)}.
\end{eqnarray*}
We reiterate that $d_\Phi$ is in general not a norm on $\R^m \times \Sm$ but only a semi-metric on the space of potentials. 
We reiterate that convergence of $\xi_M$ in a Euclidean norm will need Assumption \ref{ass: assumption conv of AM} precisely because of the ambiguity pointed out in Proposition~\ref{prop:iterates_bounded}. The following lemma shows that under Assumption \ref{ass: assumption conv of AM}, indeed convergence with respect to $d_\Phi$ is equivalent to convergence in the Euclidean norm as long as the sequence is contained in a compact set in $\R^m \times \Sm$. 
Assumption \ref{ass: assumption conv of AM}, i.e., the requirement that $v^\top \wA$ is $\mu_0$-a.s. non-constant for every non-zero $v \in \R^m$, resolves the ambiguity of $G$ on the posterior that is pointed out in \eqref{eq: posterior const wrt G}.

\begin{lemma}\label{lem: equivalence of the norms}
Let Assumptions \ref{ass: assumption} hold and let $K$ be a compact set in $\R^m \times \Sm$ and $(\xi_k, G_k) \subset K$ be a sequence. Then the convergence
\begin{eqnarray*}
    (\xi_k, G_k) \to (\xi, G) \in K
\end{eqnarray*}
with respect to the Euclidean norm on $\R^m \times \Sm$ implies convergence with respect to $d_\Phi$.
If Assumption \ref{ass: assumption conv of AM} holds, then convergence with respect to $d_\Phi$ implies convergence in Euclidean norm.
\end{lemma}
\begin{proof}
    Assume $(\xi_k, G_k) \to (\xi, G) \in K$ with respect to the Euclidean norm. Then by \cite[Lemma 2.6]{schneider2026convexapproximationframeworkneural} it holds for some $p\ge p_0 >1$ and a constant~$C$
\begin{eqnarray*}
    d_\Phi\!\bigl((\xi_n,G_n),\,(\xi,G)\bigr) \le C \norm{\norm{\wA(x) - y +\xi_n}_{G_n}^2 - \norm{\wA(x) - y +{\xi}}_{G}^2}_{L_p(\lambda)}.
\end{eqnarray*}
The dominated convergence theorem applies by the arguments used in Theorem \ref{thm: Expected KL divergence analytic} and hence convergence with respect to $d_\Phi$ holds. 

Assume now that $(\xi_k, G_k) \to (\xi, G) \in K$ with respect to $d_\Phi$. Since $K$ is compact, there exists a convergent subsequence (with respect to Euclidean norm and hence with respect to $d_\Phi$), which we denote by $(\xi_{n_k}, G_{n_k}) \to (\xi^\dagger, G^\dagger) \in K$.
It holds 
\begin{eqnarray*}
    d_\Phi \bigl((\xi^\dagger, G^\dagger),\,(\xi,G)\bigr) \le d_\Phi \bigl((\xi^\dagger, G^\dagger),\,(\xi_{n_k},G_{n_k})\bigr) +d_\Phi \bigl((\xi_{n_k}, G_{n_k}),\,(\xi,G)\bigr) \to 0 ,
\end{eqnarray*}
Hence $d_\Phi \bigl((\xi^\dagger, G^\dagger), \, (\xi, G)\bigr) = 0$. By \cite[Lemma 2.4]{schneider2026convexapproximationframeworkneural} this implies that the term
\begin{eqnarray*}
    \norm{\wA(x) - y +\xi}_G^2 - \norm{\wA(x) - y +\xi^\dagger}_{G^\dagger}^2 
\end{eqnarray*}
 depends only on $y$ almost surely under $\lambda$. By expanding, 
 \begin{eqnarray}\label{eq: square term for equiv of norms}
     \wA(x)^\top (G-G^\dagger) \wA(x) + 2 \wA(x)^\top \left( G \xi - G^\dagger \xi^\dagger + (G^\dagger-G)y\right) + y^\top (G-G^\dagger) y
 \end{eqnarray}
depends only on $y$ almost surely under $\lambda$. Since $p$ admits a density and is thus equivalent to Lebesgue measure on $\R^m$ , evaluating \eqref{eq: square term for equiv of norms} at $y$ and $y+z$ for arbitrary $z\in\R^m$ and subtracting gives that $\wA(x)^\top(G^\dagger-G)z$ is $\mu_0$-a.s.\ constant in $x$. Since $z$ is 
arbitrary, Assumption~\ref{ass: assumption conv of AM} implies $G = G^\dagger$. Substituting back into \eqref{eq: square term for equiv of norms} gives that $\wA(x)^\top G(\xi-\xi^\dagger)$ is $\mu_0$-a.s.\ constant in $x$, and Assumption~\ref{ass: assumption conv of AM} together with invertibility of $G$ yields $\xi = \xi^\dagger$. 
Since every Euclidean-convergent subsequence of $(\xi_{n}, G_{n})\subset K$ must have limit $(\xi, G)$, compactness of $K$ implies $(\xi^\dagger, G^\dagger) = (\xi, G)$ and hence $(\xi_{n}, G_{n}) \to (\xi, G)$ in Euclidean norm.
\end{proof}

Following \cite{schneider2026convexapproximationframeworkneural}, we will use the following theorem to show consistency of the double-loop estimators.
\begin{theorem}[Consistency, {\cite[Theorem 5.7]{Vaart_1998}}]\label{thm: consistency, Vaart}
    Let $(\D, d)$ be a pseudo-metric space and $\phi: \D \to \R$ with non-empty set of minimisers $\D_0 := \argmin_{f \in \D} \phi(f)$. Assume that for every $\varepsilon > 0$
    \begin{align}
        \inf_{f \in \D : \, d(f, \D_0) > \varepsilon}  \phi (f) >& \inf_{f \in \D} \phi(f) \label{eq: uniqueness condition}
    \end{align}
    and 
    \begin{align}
                \sup_{f \in \D}  \big| \phi_N(f) -& \phi(f) \big|  \rightarrow 0 \label{eq: glivenko cantelli condition} 
    \end{align}
    in probability. Then any $\hat f_N$ satisfying $\phi_N(\hat f_N) \le \inf_{f\in \D} \phi_N(f) + o_P(1)$ satisfies $d(\hat f_N, \D_0) \to 0$ in probability.
\end{theorem}
\begin{proof}
    For $f_0 \in \D_0$, \eqref{eq: glivenko cantelli condition} gives $\phi(\hat f_N) \le \phi_N(\hat f_N) + o_P(1) \le \phi_N(f_0) + o_P(1) \le \phi(f_0) + o_P(1)$, so $\phi(\hat f_N) - \inf_\D \phi = o_P(1)$, and \eqref{eq: uniqueness condition} yields the claim.
\end{proof}

\begin{theorem}\label{cor: consistency for double loop}
    Let Assumptions \ref{ass: assumption} and \ref{ass: assumption diff} hold, let $K$ be a compact subset of $\R^m \times \Sm$ and set $\Theta_0 := \argmin_{(\xi, G) \in K} \OAE_\alpha(\xi, G)$, which is non-empty by continuity of $\OAE_\alpha$. For the estimators
    \begin{eqnarray*}
        (\hat \xi_M, \hat G_M) \in \argmin_{(\xi, G) \in K} \widehat{\OAE}_\alpha^{M, L}(\xi, G),
    \end{eqnarray*}
    it holds for every $\varepsilon > 0$ that
    \begin{eqnarray*}
        \mathbb{P}\left( d_\Phi\big((\hat \xi_M, \hat G_M), \Theta_0\big) > \varepsilon \right) \rightarrow 0  \quad \text{as } M \rightarrow \infty.
    \end{eqnarray*}
    If Assumption \ref{ass: assumption conv of AM} holds, the same holds with $d_\Phi$ replaced by the Euclidean distance.
\end{theorem}
\begin{proof}[Proof of Theorem \ref{cor: consistency for double loop}]
By \cite[Theorem 3.6]{schneider2026convexapproximationframeworkneural} applied with 
    \begin{eqnarray*}
        \mathcal{F}_1:= \mathrm{conv} \left(\left\{ \tfrac{1}{2} \norm{\wA(x) - y + \xi}_G^2; (\xi, G) \in K \right\} \cup \left\{ \tfrac{1}{2 \delta^2} \norm{A(x) - y}^2 \right\} \right),
    \end{eqnarray*}
    and since the regularisation term cancels in the difference, we get
    \begin{eqnarray*}
               \sup_{(\xi, G) \in K}  \big|  \widehat{\OAE}_\alpha^{M, L}(\xi, G) - \OAE_\alpha(\xi, G) \big| \le \sup_{f \in \F_1}  \big| \phi_N(f) - \phi(f) \big|  \rightarrow 0 
    \end{eqnarray*}
    in probability, which is \eqref{eq: glivenko cantelli condition} for $\D = K$, $\phi = \OAE_\alpha$ and $\phi_N = \widehat{\OAE}_\alpha^{M, L}$. It remains to verify \eqref{eq: uniqueness condition} with $d = d_\Phi$ and $\D_0 = \Theta_0$, i.e.\ that for every $\varepsilon > 0$
    \begin{align*}
       \inf_{(\xi, G) \in K, \, d_\Phi\left((\xi, G), \Theta_0\right) > \varepsilon}  \OAE_\alpha(\xi, G)   > c^* := \min_{(\xi, G) \in K} \OAE_\alpha(\xi, G).
    \end{align*}
Suppose, for contradiction, that for some $\varepsilon>0$ there exists a sequence $(\xi_n,G_n)\subset K$ with
\begin{equation}\label{eq:contradiction_assumption}
    d_\Phi\!\bigl((\xi_n,G_n),\,\Theta_0\bigr) \;\ge\; \varepsilon
    \qquad\text{for all } n\in\mathbb{N},
\end{equation}
and $\OAE_\alpha(\xi_n,G_n)\to c^*$. By compactness of $K$, after passing to a subsequence, $(\xi_n,G_n) \to (\xi^*,G^*) \in K$ in Euclidean norm. Since $\OAE_\alpha$ is continuous (Theorem~\ref{thm: Expected KL divergence analytic} and Lemma~\ref{lem:regularisation}), $\OAE_\alpha(\xi^*, G^*) = c^*$, so $(\xi^*,G^*)\in\Theta_0$. By Lemma~\ref{lem: equivalence of the norms}, Euclidean convergence implies $d_\Phi((\xi_n, G_n), (\xi^*, G^*)) \to 0$, hence $d_\Phi((\xi_n,G_n),\Theta_0) \to 0$ along the subsequence, contradicting \eqref{eq:contradiction_assumption}. Theorem~\ref{thm: consistency, Vaart} now yields the first claim.

Under Assumption~\ref{ass: assumption conv of AM}, for every $\varepsilon > 0$ there is $\eta > 0$ such that $d_\Phi(\theta, \Theta_0) \le \eta$ implies $\dist(\theta, \Theta_0) \le \varepsilon$ for all $\theta \in K$: otherwise there would be $\theta_n \in K$ and $\vartheta_n \in \Theta_0$ with $d_\Phi(\theta_n, \vartheta_n) \to 0$ and $\dist(\theta_n, \Theta_0) > \varepsilon$; by compactness we may assume $\theta_n \to \theta^*$ and $\vartheta_n \to \vartheta^* \in \Theta_0$ (closed) in Euclidean norm, so $d_\Phi(\theta^*, \vartheta^*) = 0$ by Lemma~\ref{lem: equivalence of the norms} and the triangle inequality, hence $\theta^* = \vartheta^*$ by the identification argument in the proof of that lemma, contradicting $\dist(\theta^*, \Theta_0) \ge \varepsilon$. The second claim follows.
\end{proof}

\begin{remark}
The compactness of $K$ in Theorem~\ref{cor: consistency for double loop} cannot easily be dispensed with by the argument used for the single-loop estimator in Lemma~\ref{lem: emp opt in compact set}. In the single-loop case, the Gaussian-linear structure makes the normalising constant $\hat{Z}$ available in closed form with respect to the same empirical measure $\tilde{\lambda}_M$, so the decomposition of Lemma~\ref{lem: decompositon of obj fct} applies and confinement of $G$ can be established independently of $\xi$. 
For the double-loop estimator, we might want to use the identity~\eqref{eq: from KL divergence to our objective}. However, the inner samples $(x_{i,j})$ approximating $\hat{Z}(y_j;\xi,G)$ are drawn from $\mu_0$, while the outer samples $(x_j,y_j)$ are drawn from $\lambda$. 
These two empirical measures do not share a common marginal, so~\eqref{eq: from KL divergence to our objective} does not hold for the double-loop objective, and the confinement argument for $G$ cannot be decoupled from $\xi$. 
\end{remark}

\section{Numerical experiments}
\label{sec: numerics}

\subsection{Computed tomography}

\begin{figure}[h]
\centering

\begin{minipage}{0.25\textwidth}
    \centering
    \includegraphics[width=\linewidth]{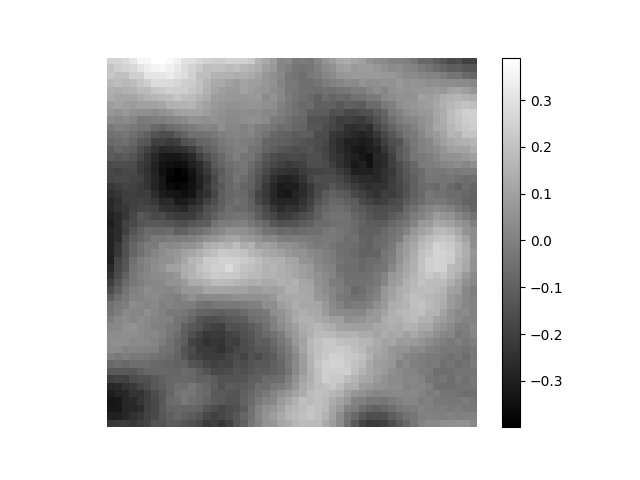}
    \caption{Ground truth $x_0$ for the CT problem, a sample from the prior $\mu_0$.}
    \label{fig:ground_truth_CT}
\end{minipage}
\hfill
\begin{minipage}{0.7\textwidth}
    \centering

    \includegraphics[width=0.3\linewidth]{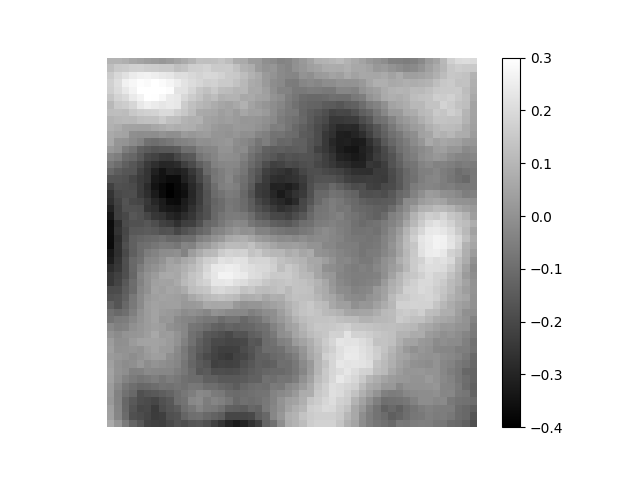}%
    \includegraphics[width=0.3\linewidth]{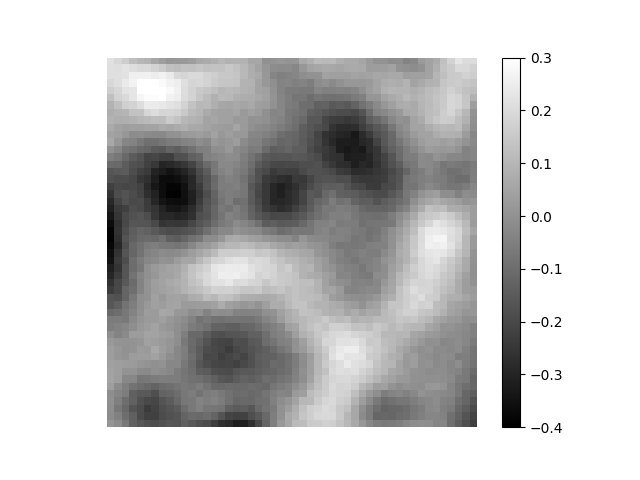}%
    \includegraphics[width=0.3\linewidth]{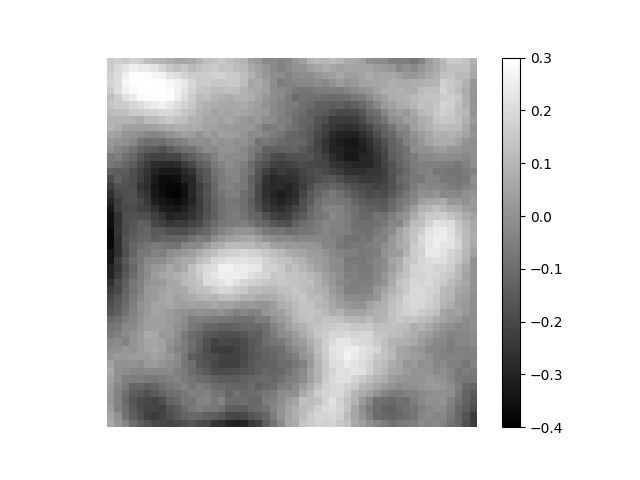}\\

    \includegraphics[width=0.3\linewidth]{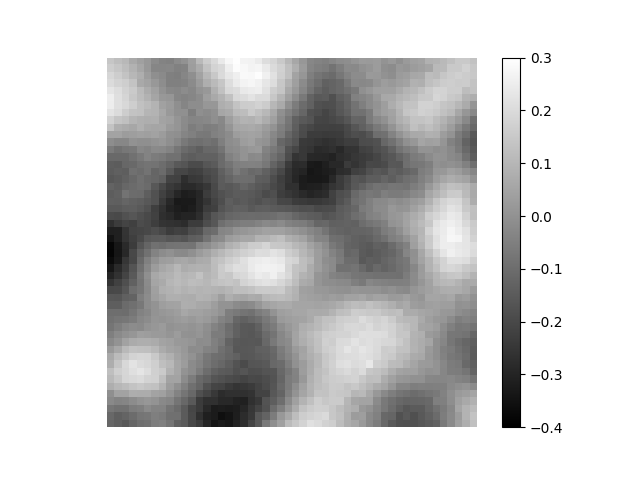}%
    \includegraphics[width=0.3\linewidth]{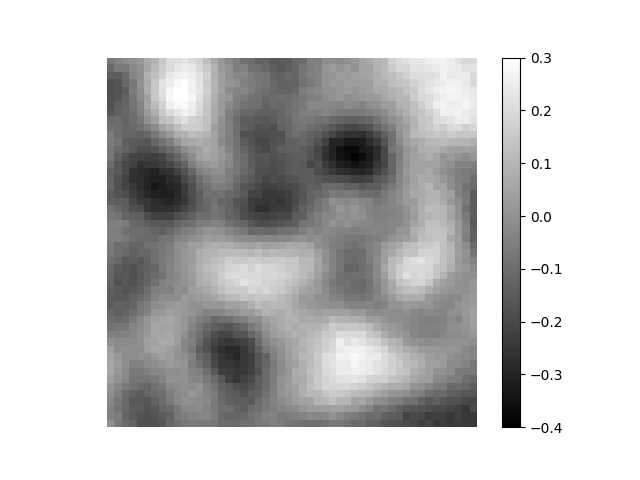}%
    \includegraphics[width=0.3\linewidth]{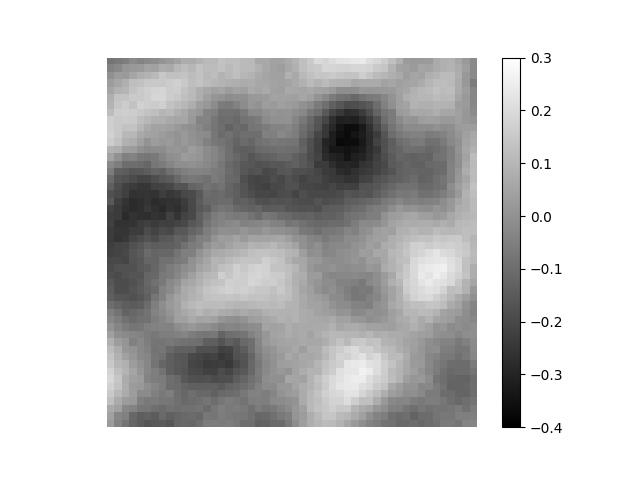}\\

    \includegraphics[width=0.3\linewidth]{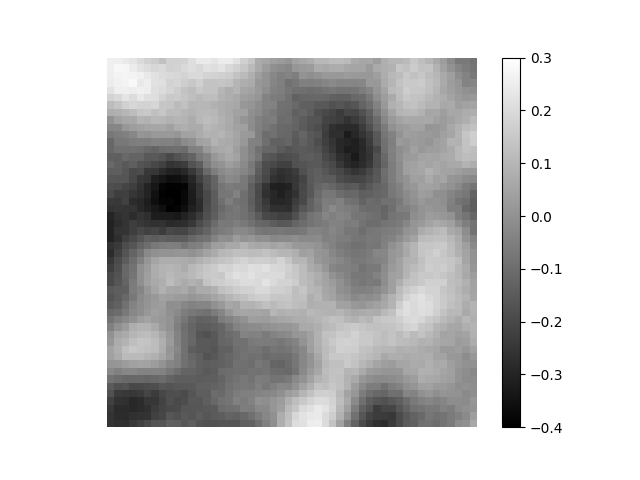}%
    \includegraphics[width=0.3\linewidth]{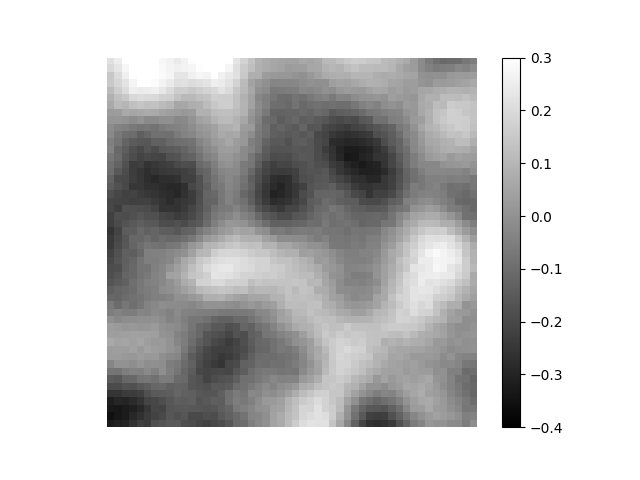}%
    \includegraphics[width=0.3\linewidth]{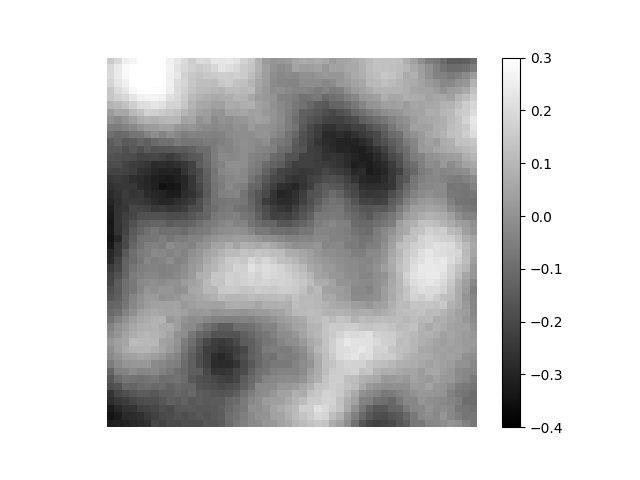}\\

    \includegraphics[width=0.3\linewidth]{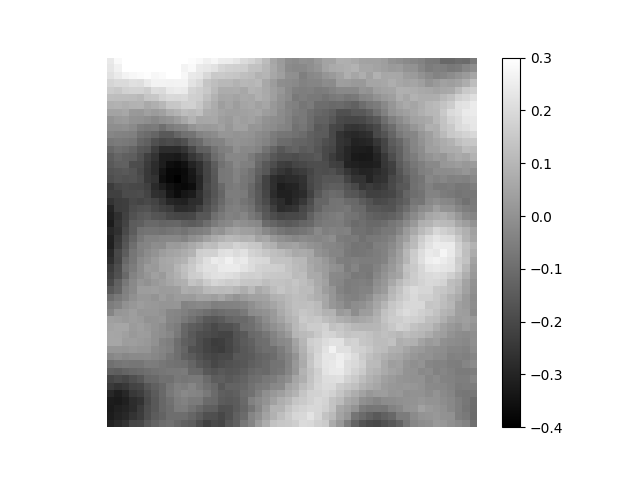}%
    \includegraphics[width=0.3\linewidth]{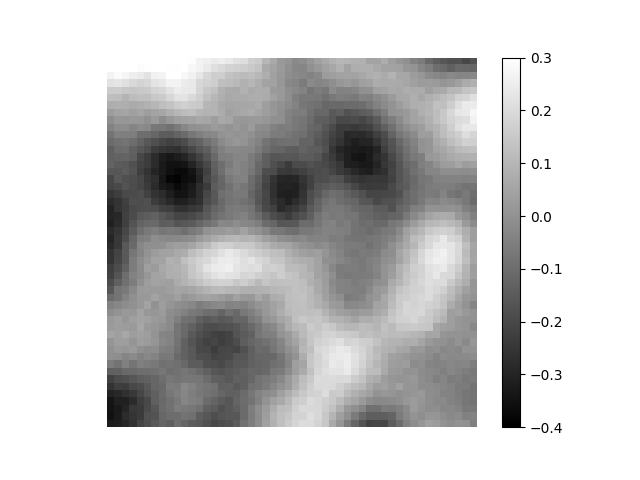}%
    \includegraphics[width=0.3\linewidth]{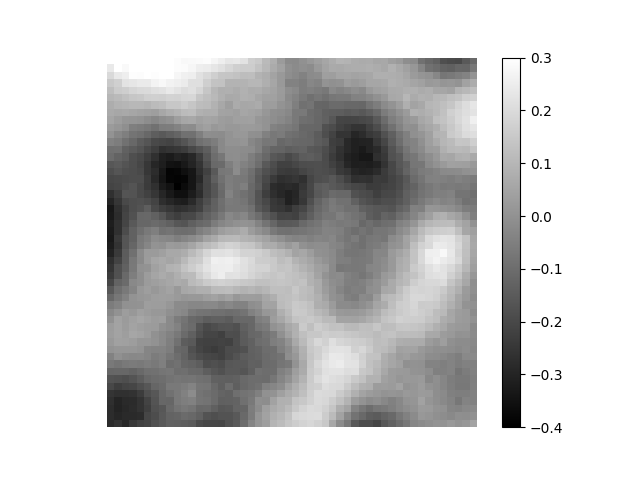}

    \caption{Posterior samples for the CT problem. Rows correspond, from top to bottom, to None, EEA, OAE and the true posterior.}
    \label{fig:CT_samples}
\end{minipage}

\end{figure}

We model the forward map $A$ as a discrete Radon transform of a $50 \times 50$ pixel image, so that $n = 2500$, with $47$ equispaced projection angles over $0$--$170$ degrees. The detector is assumed to have $50$ apertures spanning the width of the image, i.e.\ $50$ parallel line integrals are measured per angle and $m = 47 \cdot 50 = 2350$. The measurement is corrupted with additive Gaussian noise of level $\delta = 0.005$. The prior is Gaussian, $\mu_0 = \N(0, C_0)$, with $C_0$ corresponding to squared exponential covariance 
\begin{eqnarray*}
    k(x, y) = \frac{1}{\sigma^2} \exp\left(-\frac{\norm{x-y}^2_2}{2 \, l^2} \right),
\end{eqnarray*}
with length-scale $l = 0.1 (50)^{-1/2}$  and variance $\sigma^2 = 50$. A sample is shown in Figure~\ref{fig:ground_truth_CT}.

The surrogate $\wA$ is obtained from a misspecified geometry: the same number of equispaced projections is taken over $5$--$177$ degrees instead of $0$--$170$ degrees, so that the projection angles of $\wA$ are shifted and slightly stretched relative to those of $A$. Note that $\wA$ is surjective, so Assumption~\ref{ass: assumption conv of AM} holds by Remark~\ref{rem:linear_gaussian}.

We compare three choices of the Gaussian likelihood: the surrogate likelihood without any correction, i.e.\ $\xi = 0$ and $G = \delta^{-2} I_m$ (None); the enhanced error approach \eqref{eq: enhanced error approach} (EEA); and the optimal approximation error approach \eqref{eq:AEM_optimization} (OAE). Since the setting is Gaussian-linear, OAE uses closed-form objective function given in Section \ref{sec: Gaussian linear} and $\alpha = 0$, minimised over $G$ by the plain gradient descent, initialised at the EEA precision $G_*$ and stopped after $15,000$ iterations. All posterior quantities are Gaussian and available in closed form via \eqref{eq: fin dim lin post}. The expected KL divergence is evaluated in closed form..

\begin{figure}[h]
\centering
\begin{minipage}{0.23\textwidth}
  \centering
\includegraphics[width=1.15\linewidth]{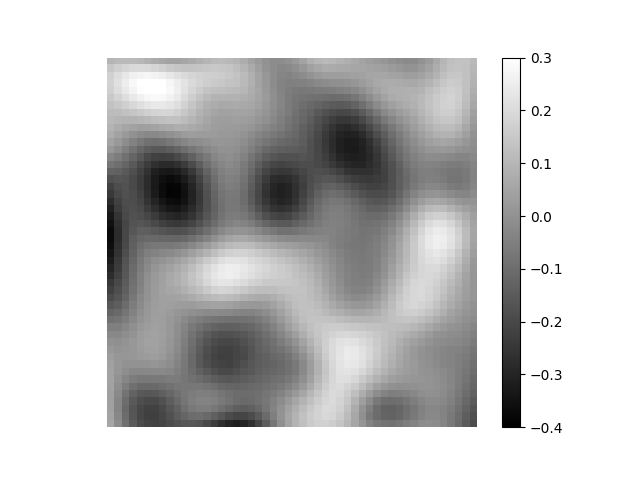}
\end{minipage}%
\begin{minipage}{0.23\textwidth}
  \centering
  \includegraphics[width=1.15\linewidth]{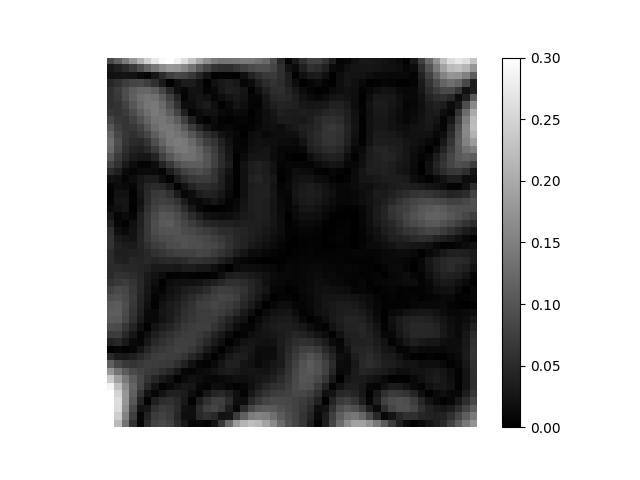}
\end{minipage}%
\begin{minipage}{0.23\textwidth}
  \centering
\includegraphics[width=1.15\linewidth]{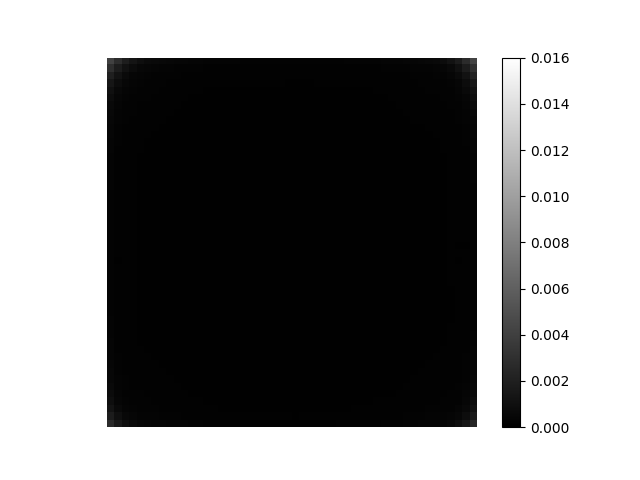}
\end{minipage}
\begin{minipage}{0.23\textwidth}
  \centering
\includegraphics[width=1.15\linewidth]{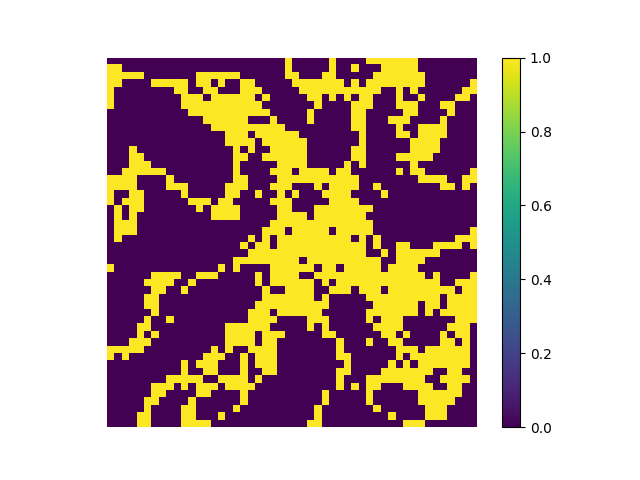}
\end{minipage} \\
\begin{minipage}{0.23\textwidth}
  \centering
\includegraphics[width=1.15\linewidth]{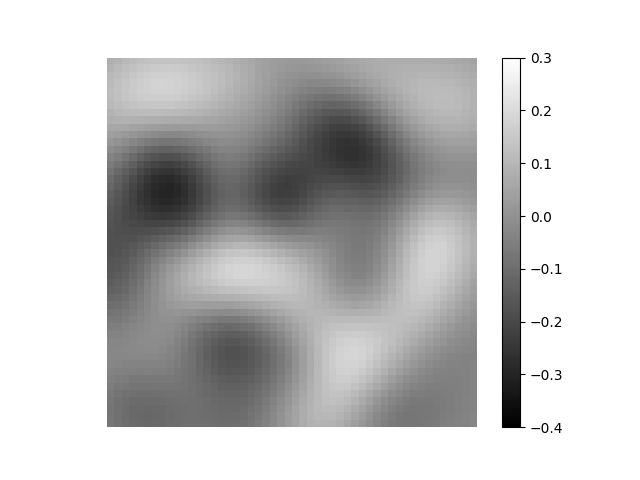}
\end{minipage}%
\begin{minipage}{0.23\textwidth}
  \centering
  \includegraphics[width=1.15\linewidth]{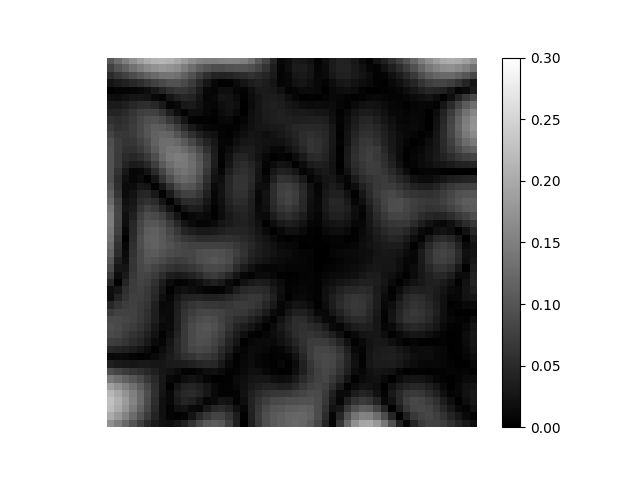}
\end{minipage}%
\begin{minipage}{0.23\textwidth}
  \centering
\includegraphics[width=1.15\linewidth]{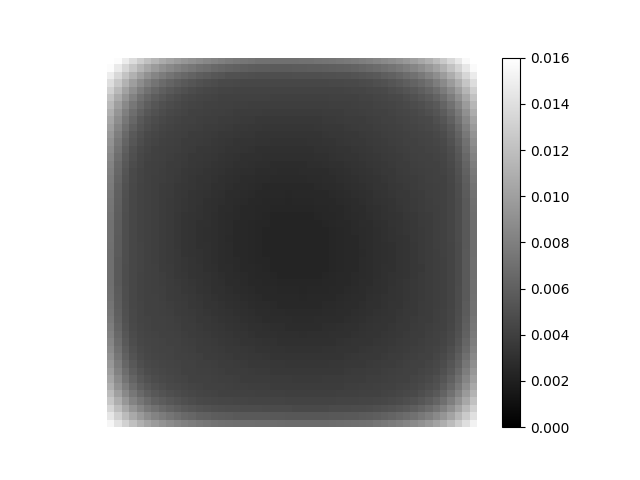}
\end{minipage}
\begin{minipage}{0.23\textwidth}
  \centering
  \includegraphics[width=1.15\linewidth]{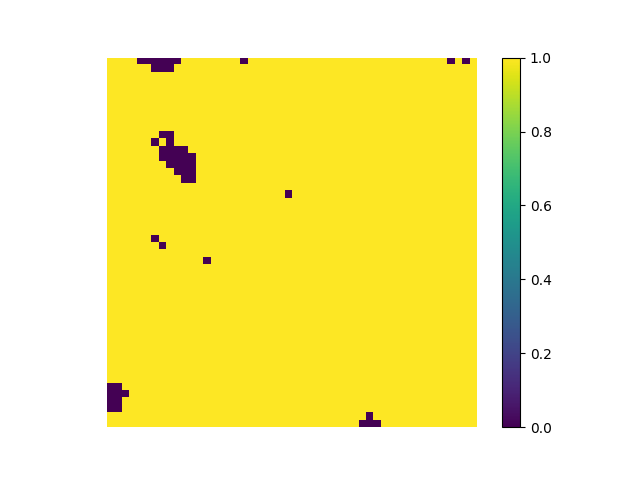}
\end{minipage} \\
\begin{minipage}{0.23\textwidth}
  \centering
  \includegraphics[width=1.15\linewidth]{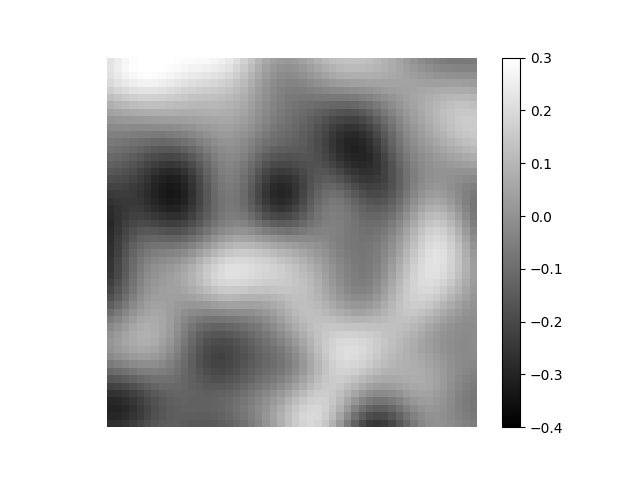}
\end{minipage}%
\begin{minipage}{0.23\textwidth}
  \centering
  \includegraphics[width=1.15\linewidth]{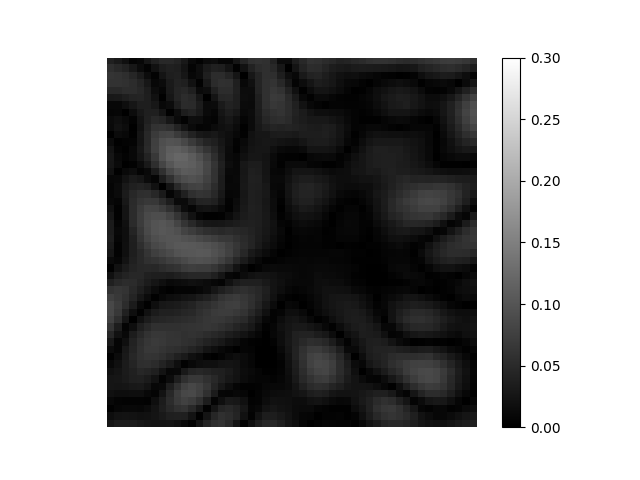}
\end{minipage}%
\begin{minipage}{0.23\textwidth}
  \centering
\includegraphics[width=1.15\linewidth]{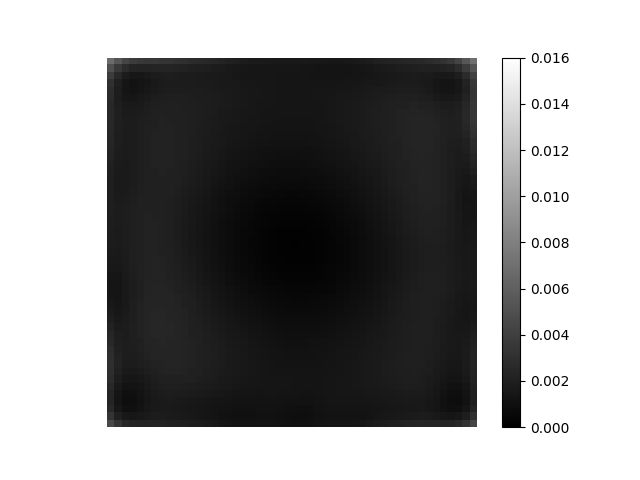}
\end{minipage}
\begin{minipage}{0.23\textwidth}
  \centering
  \includegraphics[width=1.15\linewidth]{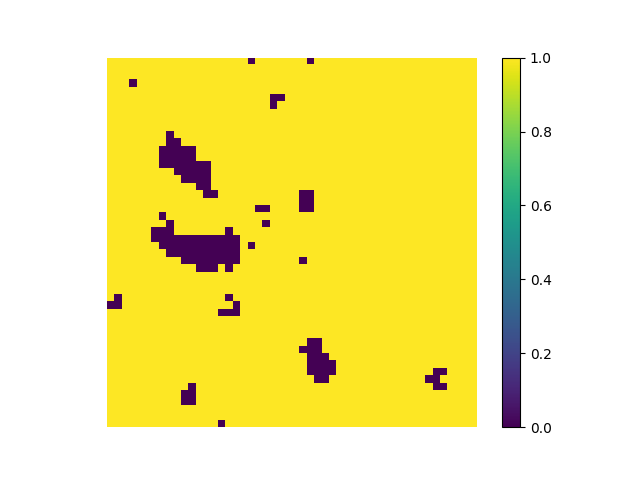}
\end{minipage}
\begin{minipage}{0.23\textwidth}
  \centering
  \includegraphics[width=1.15\linewidth]{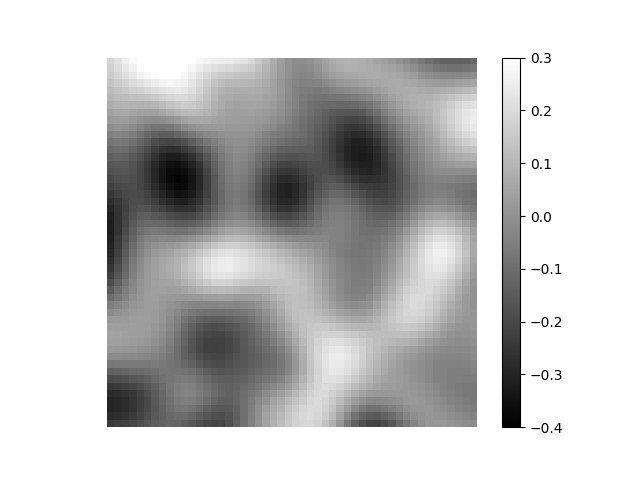}
\end{minipage}%
\begin{minipage}{0.23\textwidth}
  \centering
  \includegraphics[width=1.15\linewidth]{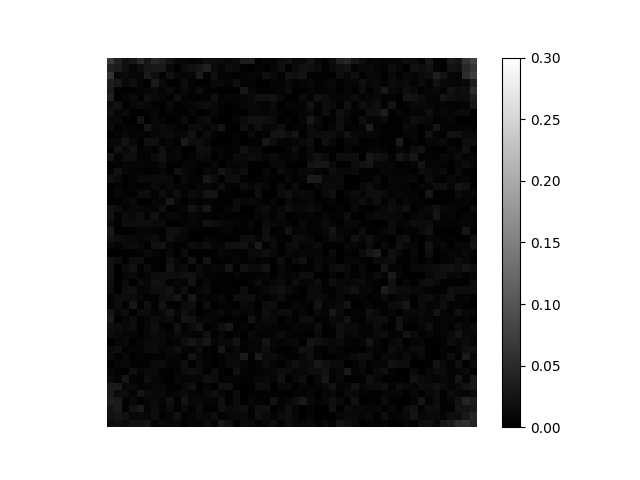}
\end{minipage}%
\begin{minipage}{0.23\textwidth}
  \centering
\includegraphics[width=1.15\linewidth]{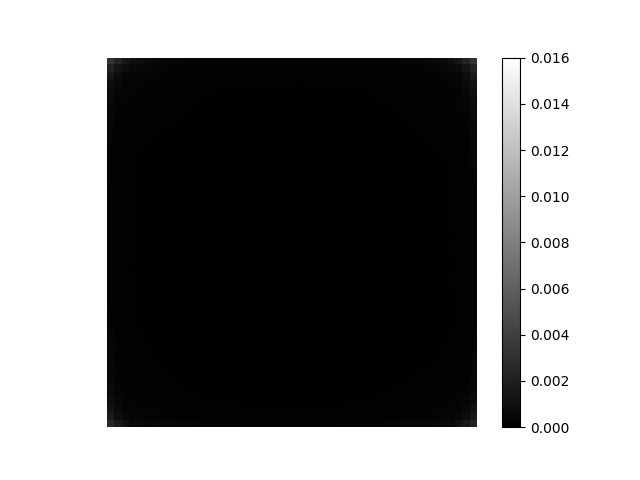}
\end{minipage}
\begin{minipage}{0.23\textwidth}
  \centering
  \includegraphics[width=1.15\linewidth]{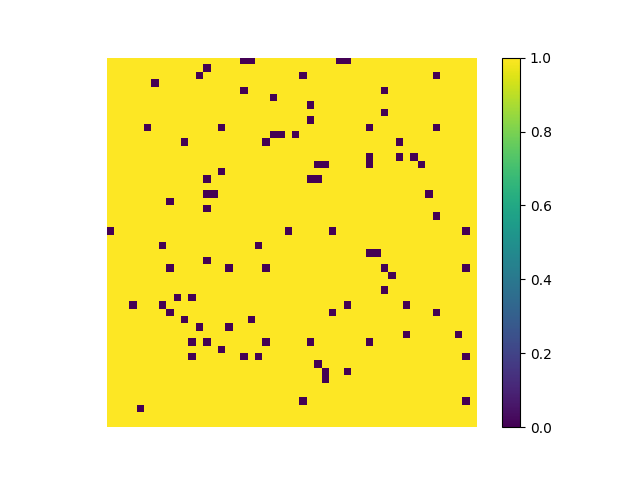}
\end{minipage}
\caption{For all methods, we report the posterior mean, difference posterior mean to the ground truth, diagonal of the posterior covariance, and the pixels contained in a confidence interval based on the posterior mean and variance. Top row corresponds to None, second row to EEA, third row to OAE and final row to the true posterior}
\label{fig:CT grf}
\end{figure}

Figure~\ref{fig:CT_samples} shows posterior samples for the three methods and for the true posterior, Figure~\ref{fig:CT grf} shows the posterior mean, its error with respect to the ground truth of Figure~\ref{fig:ground_truth_CT}, the pointwise posterior standard deviation and the pixels covered by a $2\sigma$ credible band, and Table~\ref{tab: CT} collects summary metrics.

Table~\ref{tab: CT} shows that OAE reduces the expected posterior KL divergence from $342$ (EEA) to $239$, and the error of the posterior mean by a factor of about $2.5$, both in expectation and for the single draw $x_0$. The uncorrected likelihood performs far worse on both counts, and its credible band covers only $37\%$ of the pixels, reflecting its overconfident covariance. The one metric on which None is best, $\norm{C - \widehat C}_F$, is expected: None retains the noise precision $\delta^{-2} I_m$, so its posterior covariance is closest to that of the true model, but this is of little use when the mean is biased. OAE inflates the covariance slightly more than EEA, yet its coverage ($94\%$) is as close to that of the true posterior ($96\%$) as EEA ($98\%$). 

\begin{table}[h]
    \centering
        \caption{Error metrics for the CT-imaging problem. All norms are averaged over the image pixels. In the quantities in the two last rows, $x_0$ corresponds to the unknown in Figure \ref{fig:ground_truth_CT} and $y_0$ to a simulated corresponding measurement.}
    \label{tab: CT}
    \begin{tabular}{c|c c c c}
         \hline
          & \textbf{None} & \textbf{EEA} & \textbf{OAE} & \textbf{Posterior} \\
         \hline
         $\norm{C - \widehat{C}}_F$ & \textbf{2.41e-09} & 2.25e-07  & 8.24e-07 & 0 \\ 
         $\E^{p} \norm{m(y) - \hat{m}(y)}^2$   & 0.0044 & 0.0037  & \textbf{0.0014} & 0 \\ 
         $\E^{p} \, \KL{\mu(\cdot | y)}{\hat \mu(\cdot | y)}$ & 13734 & 342  &  \textbf{239} & 0 \\ 
         $\E^{\lambda} \norm{x - \hat{m}(y)}^2$ & 0.0046 & 0.0039 & \textbf{0.0015} & 0.0002  \\
         \hline
         $ \norm{x_0 - \hat{m}(y_0)}^2$ & 0.0042 & 0.0038 & 0.0016 & 0.0002 \\
         $\#  \left\{ |x_0-\hat{m}(y_0)| \le 2 \, \sqrt{\text{diag} \widehat{C}} \right\}$ &37 \% & 98 \% & 94 \% & 96 \%
    \end{tabular}
\end{table}

\subsection{Darcy Flow}

We seek to recover the permeability $\sigma \in C^1 (\Omega)$ on the domain $\Omega = [0, 1]^2$ from the groundwater filtration $u \in C^2 (\Omega)$ for a fixed source $g \in C (\Omega)$.
The forward map $F: \sigma \mapsto u$ corresponds to solving the PDE
\begin{eqnarray*}
    - \nabla \cdot (\sigma \nabla u) &=& g \qquad  z \in \Omega \\
    u & = & 0 \qquad  z \in \partial \Omega.
\end{eqnarray*}

The surrogate $\wA$ is the affine map obtained by a numerical first-order Taylor expansion of $F$ about the prior mean. Training data consist of $M = 5000$ pairs $(x_j, y_j) \sim \lambda$. We use the double-loop estimator $\widehat{\OAE}_\alpha^{M, L}$ with $\alpha = 0.001$, where the inner expectation over $\mu_0$ in $\hat Z(y_j; \xi, G)$ is approximated $L= 100$ data points. The objective is minimised by the Adam optimizer with learning rate $0.01$, initialised at the optimal values for EEA.

\begin{figure}[h]    
        \centering
        \includegraphics[width=0.3\linewidth]{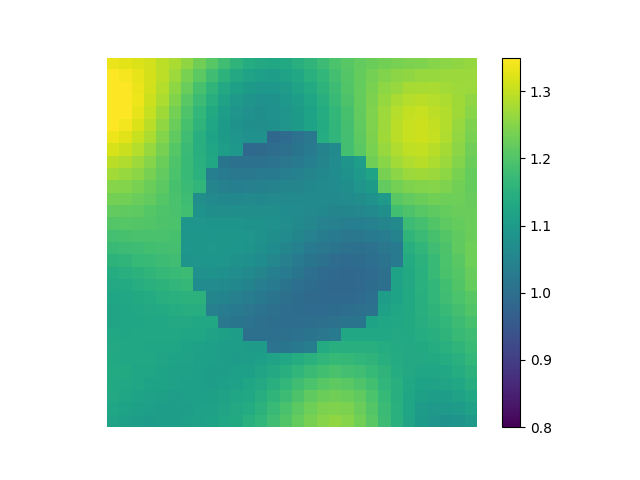}
    \caption{Ground truth sample from the prior}
    \label{fig:2d darcy; prior sample}

\begin{minipage}{0.2\textwidth}
  \centering
\includegraphics[width=1.\linewidth]{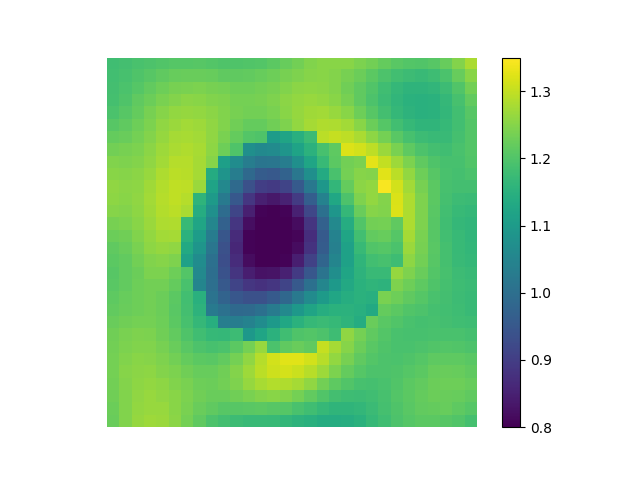}
\end{minipage}%
\begin{minipage}{0.2\textwidth}
  \centering
\includegraphics[width=1.\linewidth]{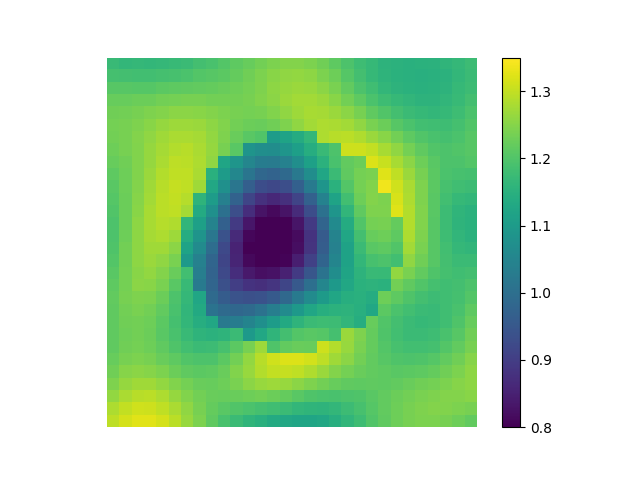}
\end{minipage}%
\begin{minipage}{0.2\textwidth}
  \centering
\includegraphics[width=1.\linewidth]{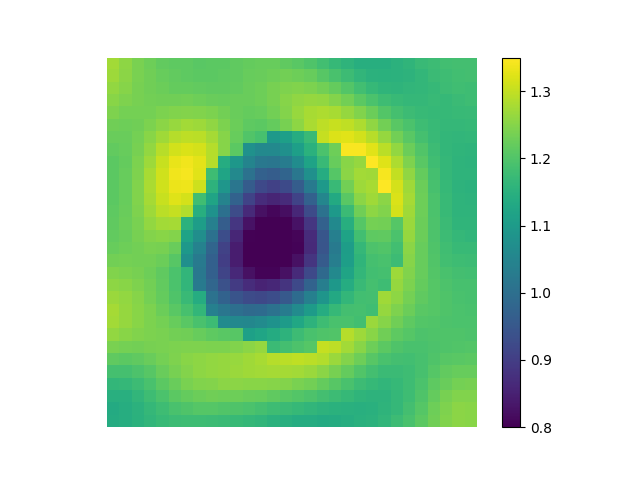}
\end{minipage}%
\begin{minipage}{0.2\textwidth}
  \centering
\includegraphics[width=1.\linewidth]{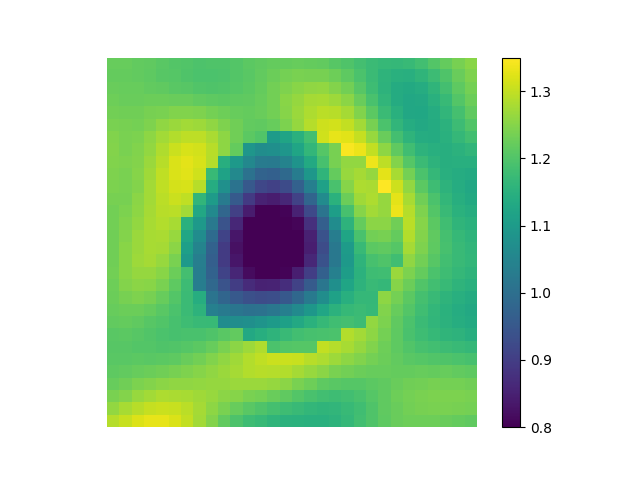}
\end{minipage}%
\begin{minipage}{0.2\textwidth}
  \centering
\includegraphics[width=1.\linewidth]{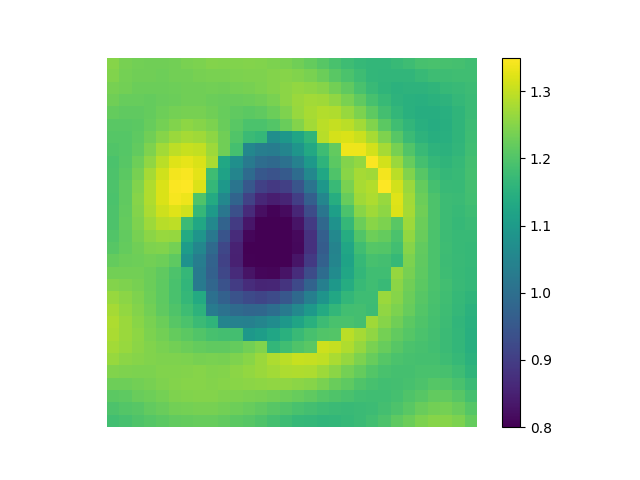}
\end{minipage} \\
\begin{minipage}{0.2\textwidth}
  \centering
\includegraphics[width=1.\linewidth]{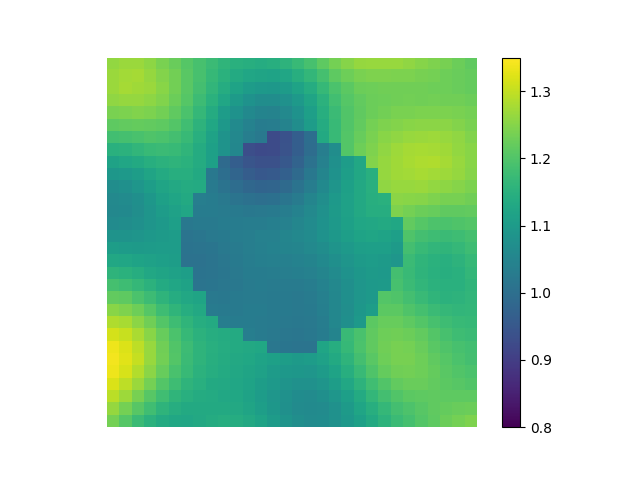}
\end{minipage}%
\begin{minipage}{0.2\textwidth}
  \centering
\includegraphics[width=1.\linewidth]{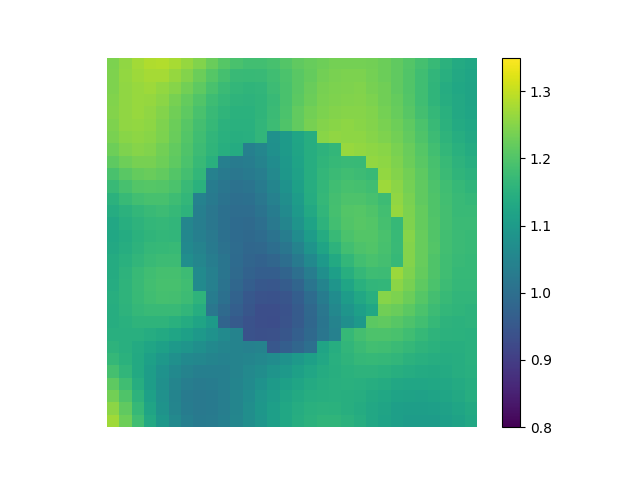}
\end{minipage}%
\begin{minipage}{0.2\textwidth}
  \centering
\includegraphics[width=1.\linewidth]{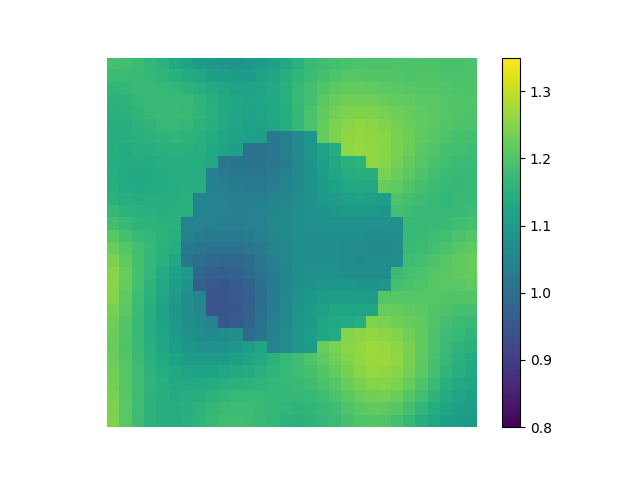}
\end{minipage}%
\begin{minipage}{0.2\textwidth}
  \centering
\includegraphics[width=1.\linewidth]{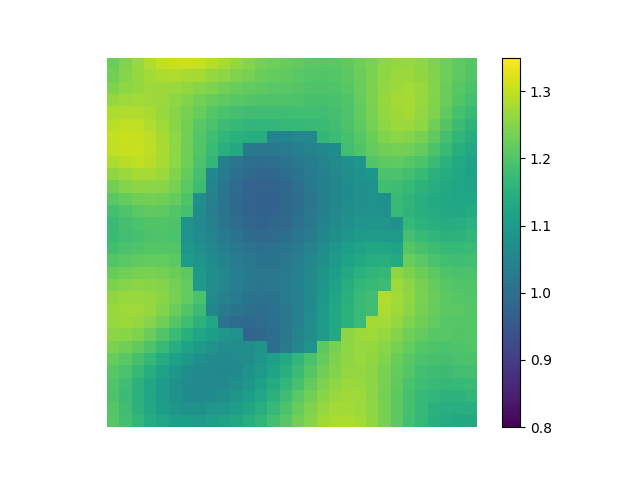}
\end{minipage}%
\begin{minipage}{0.2\textwidth}
  \centering
\includegraphics[width=1.\linewidth]{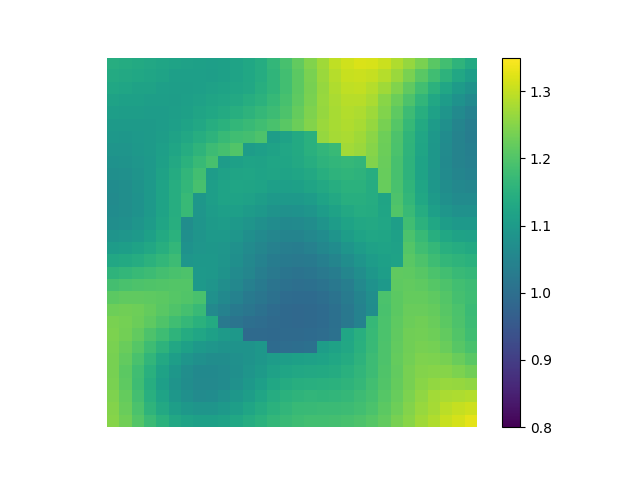}
\end{minipage} \\
\begin{minipage}{0.2\textwidth}
  \centering
\includegraphics[width=1.\linewidth]{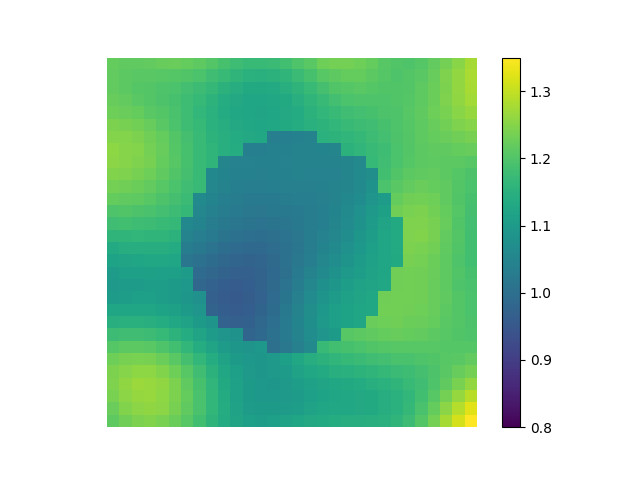}
\end{minipage}%
\begin{minipage}{0.2\textwidth}
  \centering
\includegraphics[width=1.\linewidth]{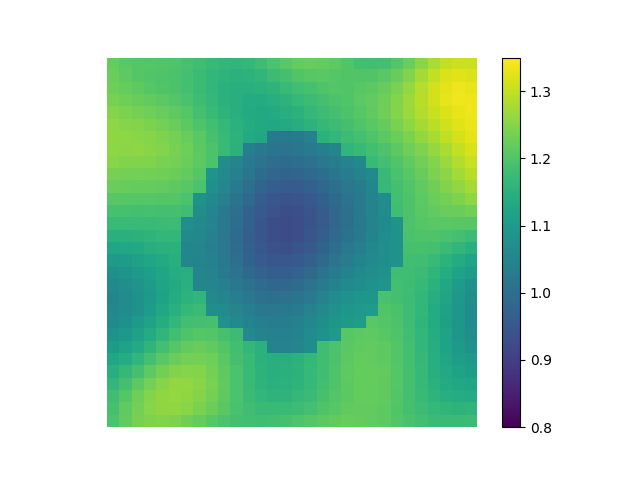}
\end{minipage}%
\begin{minipage}{0.2\textwidth}
  \centering
\includegraphics[width=1.\linewidth]{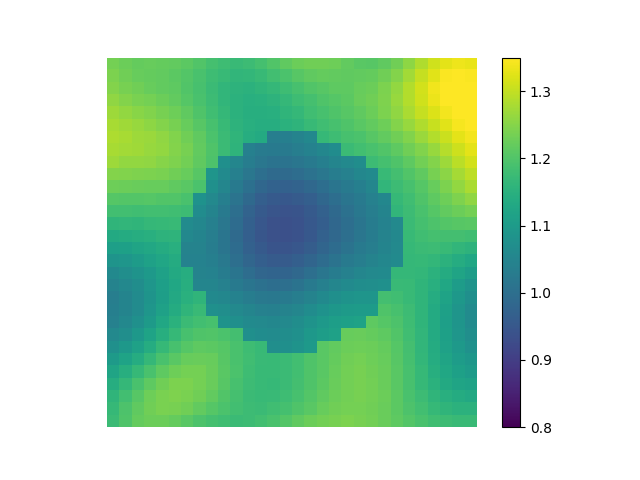}
\end{minipage}%
\begin{minipage}{0.2\textwidth}
  \centering
\includegraphics[width=1.\linewidth]{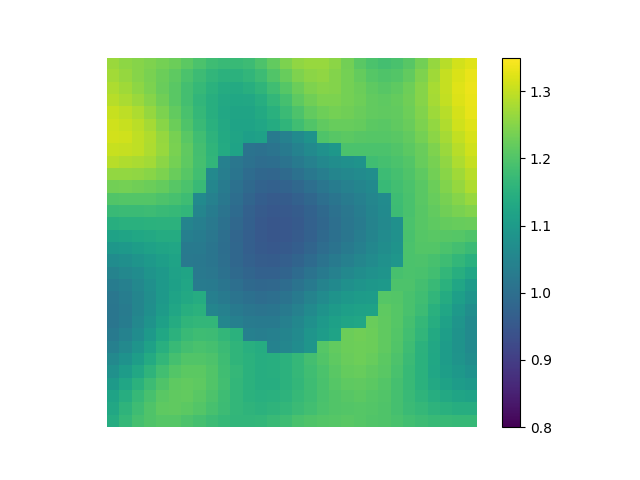}
\end{minipage}%
\begin{minipage}{0.2\textwidth}
  \centering
\includegraphics[width=1.\linewidth]{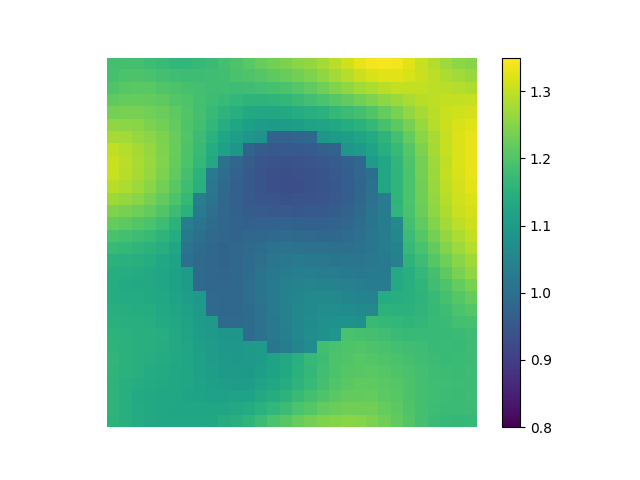}
\end{minipage} \\
\begin{minipage}{0.2\textwidth}
  \centering
\includegraphics[width=1.\linewidth]{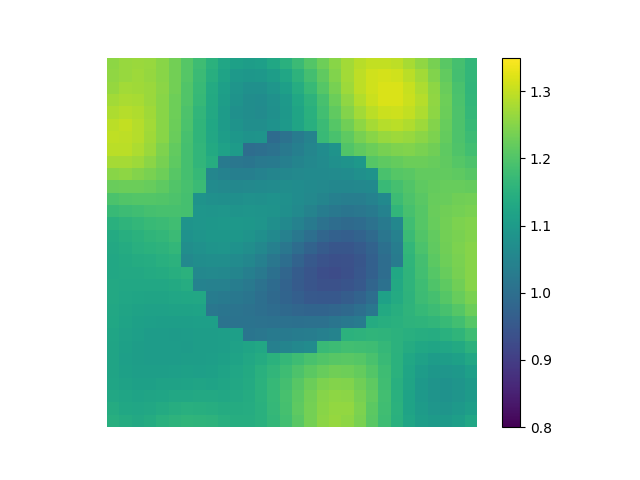}
\end{minipage}%
\begin{minipage}{0.2\textwidth}
  \centering
\includegraphics[width=1.\linewidth]{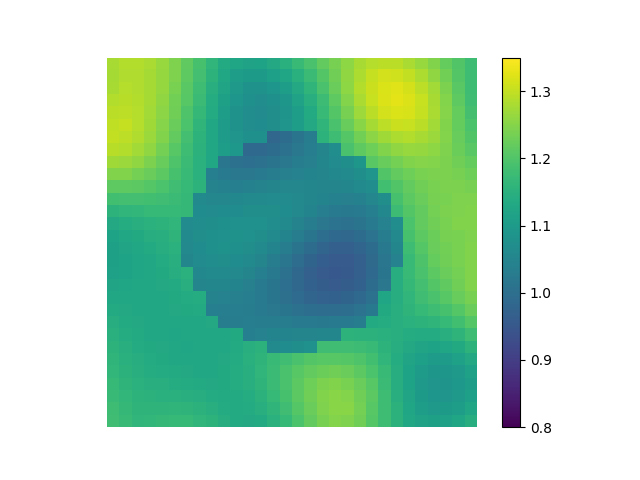}
\end{minipage}%
\begin{minipage}{0.2\textwidth}
  \centering
\includegraphics[width=1.\linewidth]{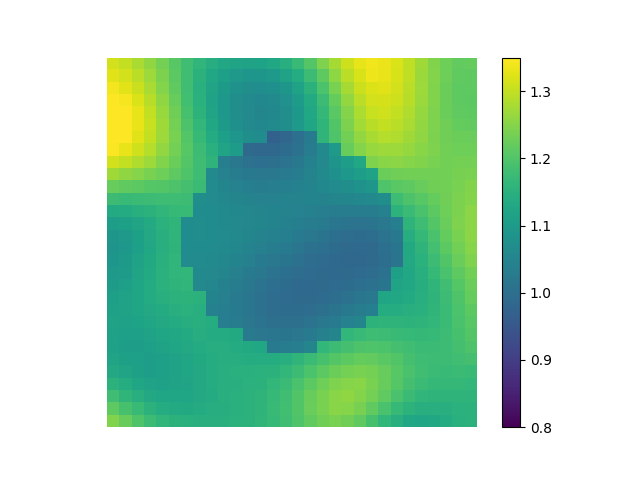}
\end{minipage}%
\begin{minipage}{0.2\textwidth}
  \centering
\includegraphics[width=1.\linewidth]{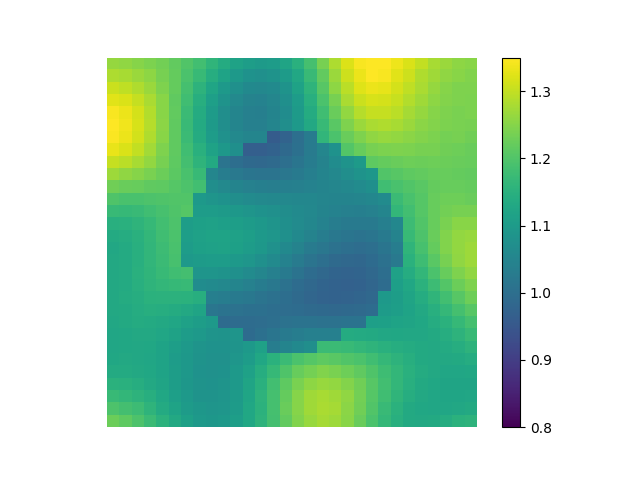}
\end{minipage}%
\begin{minipage}{0.2\textwidth}
  \centering
\includegraphics[width=1.\linewidth]{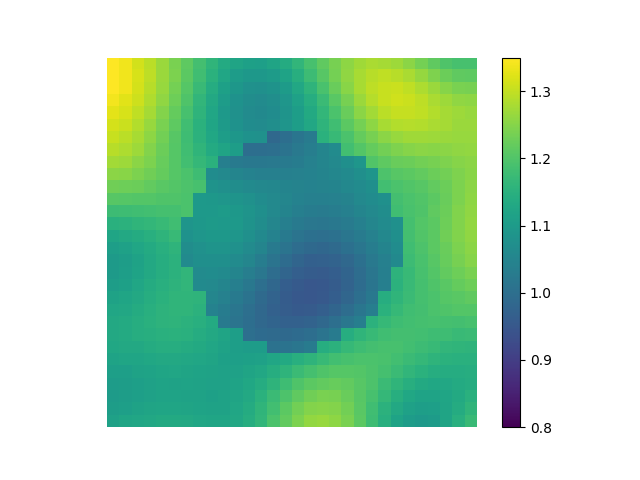}
\end{minipage} \\

\caption{Posterior samples for all methods. Top row corresponds to None, second row to EEA third row to OAE, and final row to the true posterior}
\label{fig:Darcy samples}
\end{figure}

The surrogate and true posterior are sampled by the same random-walk Metropolis--Hastings algorithm with $6000$ samples after a burn-in of $20,000$ and a thinning with the factor $t$. The posterior means and covariances $m(y_0)$, $C(y_0)$, $\hat m(y_0)$, $\widehat C(y_0)$ and the samples in Figure~\ref{fig:Darcy samples} are computed from these chains. 

\begin{figure}[h]
\centering
\begin{minipage}{0.23\textwidth}
  \centering
\includegraphics[width=1.1\linewidth]{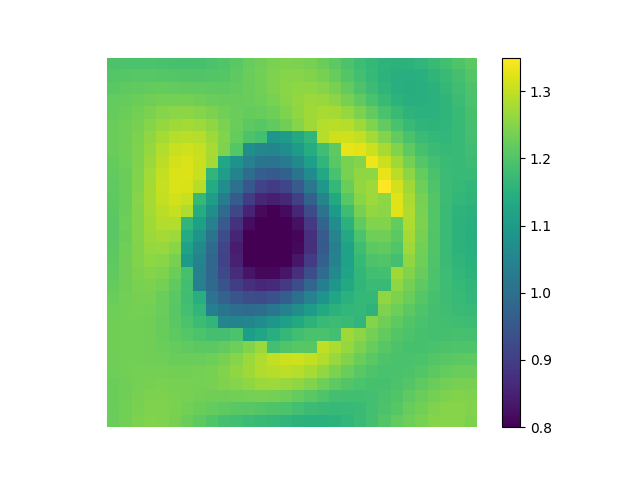}
\end{minipage}%
\begin{minipage}{0.23\textwidth}
  \centering
\includegraphics[width=1.1\linewidth]{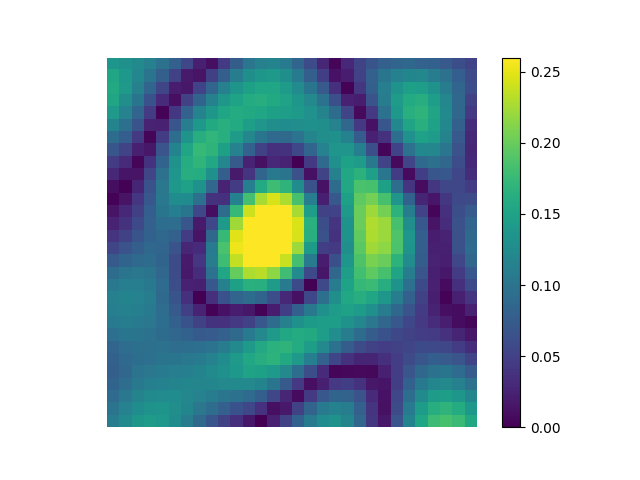}
\end{minipage}%
\begin{minipage}{0.23\textwidth}
  \centering
\includegraphics[width=1.1\linewidth]{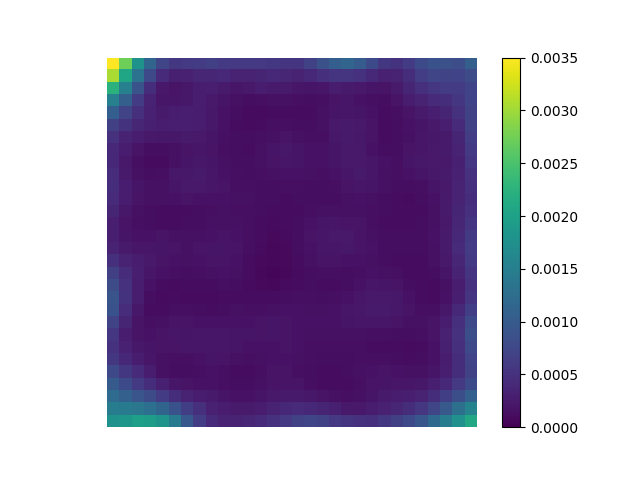}
\end{minipage}
\begin{minipage}{0.23\textwidth}
  \centering
\includegraphics[width=1.1\linewidth]{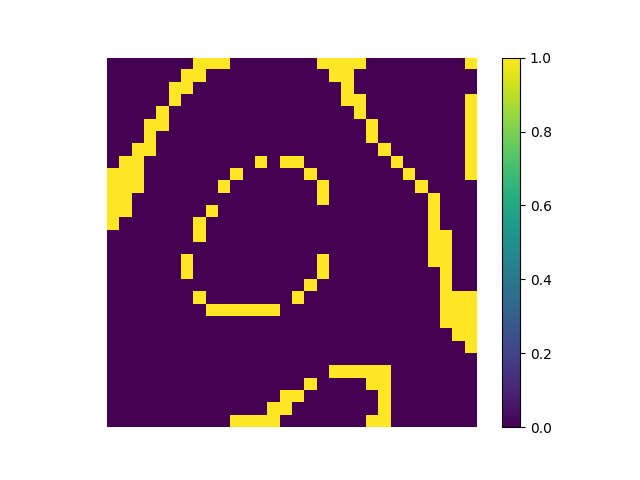}
\end{minipage} \\
\begin{minipage}{0.23\textwidth}
  \centering
\includegraphics[width=1.1\linewidth]{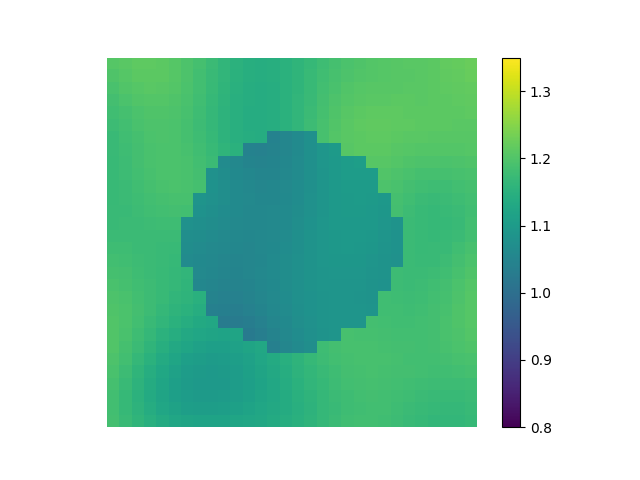}
\end{minipage}%
\begin{minipage}{0.23\textwidth}
  \centering
  \includegraphics[width=1.1\linewidth]{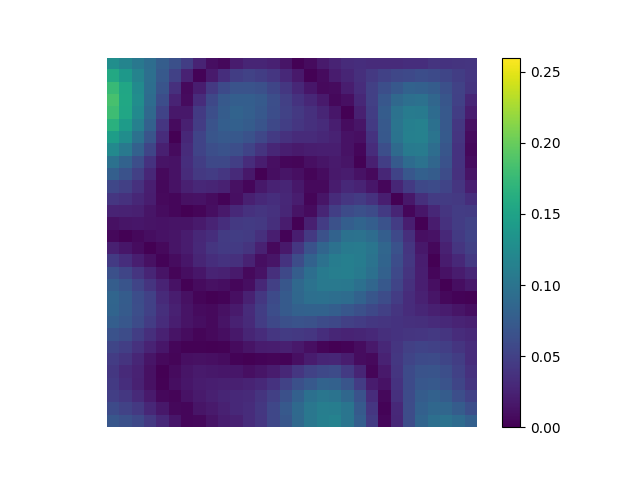}
\end{minipage}%
\begin{minipage}{0.23\textwidth}
  \centering
  \includegraphics[width=1.1\linewidth]{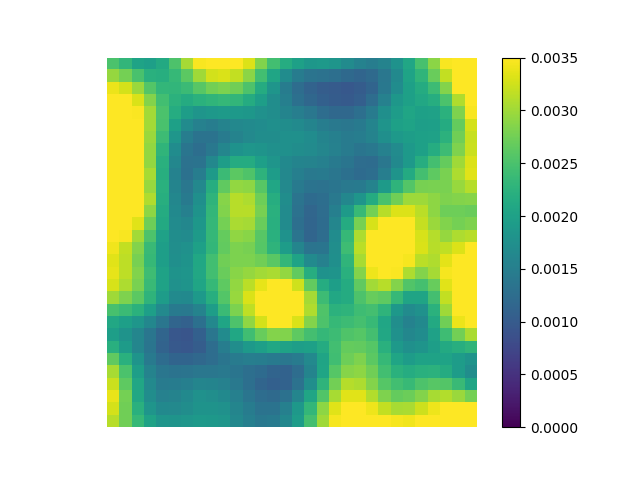}
\end{minipage}
\begin{minipage}{0.23\textwidth}
  \centering
  \includegraphics[width=1.1\linewidth]{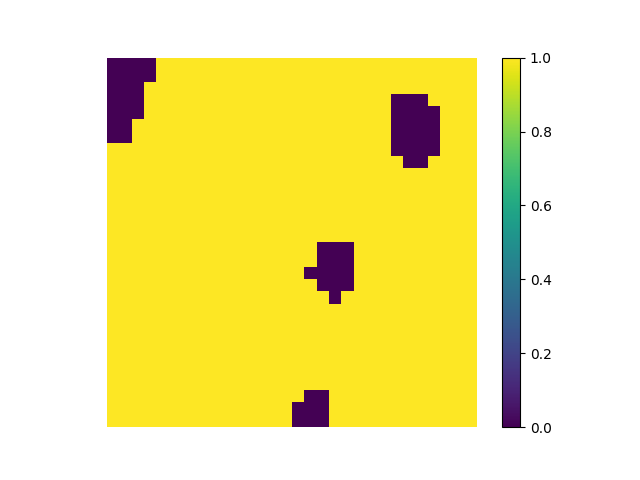}
\end{minipage} \\
\begin{minipage}{0.23\textwidth}
  \centering
  \includegraphics[width=1.1\linewidth]{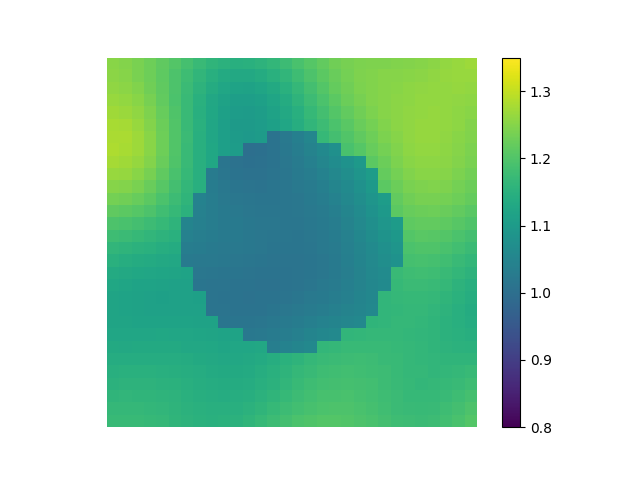}
\end{minipage}%
\begin{minipage}{0.23\textwidth}
  \centering
  \includegraphics[width=1.1\linewidth]{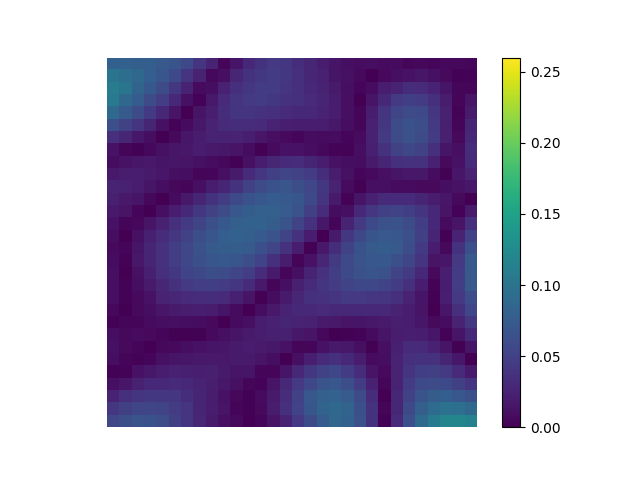}
\end{minipage}%
\begin{minipage}{0.23\textwidth}
  \centering
  \includegraphics[width=1.1\linewidth]{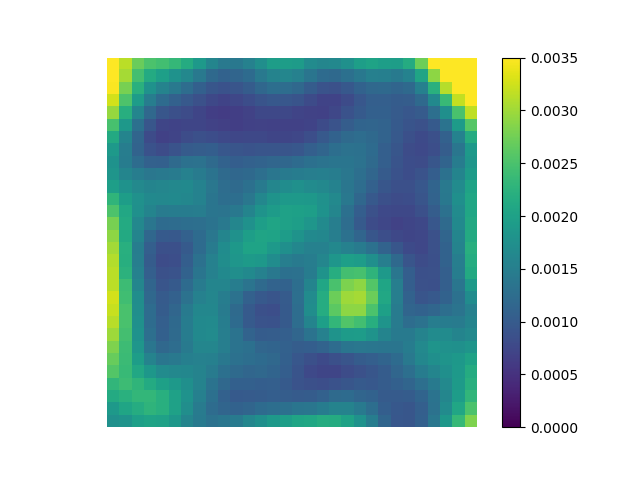}
\end{minipage}
\begin{minipage}{0.23\textwidth}
  \centering
  \includegraphics[width=1.1\linewidth]{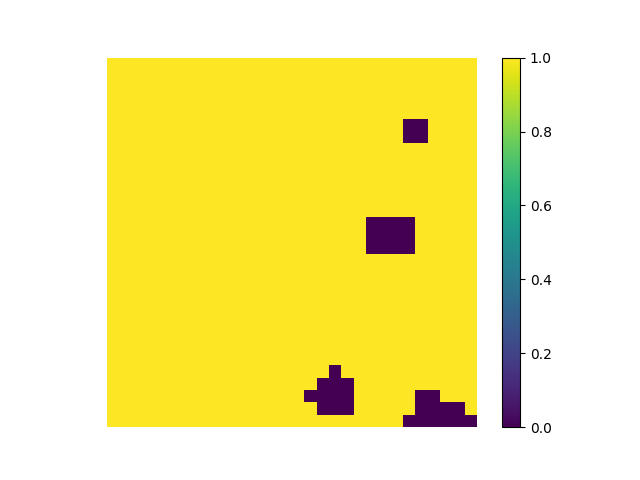}
\end{minipage}
\\
\begin{minipage}{0.23\textwidth}
  \centering
  \includegraphics[width=1.1\linewidth]{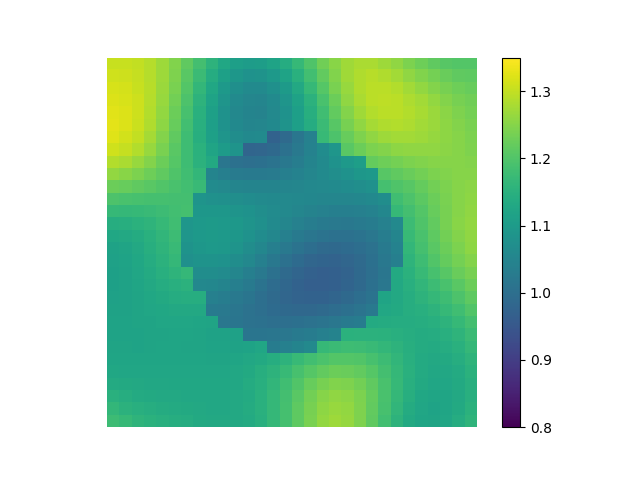}
\end{minipage}%
\begin{minipage}{0.23\textwidth}
  \centering
  \includegraphics[width=1.1\linewidth]{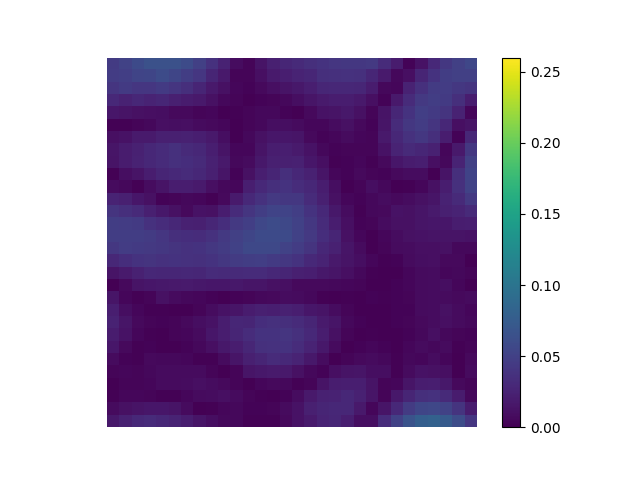}
\end{minipage}%
\begin{minipage}{0.23\textwidth}
  \centering
  \includegraphics[width=1.1\linewidth]{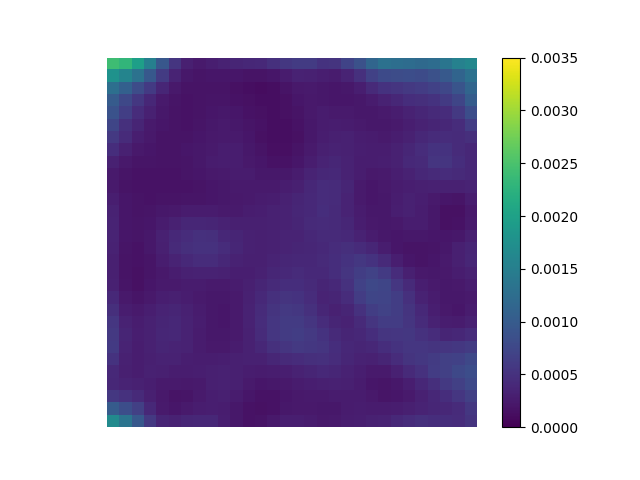}
\end{minipage}
\begin{minipage}{0.23\textwidth}
  \centering
  \includegraphics[width=1.1\linewidth]{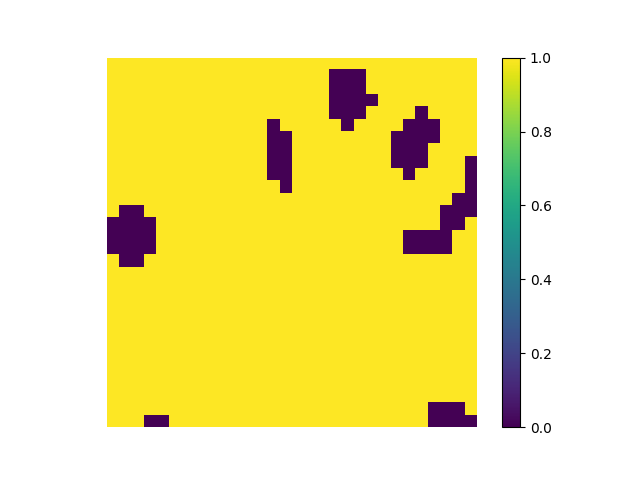}
\end{minipage}

\caption{For all methods, we report the posterior mean, difference posterior mean to the ground truth, diagonal of the posterior covariance, and the pixels contained in a confidence interval based on the posterior mean and variance. Rows  correspond to (top to bottom) None, EEA, OAE, true posterior}
\label{fig:Darcy}
\end{figure}

Table~\ref{tab: darcy law 2d table} collects the error metrics, Figure~\ref{fig:Darcy samples} shows posterior samples and Figure~\ref{fig:Darcy} compares the posterior means, their errors, the pointwise posterior standard deviation and the credible-band coverage with the ground truth.

The picture is consistent with the CT experiment. OAE halves the error of the posterior mean relative to EEA ($0.0015$ versus $0.0027$) and reduces the expected KL divergence, though here the gain over EEA is modest ($454$ versus $475$), while both improve by an order of magnitude on the uncorrected likelihood. Coverage of the $2\sigma$ band is $96\%$ for OAE and $93\%$ for EEA, compared with $91\%$ for the true posterior, so both corrections are, if anything, slightly conservative, whereas the uncorrected likelihood covers only $13\%$ of the pixels. As in the CT case, None has the smallest covariance error but a strongly biased mean.

\begin{table}[h]
    \centering
    \caption{Error metrics for the 2D Darcy-flow problem. All norms are averaged over the image pixels. In the quantities in the two last rows, $x_0$ corresponds to the unknown in Figure \ref{fig:2d darcy; prior sample} and $y_0$ to a simulated corresponding measurement.}
    \label{tab: darcy law 2d table}
    \begin{tabular}{c|c c c c}
         \hline
          & \textbf{None} & \textbf{EEA} & \textbf{OAE} & \textbf{Posterior} \\
         \hline
         $\norm{C(y_0) - \widehat{C}(y_0)}_F$ & \textbf{1.24e-8} &  5.12e-7 & 1.64e-7 & 0 \\ 
         $ \norm{m(y_0) - \hat{m}(y_0)}^2$   & 0.0132 &  0.0027 & \textbf{0.0015} & 0 \\ 
         $\E^{p} \, \KL{\mu(\cdot | y)}{\hat \mu(\cdot | y)}$ & 3387 & 475  &  \textbf{454} & 0 \\ 
         \hline
         $\norm{x_0 - \hat{m}(y_0)}^2$ & 0.0132 &  0.0029 &   0.0015 &0.0004 \\
         $\#  \left\{ |x_0-\hat{m}(y_0)| \le 2 \, \sqrt{\text{diag} \widehat{C}(y_0)} \right\}$ & 13 \% & 93 \% & 96 \% & 91 \% 
    \end{tabular}
\end{table}

\section{Conclusions}\label{sec: conclusions}
We introduced the optimal approximation error (OAE) framework for calibrating the noise parameters of the approximation error model by directly minimising the expected KL divergence between the true and surrogate posteriors. We showed that OAE subsumes the classical enhanced error approach as the variant that instead targets the expected divergence between likelihoods. We established smoothness, bi-convexity, and, after regularisation, coercivity and the Kurdyka–Łojasiewicz property, yielding a proximal alternating minimisation scheme with guaranteed convergence to a critical point. In the Gaussian linear case we obtained the objective function and the optimal mean shift in closed form. More generally, we introduced an empirical estimator and proved its consistency as the training sample size grows. The framework applies in a discretisation-invariant manner to unknowns in infinite-dimensional Banach spaces, and our experiments on a linear CT problem and a nonlinear Darcy-flow problem show that OAE improves the posterior mean and the expected KL divergence relative to the uncorrected and EEA-corrected likelihoods, with credible-band coverage comparable to that of EEA and of the true posterior.

Several directions remain open. The present framework parametrises the noise covariance by a full precision matrix $G \in \Sm$, which becomes impractical when $m$ is large; structured or low-rank parametrisations would reduce the parameter count while preserving the OAE optimality criterion. 
Extending OAE to non-Gaussian likelihood families would broaden its applicability to problems with heavy-tailed or count-valued observations. Finally, an online or sequential variant of OAE, in which the parameters $(\xi, G)$ are updated as new training 
pairs arrive, would be natural in settings where the high-fidelity model is evaluated adaptively.

\section*{Acknowledgments}
Part of this research was carried out while F.S. was at TU Wien, Austria.
Funding was received from the European Research Council (ERC) under the Horizon 2020 research and innovation program of the European Union (Grant agreement No. 101125225). This work was supported by the Research Council of Finland (353094, 348504 and 359183). Support by the Austrian Science Fund (FWF) grant I6667-N is acknowledged 

\subsection*{AI disclosure}
During the preparation of this manuscript, the authors used Claude Fable 5.1 to support proofreading and editing of the text and to improve its clarity. The authors reviewed, validated and critically assessed all AI-assisted outputs and take full responsibility for the content of the manuscript.


\printbibliography

\end{document}